\documentstyle[10pt]{article}

\title{A Classification Theorem for Nuclear Purely Infinite Simple
   C$^*$-Algebras
  \thanks{Research partially supported by NSF grant
     DMS 94-00904, and by the Fields
     Institute for Research in Mathematical Sciences. \protect\\
   AMS 1991 subject classification numbers: Primary 46L35;
   Secondary 19K99, 46L80. \protect\\
 }
  }
\author{N. Christopher Phillips \\ \\
  Department of Mathematics \\ University
    of Oregon \\ Eugene OR 97403-1222  \\  U.S.A.
   \\  and  \\
 Fields Institute for Research in Mathematical Sciences \\
 185 Columbia St.\  West \\
 Waterloo ON N2L 5Z5 \\ Canada \\  }
\date{9 December 1997}
\begin{document}
\maketitle

\newcommand{\beq}{\begin{equation}}
\newcommand{\eeq}{\end{equation}}
\newcommand{\bd}{\begin{displaymath}}
\newcommand{\ed}{\end{displaymath}}
\newcommand{\ben}{\begin{enumerate}}
\newcommand{\een}{\end{enumerate}}
\newcommand{\bde}{\begin{description}} 
\newcommand{\ede}{\end{description}}
\newcommand{\beqr}{\begin{eqnarray*}}
\newcommand{\eeqr}{\end{eqnarray*}}
\newcommand{\hs}[1]{\hspace{#1}}
\newcommand{\vs}[1]{\vspace{#1}}
\newcommand{\bc}{\begin{center}}
\newcommand{\ec}{\end{center}}
\newcommand{\bv}{\begin{verbatim}}
\newcommand{\ev}{\end{verbatim}}


\newcommand{\half}{\frac{1}{2}}
\newcommand{\p}{\partial}
\newcommand{\limi}[1]{\lim_{{#1} \to \infty}}

\newcommand{\af}{\alpha}
\newcommand{\bt}{\beta}
\newcommand{\gm}{\gamma}
\newcommand{\dt}{\delta}
\newcommand{\ep}{\varepsilon}
\newcommand{\zt}{\zeta}
\newcommand{\et}{\eta}
\newcommand{\ch}{\chi}
\newcommand{\io}{\iota}
\newcommand{\th}{\theta}
\newcommand{\ld}{\lambda}
\newcommand{\sm}{\sigma}
\newcommand{\kp}{\kappa}
\newcommand{\ph}{\varphi}
\newcommand{\ps}{\psi}
\newcommand{\rh}{\rho}
\newcommand{\om}{\omega}
\newcommand{\ta}{\tau}

\newcommand{\Gm}{\Gamma}
\newcommand{\Dt}{\Delta}
\newcommand{\Et}{\Eta}
\newcommand{\Th}{\Theta}
\newcommand{\Ld}{\Lambda}
\newcommand{\Sm}{\Sigma}
\newcommand{\Ph}{\Phi}
\newcommand{\Ps}{\Psi}
\newcommand{\Om}{\Omega}

\newcommand{\Q}{{\bf Q}}
\newcommand{\Z}{{\bf Z}}
\newcommand{\R}{{\bf R}}

\newcommand{\tta}{\theta}
\newcommand{\vph}{\varphi}

\renewcommand{\i}{\subset}

\newcommand{\s}[1]{\noindent {\bf EXERCISE~#1}{}}

\setcounter{section}{-1}  
\newtheorem{thm}{Theorem}[section]
\newtheorem{lem}[thm]{Lemma}
\newtheorem{prop}[thm]{Proposition}
\newtheorem{dfn}[thm]{Definition}
\newtheorem{cor}[thm]{Corollary}
\newtheorem{conj}[thm]{Conjecture}
\newtheorem{conv}[thm]{Convention}
\newtheorem{rmk}[thm]{Remark}
\newtheorem{ntt}[thm]{Notation}

\newcommand{\qed}{\hspace*{\fill}Q.E.D.}  
\newcommand{\cd}{\cdots}
\newcommand{\lr}{\longrightarrow}

\pagenumbering{arabic}


\newcommand{\OI}{{\cal O}_{\infty}}
\newcommand{\OIA}[1]{\OI \otimes {#1}}
\newcommand{\OA}[1]{{\cal O}_{#1}}
\newcommand{\Hom}{{\rm Hom}}
\newcommand{\Class}{\cal C}
\newcommand{\Boot}{\cal N}
\newcommand{\diag}{{\rm diag}}
\newcommand{\Ad}{{\rm Ad}}
\newcommand{\id}{{\rm id}}
\newcommand{\Cb}{C_{\rm b}}
\newcommand{\Kt}{K \otimes}
\newcommand{\KOI}[1]{\Kt \OIA{#1}}
\newcommand{\KOT}[1]{\Kt \OA{2} \otimes {#1}}
\newcommand{\OT}[1]{\OA{2} \otimes {#1}}
\newcommand{\ET}{\tilde{E}}
\newcommand{\EH}{\hat{E}}
\newcommand{\cel}{{\rm cel}}


\newcommand{\ca}{$C^*$-algebra}
\newcommand{\aab}{approximately absorbing}
\newcommand{\wab}{weakly approximately absorbing}
\newcommand{\wue}{asymptotic unitary equivalence}
\newcommand{\wyue}{asymptotically unitarily equivalent}
\newcommand{\aue}{approximate unitary equivalence}
\newcommand{\ayue}{approximately unitarily equivalent}
\newcommand{\mops}{mutually orthogonal projections}
\newcommand{\hm}{homomorphism}
\newcommand{\pisca}{purely infinite simple \ca}
\newcommand{\sep}{separable}
\newcommand{\kfalg}{separable nuclear unital \pisca}
\newcommand{\snus}{separable, nuclear, unital, and simple}
\newcommand{\kfive}
        {separable, nuclear, unital, purely infinite, and simple}
\newcommand{\tfn}{trivializing factorization}
\newcommand{\wtf}{asymptotically trivially factorizable}
\newcommand{\amm}{asymptotic morphism}
\newcommand{\wam}{local asymptotic morphism}
\newcommand{\ct}{continuous}
\newcommand{\mvn}{Murray-von Neumann equivalent}
\newcommand{\sctb}{second countable}
\newcommand{\lcp}{locally compact Hausdorff}
\newcommand{\cpt}{compact Hausdorff}
\newcommand{\sye}{asymptotically equal}
\newcommand{\dpj}{full projection}
\newcommand{\pj}{projection}
\newcommand{\wolog}{without loss of generality}
\newcommand{\tfae}{the following are equivalent}
\newcommand{\Tfae}{The following are equivalent}
\newcommand{\ifo}{if and only if}
\newcommand{\exl}{exponential length}
\newcommand{\andeqn}{\,\,\,\,\,\, {\rm and} \,\,\,\,\,\,}


\newcommand{\arrow}{\rightarrow}
\newcommand{\tdsum}{\widetilde{\oplus}}
\newcommand{\aueeps}[1]{\stackrel{#1}{\sim}}
\newcommand{\aeps}[1]{\stackrel{#1}{\approx}}
\newcommand{\dirlim}{\displaystyle \lim_{\longrightarrow}}
\newcommand{\QED}{\rule{1.5mm}{3mm}}
\newcommand{\smspace}{\vspace{0.6\baselineskip}}
\newcommand{\bigspace}{\vspace{\baselineskip}}


\begin{abstract}

Starting from Kirchberg's theorems announced at the operator
algebra conference in Gen\`eve in 1994, namely
$\OA{2} \otimes A \cong \OA{2}$
for separable unital nuclear simple $A$
and $\OIA{A} \cong A$ for
separable unital nuclear purely infinite simple $A,$ we prove that
$KK$-equivalence implies isomorphism for
nonunital separable nuclear \pisca s. It follows that if $A$ and $B$
are unital separable nuclear \pisca s
which satisfy the Universal Coefficient Theorem,
and if there is a graded isomorphism from $K_* (A)$ to $K_* (B)$
which preserves the $K_0$-class of the identity, then $A \cong B.$

Our main technical results are, we believe, of independent interest.
We say that two \amm s $t \mapsto \ph_t$ and $t \mapsto \ps_t$
from $A$ to $B$ are \wyue\  if there exists a continuous unitary path
$t \mapsto u_t$ in the unitization $B^+$ such that
$\| u_t \ph_t (a) u_t^* - \ps_t (a) \| \to 0$ for all $a$ in $A.$
We prove the following two results on deformations and unitary
equivalence. Let $A$
be \snus, and let $D$ be unital. Then any \amm\  from $A$ to
$\Kt \OIA{D}$ is \wyue\  to a \hm, and two homotopic \hm s from
$A$ to $\Kt \OIA{D}$ are necessarily \wyue.

We also give some nonclassification results for the nonnuclear case.
\end{abstract}

\vspace{\baselineskip}
\vspace{\baselineskip}

\section{Introduction}

\vspace{\baselineskip}

We prove that the isomorphism class of a \kfalg\  satisfying the
Rosenberg-Schochet Universal Coefficient Theorem is completely
determined by its $K$-theory. More precisely, let $A$ and $B$ be
\kfalg s which satisfy the Universal Coefficient Theorem, and suppose
that there is a graded isomorphism $\af : K_* (A) \to K_* (B)$ such
that $\af ([1_A]) = [1_B]$ in $K_0 (B).$ Then there is an isomorphism
$\ph : A \to B$ such that $\ph_* = \af.$ This theorem follows from
a result asserting that whenever $A$ and $B$ are \kfalg s (not
necessarily satisfying the Universal Coefficient Theorem) which are
$KK$-equivalent via a class in $KK$-theory which respects the classes of
the identities, then there is an isomorphism from $A$ to $B$ whose class
in $KK$-theory is the given one.

As intermediate results, we prove some striking facts about \hm s and
\amm s from a separable nuclear unital simple \ca\  to a the tensor
product of a unital \ca\  and the Cuntz algebra $\OI.$ If $A$ and $D$
are any two \ca s, we say that two \hm s $\ph, \, \ps : A \to D$ are
{\em \wyue} if there is a \ct\  unitary path $t \mapsto u_t$ in
$D$ such that
$\lim_{t \to \infty} u_t \ph (a) u_t* = \ps (a)$ for all $a \in A.$
(Here $\tilde{D} = D$ if $D$ is unital, and $\tilde{D}$ is the
unitization $D^+$ if $D$ is not unital.) Note that \wue\  is a
slightly strengthened form of \aue, and is an approximate form of
unitary equivalence. Our results show that if $A$ is \snus,  and
$D$ is separable and unital, then $KK^0 (A, D)$ can be computed as the
set of \wue\  classes of full \hm s from $A$ to $\KOI{D},$ with direct
sum as the operation. Note that we use something close to unitary
equivalence, and that there is no need to use \amm s, no need
to take suspensions, and (essentially because $\OI$ is purely
infinite) no need to form formal differences of classes. We can
furthermore
replace $A$ by $\KOI{A},$ in which case the Kasparov product reduces
exactly to composition of \hm s. These results can be thought of as
a form of unsuspended $E$-theory. (Compare with \cite{DL}, but note
that we don't even need to use \amm s.)
There are also perturbation results: any \amm\  is in fact
\wyue\  (with a suitable definition) to a \hm.

We also present what is now known about how badly the classification
fails in the nonnuclear case.
There are separable \pisca s $A$ with $\OIA{A} \not\cong A$
(Dykema--R\o rdam), there are infinitely many nonisomorphic
separable exact \pisca s $A$ with $\OIA{A} \cong A$ and
$K_* (A) = 0$ (easily obtained from results of Haagerup and
Cowling--Haagerup), and for given K-theory there are uncountably
many nonisomorphic
separable nonexact \pisca s with that K-theory.

Classification of \ca s started with Elliott's classification
\cite{Ell1} of AF algebras up to isomorphism by their $K$-theory.
It received new impetus with his successful classification of
certain \ca s of real rank zero with nontrivial $K_1$-groups. We refer
to \cite{Ell3} for a recent comprehensive list of work in this area.
The initial step toward classification in the infinite case was taken
in \cite{BKRS}, and was quickly followed by a number of papers
\cite{Rr1}, \cite{Rr2}, \cite{Ln5}, \cite{LP1}, \cite{ER}, \cite{Rr3},
\cite{LP2}, \cite{BEEK}, \cite{Rr4}, \cite{KK}, \cite{LP3}.
In July 1994, Kirchberg announced \cite{Kr1}
a breakthrough: proofs that if $A$ is a \kfalg, then
$\OA{2} \otimes A \cong \OA{2}$ and $\OIA{A} \cong A.$
(The proofs, closely following Kirchberg's original methods, are
in \cite{KP}.)
This quickly led to two more papers \cite{Ph2}, \cite{Rr5}.
Here, we use Kirchberg's results to nearly solve the classification
problem for \kfalg s; the only difficulty that remains is the Universal
Coefficient Theorem. The method is a great generalization of that of
\cite{Ph2}, in which we replace \hm s by \amm s and
\aue\  by \wue. We also need a form of unsuspended $E$-theory, as
alluded to above. The most crucial step is done in Section 2,
where we show that, in a particular context, homotopy implies \wue.
We suggest reading \cite{Ph2} to understand the
basic structure of Section 2.

Kirchberg has in \cite{Kr2} independently derived the same
classification theorem we have. His methods are somewhat different,
and mostly independent of the proofs in \cite{KP}.
He proves that homotopy implies a form of unitary equivalence in a
different context, and does so by eventually reducing the problem to
a theorem of this type in Kasparov's paper \cite{Ksp}. By contrast,
the main machinery in our proof is simply the repeated use of
Kirchberg's earlier results as described above.

This paper is organized as follows. In Section 1, we present some
important facts about \amm s, and introduce \wue. In Section 2,
we prove our main technical results: under suitable conditions,
homotopic \amm s are \wyue\  and  \amm s are \wyue\  to \hm s. These
results are given at the end of the section. In Section 3, we prove the
basic form (still using \amm s) of our version of unsuspended
$E$-theory.  Finally, Section 4 contains the classification theorem
and some corollaries, as well as the nicest forms of the intermediate
results discussed above. It also contains the nonclassification results.

Most of
this work was done during a visit to the Fields Institute for Research
in Mathematical Sciences during Fall 1994, and I would like to thank the
Institute for its support and for the stimulating research environment
it provided. I would also like to thank a number of people for
useful discussions, either in person or by electronic mail, including
Marius D\v{a}d\v{a}rlat, George Elliott, Uffe Haagerup,
Eberhard Kirchberg,
Alex Kumjian, Huaxin Lin,
Mikael R\o rdam, Jonathan Samuel, Claude Schochet,
and Shuang Zhang. These discussions have led me to
considerable simplification of the arguments and improvement of the
terminology.

Throughout this paper, $U (D)$ denotes the unitary group of a unital
\ca\  $D,$ and $U_0 (D)$ denotes the connected component of $U (D)$
containing 1. We will use repeatedly and without comment Cuntz's
result that $K_1 (D) = U (D) / U_0 (D)$ for a unital \pisca\  $D,$
as well as his corresponding result that $K_0 (D)$ is the set of
Murray-von Neumann equivalence classes of nonzero projections
\cite{Cu2}. We similarly use Kasparov's $KK$-theory \cite{Ksp}, and we
recall here (and do not mention again) that every separable nonunital
\pisca\  has the form $\Kt D$ for a unital \pisca\  $D$ \cite{Zh0}.

\vspace{\baselineskip}
\section{Asymptotic morphisms and \wue}
\vspace{\baselineskip}

The basic objects we work with in this paper are \amm s.
In the first subsection, we state for convenient reference some of
the facts we need about \amm s, and establish notation concerning them.
In the second subsection, we define and discuss full \amm s; fullness
is used as a nontriviality condition later in the paper.
In the third subsection, we introduce \wue\  of \amm s. This relation
is the appropriate version of unitary equivalence in the context of
\amm s, and will play a fundamental role in Sections 2 and 3.

\vspace{0.6\baselineskip}
\subsection{Asymptotic morphisms and \wue}
\vspace{0.6\baselineskip}

Asymptotic morphisms were introduced by Connes and Higson \cite{CH}
for the purpose of defining $E$-theory, a simple construction of
$KK$-theory (at least if the first variable is nuclear).
In this subsection, we recall the definition and some of the basic
results on \amm s, partly to establish our notation and partly for
ease of reference. We also prove a few facts that are well known but
seem not to have been published. We refer to \cite{CH}, and the
much more detailed paper \cite{Sm},
for the details of the rest of the development of $E$-theory.

If $X$ is a \cpt\  Hausdorff space, then $C (X, D)$ denotes the
\ca\  of all \ct\  functions from $X$ to $D,$ while if $X$ is
\lcp\  Hausdorff, then $C_0 (X, D)$ denotes the \ca\  of all
\ct\  functions from $X$ to $D$ which vanish at infinity, and
$\Cb (X, D)$ denote the \ca\  of all bounded \ct\  functions from
$X$ to $D.$

We begin by recalling the definition of an \amm.

\vspace{0.6\baselineskip}

{\bf  1.1.1 Definition.}
Let $A$ and $D$ be \ca s, with $A$ separable.
An {\em \amm} $\ph : A \to D$
is a family $t \to \ph_t$ of functions from $A$ to $D,$
defined for $t \in [0, \infty),$ satisfying the following conditions:

(1) For every $a \in A,$ the function $t \mapsto \ph_t (a)$ is
\ct\  from $[0, \infty)$ to $D.$

(2) For every $a, b \in A$ and $\af, \bt \in {\bf C},$ the limits
\[
\lim_{t \to \infty} (\ph_t (\af a + \bt b) - \af \ph_t (a)
                                                - \bt \ph_t (b)),
\]
\[
\lim_{t \to \infty} (\ph_t (ab) - \ph_t (a) \ph_t (b)), \andeqn
\lim_{t \to \infty} (\ph_t (a^*) - \ph_t (a)^*)
\]
are all zero.

\vspace{0.6\baselineskip}

{\bf  1.1.2 Definition.} (\cite{CH}) Let $\ph$ and $\ps$ be
\amm s  from $A$ to $D.$

(1) We say that $\ph$ and $\ps$ are {\em asymptotically equal} (called
``equivalent'' in \cite{CH}) if for all $a \in A$, we have
$\lim_{t \to \infty} (\ph_t (a) - \ps_t (a)) = 0.$

(2) We say that $\ph$ and $\ps$ are {\em homotopic} if
there is an \amm\  $\rh : A \to C ([0, 1], D)$ whose
restrictions to $\{0\}$ and $\{1\}$ are
$\ph$ and $\ps$ respectively. In this case, we refer to
$\af \mapsto \rh^{(\af)} = {\rm ev}_{\af} \circ \rh$ (where
${\rm ev}_{\af} : C ([0, 1], D) \to D$ is evaluation at $\af$)
as a homotopy from $\ph$ to $\ps,$ or as a continuous path of
\amm s from $\ph$ to $\ps.$

The set of homotopy classes of \amm s from $A$ to $D$ is denoted
$[[A, D]],$ and the homotopy class of an \amm\ $\ph$ is denoted
$[[\ph]].$

\vspace{0.6\baselineskip}

It is easy to check that asymptotic equality implies homotopy
(\cite{Sm}, Remark 1.11).

\vspace{0.6\baselineskip}

{\bf  1.1.3 Definition.}
Let $\ph, \ps : A \to \Kt D$ be \amm s. The direct sum
$\ph \oplus \ps,$ well defined up to unitary equivalence (via
unitaries in $M(\Kt D)$), is defined as follows.
Choose any isomorphism $\dt : M_2 (K) \to K,$ let
$\overline{\dt} : M_2 (\Kt D) \to \Kt D$ be the induced map, and define
\[
(\ph \oplus \ps)_t (a) = \overline{\dt} \left( \left(
     \begin{array}{cc} \ph_t (a) & 0 \\ 0 & \ps_t (a) \end{array}
  \right) \right).
\]
Note that any two choices for $\dt$ are unitarily equivalent (and
hence homotopic).

\vspace{0.6\baselineskip}

The individual maps $\ph_t$ of an \amm\  are not assumed bounded
or even linear.

\vspace{0.6\baselineskip}

{\bf  1.1.4 Definition.}
Let $\ph : A \to D$ be an \amm .

(1) We say that $\ph$ is {\em completely positive contractive}
if each $\ph_t$ is a linear completely positive contraction.

(2) We say that $\ph$ is {\em bounded} if each $\ph_t$ is linear
and $\sup_{t} \| \ph_t \| $ is finite.

(3) We say that $\ph$ is {\em selfadjoint} if
$\ph_t (a^*) = \ph_t (a)^*$ for all $t$ and $a.$

Unless otherwise specified, homotopies of \amm s  from $A$ to $D$
satisfying one or more of these conditions will be assumed to
satisfy the same conditions as \amm s  from $A$ to $C ([0, 1], D).$

\vspace{0.6\baselineskip}

Note that if $\ph$ is bounded, then the formula
$\ps_t (a) = \half (\ph_t (a) + \ph_t (a^*)^*)$ defines a
selfadjoint bounded asymptotic morphism which is \sye\  to $\ph.$
We omit the easy verification that $\ps$ is in fact an \amm .

\vspace{0.6\baselineskip}

{\bf    1.1.5 Lemma.} (\cite{Sm}, Lemma 1.6.)
Let $A$ and $D$ be \ca s, with $A$ separable and nuclear.
Then every
\amm\  from $A$ to $D$ is asymptotically equal to
a completely positive contractive \amm . Moreover, the obvious
map defines a bijection between the sets of homotopy classes
of completely positive contractive \amm s
and arbitrary \amm s. (Homotopy classes are as in the
convention  in Definition 1.1.4.)

\vspace{0.6\baselineskip}

{\bf  1.1.6 Lemma.}
Let $\ph : A \to D$ be an \amm .  Define $\ph^+ : A^+ \to D^+$ by
$\ph_{t} (a + \ld \cdot 1) = \ph_{t} (a) + \ld \cdot 1$ for $a \in A$
and $\ld \in {\bf C}.$ Then $\ph^+$ is an \amm\  from $A^+$ to $D^+,$
and is completely positive contractive, bounded, or selfadjoint whenever
$\ph$ is.

\vspace{0.6\baselineskip}

The proof of this is straightforward, and is omitted.

The following result is certainly known, but we know of no reference.

\vspace{0.6\baselineskip}

{\bf   1.1.7 Proposition.} Let $A$ be a \ca\  which is given by
exactly stable (in the sense of Loring \cite{Lr1}) generators
and relations $(G, R),$ with  both $G$ and $R$ finite. Let
$D$ be a \ca . Then any \amm\   from $A$ to $D$ is asymptotically equal
to a continuous family of \hm s from $A$ to $D$ (parametrized by
$[0, \infty)$). Moreover, if  $\ph^{(0)}$ and $\ph^{(1)}$
are two homotopic \amm s from $A$ to $D,$ such that each $\ph^{(0)}_t$
and each $\ph^{(1)}_t$ is a \hm, then there
is a homotopy $\af \mapsto \ph^{(\af)}$ which is \sye\ to the given
homotopy and such that each $\ph^{(\af)}_{t}$ is a \hm.

\vspace{0.6\baselineskip}

Note that it follows from Theorem 2.6 of \cite{Lr2}
that exact stability of $(G, R)$ depends only on $A$, not on the
specific choices of $G$ and $R.$

\vspace{0.6\baselineskip}

{\em Proof of Proposition   1.1.7:}
Theorem 2.6 of \cite{Lr2} implies that
the algebra $A$ is semiprojective in the sense of
Blackadar \cite{Bl1}. (Also see Definition 2.3 of \cite{Lr2}.)
We will use semiprojectivity instead of exact stability.

We prove the first statement.
Let $\ph: A \to D$ be a \amm . Then $\ph$ defines in a standard way
(see Section 1.2 of \cite{Sm}) a
\hm\  $\ps : A \to \Cb ([0, \infty), D)/C_0 ([0, \infty), D).$ Let
\[
I_n (D) = \{ f \in \Cb ([0, \infty), D) :
      f(t) = 0 \,\, {\rm for} \,\, t \geq n\}.
\]
Then
$C_0 ([0, \infty), D) = \overline{\bigcup_{n = 1}^{\infty} I_n (D)}.$
Semiprojectivity of $A$ provides an $n$ and a
\hm\  $\sm : A \to  \Cb ([0, \infty), D)/I_n (D)$ such that
the composite of $\sm$ and the quotient map
\[
\Cb ([0, \infty), D)/I_n (D) \to
           \Cb ([0, \infty), D)/C_0 ([0, \infty), D)
\]
is $\ps.$ Now $\sm$ can be viewed as a continuous
family of \hm s $\sm_{t}$ from $A$ to $D,$ parametrized by
$[n, \infty).$ Define $\sm_{t} = \sm_{n}$ for
$0 \leq t \leq n.$ This gives the required continuous family
of \hm s.

The proof of the statement about homotopies is essentially the
same. We use $\Cb ([0, 1] \times [0, \infty), D)$ in place
of $\Cb ([0, \infty), D),$
\[
J = \{ f \in C_0 ([0, 1] \times [0, \infty), D):
  f(\af, t) = 0 \,\, {\rm for} \,\, \af = 0, 1 \}
\]
in place of $C_0 ([0, \infty), D),$ and $J \cap I_n ([0, 1], D)$
in place of $I_n (D).$
We obtain $\ph^{(\af)}_{t}$ for all $t$ greater than or equal to
some $t_0,$ and for all $t$ when $\af = 0$ or $1.$ We then
extend over $(0, 1) \times [0, t_0)$ via a \ct\  retraction
\[
[0, 1] \times [0, \infty) \to  ([0, 1] \times [t_0, \infty)) \cup
   (\{0\} \times [0, \infty)) \cup (\{1\} \times [0, \infty)).
\]
\QED

\vspace{0.6\baselineskip}

We refer to
\cite{CH} (and to \cite{Sm} for more detailed proofs) for the
definition of $E (A, B)$ as the abelian group of homotopy classes
of \amm s from $\Kt SA$ to $\Kt SB,$ for the construction of the
composition of \amm s (well defined up to homotopy), and for the
construction of the natural map $KK^0 (A, B) \to E(A, B)$ and the
fact that it is an isomorphism if $A$ is nuclear. We do state here for
reference the existence of the tensor product of \amm s. For the proof,
see Section 2.2 of \cite{Sm}.

\vspace{0.6\baselineskip}

{\bf   1.1.8 Proposition.} (\cite{CH})
Let $A_1,$ $A_2,$ $B_1,$ and $B_2$ be
separable \ca s, and let $\ph^{(i)} : A_i \to B_i$ be \amm s. Then
there exists an \amm\  $\ps : A_1 \otimes A_2 \to B_1 \otimes B_2$
(maximal tensor products) such that
$\ps_t (a_1 \otimes a_2) - \ph_t^{(1)} (a_1) \otimes \ph_t^{(2)} (a_2)
                                        \to 0$
as $t \to \infty,$ for all $a_1 \in A_1$ and $a_2 \in A_2.$ Moreover,
$\ps$ is unique up to asymptotic equality.

\vspace{0.6\baselineskip}
\subsection{Full \amm s}
\vspace{0.6\baselineskip}

In this subsection, we define full   \amm s.
Fullness will be used as a nontriviality condition on \amm s in
Section 3.
It will also be convenient (although not, strictly speaking,
necessary) in Section 2.

We make our definitions in terms of projections, because the behavior
of \amm s on projections can be reasonably well controlled. We do
not want to let the \amm\  $\ph : C_0 ({\bf R}) \to C_0 ({\bf R}),$
defined by $\ph_t (f) = t f,$ be considered to be full, since it
is \sye\  to the zero \amm, but in the absence of \pj s it is not
so clear how to rule it out. Fortunately, in the present paper this
issue does not arise.

We start with a useful definition and some observations related to the
evaluation of \amm s on \pj s.

\vspace{0.6\baselineskip}

{\bf  1.2.1 Definition.}
Let $A$ and $D$ be \ca s, with $A$ separable.
Let $p \in A$ be a projection, and let
$\ph : A \to D$ be an \amm . A {\em tail \pj}\  for
$\ph (p)$ is a continuous function $t \mapsto q_{t}$
from $[0, \infty)$ to the projections in $D$ which, thought of
as an \amm\  $\ps : {\bf C} \to D$ via
$\ps_{t} (\ld) = \ld q_{t},$ is asymptotically equal to the
\amm\  $\ps'_{t} (\ld) = \ld \ph_{t} (p).$

\vspace{0.6\baselineskip}

{\bf  1.2.2 Remark.}
(1) Tail projections always exist: Choose a suitable $t_0,$
apply functional calculus
to $\half (\ph_{t} (p) + \ph_{t} (p)^*)$ for $t \geq t_0,$
and take the value at $t$ for $t \leq t_0$ to be
the value at $t_0.$ (Or use Proposition   1.1.7.)

(2) If $\ph$ is an \amm\  from $A$ to $D,$ then a tail \pj\  for
$\ph (p),$ regarded as an \amm\  from ${\bf C}$ to $D,$
is a representative of the product homotopy class of
$\ph$ and the \amm\  from ${\bf C}$ to $A$ given by $p.$

(3) A {\em homotopy} of tail \pj s is defined in the obvious way: it is
a continuous family of projections
$( \af, t) \to q_{t}^{(\af)}$ with given values at
$\af = 0$ and $\af = 1$.

(4) If $\ph$ is an \amm, then
it makes sense to say that a tail \pj\  is (or is not) full
(that is, generates a full hereditary subalgebra),
since  fullness depends only on the homotopy class of a \pj.

\vspace{0.6\baselineskip}

{\bf  1.2.3 Lemma.}
Let $A$ and $D$ be \ca s, with $A$ separable.
Let $\ph : A \to D$ be
an \amm, and let $p_1$ and $p_2$ be \pj s in $A$.
If $p_1$ is \mvn\  to a sub\pj\  of $p_2,$ then
a tail projection for $\ph (p_1)$ is  \mvn\  to a sub\pj\  of
a tail projection for $\ph (p_2)$.

\vspace{0.6\baselineskip}

{\em Proof:}
Let $t \mapsto q_t^{(1)}$ and $t \mapsto q_t^{(2)}$  be
tail \pj s for $\ph (p_1)$ and $\ph (p_2)$ respectively. Let $v$ be a
partial isometry with $v^* v = p_1$ and $v v^* \leq p_2.$
Using asymptotic multiplicativity and the definition of a
tail projection, we have
\[
\lim_{t \to \infty} (\ph_t (v)^* \ph_t (v) - q_t^{(1)}) = 0
       \andeqn
\lim_{t \to \infty}
        (q_t^{(2)} \ph_t (v) \ph_t (v)^* q_t^{(2)} - q_t^{(2)}) = 0.
\]
It follows that for $t$ sufficiently large,
$q_t^{(1)}$ is \mvn\  to a sub\pj\  of $q_t^{(2)},$
with the Murray-von Neumann equivalence depending \ct ly on $t.$
It is easy to extend it from an interval $[t_0, \infty)$ to
$[0, \infty)$. \QED

\vspace{0.6\baselineskip}

{\bf  1.2.4 Lemma.} Let $A$ and $D$ be as in Definition  1.2.1,
let $\af \mapsto \ph^{(\af)}$  be a homotopy of \amm s from $A$ to
$D,$ and let $p_0, \, p_1 \in A$ be homotopic projections.
Let $q^{(0)}$ and $q^{(1)}$ be tail \pj s for $\ph^{(0)} (p_0)$ and
$\ph^{(1)} (p_1)$ respectively. Then $q^{(0)}$ is homotopic to
$q^{(1)}$ in the sense of Remark  1.2.1 (3).

\vspace{0.6\baselineskip}

{\em Proof:}
This can be proved directly, but also follows by combining Remark
1.2.2 (2), Proposition  1.1.7, and the fact that products of homotopy
classes of \amm s are well defined.
\QED

\vspace{0.6\baselineskip}

{\bf  1.2.5 Definition.}
Let $A$ be a separable \ca\  which contains a full \pj, and let
$D$ be any \ca. Then an \amm\  $\ph : A \to D$ is {\em full} if there
is a full \pj\  $p \in A$ such that some (equivalently, any) tail
\pj\  for $\ph (p)$ is full in $D.$

\vspace{0.6\baselineskip}

This definition rejects, not only the identity map of
$C_0 ({\bf R}),$ but also the identity map of $C_0 ({\bf Z}).$
(The algebra $C_0 ({\bf Z})$ has no full \pj s.)
However, it will do for our purposes.

Note that, by Lemma  1.2.3, if a tail \pj\  for $\ph (p)$ is full,
then so is a tail \pj\  for $\ph (q)$ whenever $p$ is \mvn\  to a
sub\pj\  of $q.$

We now list the relevant properties of full \amm s. We omit the proofs;
they are mostly either immediate or variations on the proof of
Lemma  1.2.3.

\vspace{0.6\baselineskip}

{\bf  1.2.6 Lemma.}
(1) Fullness of an \amm\  depends only on its homotopy class.

(2) If $\ph, \ps : A \to D$ are \amm s, and if $\ph$ is
full, then so is the \amm\  $\ph \oplus \ps : A \to M_2 (D).$

(3) Let $B$ be separable, and have a full \pj, and further
assume that given two full \pj s in $B,$ each is \mvn\  to a
subprojection of the other. Then any \amm\  representing the product of
full \amm s from $A$ to $B$ and from $B$ to $D$ is again full.

\vspace{0.6\baselineskip}

The extra assumption in part (3) is annoying, but we don't see an easy
way to avoid it. This suggests that we don't quite have the right
definition. However, in this paper $B$ will almost always have the form
$\KOI{D}$ with $D$ unital. Lemma  2.1.8 (1) below will ensure that
the assumption holds in this case.

\vspace{0.6\baselineskip}
\subsection{Asymptotic unitary equivalence}
\vspace{0.6\baselineskip}

Approximately unitarily equivalent \hm s have the same class in
R\o rdam's $KL$-theory (Proposition 5.4 of \cite{Rr4}), but need not
have the same class in $KK$-theory. (See Theorem 6.12 of \cite{Rr4},
and note that $KL (A, B)$ is in general a proper quotient of
$KK^0 (A, B).$)
Since the theorems we prove in Section 3 give information about
$KK$-theory rather than about R\o rdam's $KL$-theory, we introduce
and use the notion of \wue\  instead.  We give the definition
for \amm s because we will make extensive technical use of it in this
context, but, for reasons to be explained below, it is best suited to
\hm s.

\vspace{0.6\baselineskip}

{\bf  1.3.1 Definition.} Let $A$ and $D$ be \ca s, with $A$ separable.
Let $\ph, \ps : A \to D$ be two \amm s.
Then $\ph$ is {\em \wyue}\  to $\ps$ if there is a \ct\  family of
unitaries $t \mapsto u_{t}$ in $\tilde{D},$ defined
for $t \in [0, \infty),$ such that
\[
\lim_{t \to \infty} \|u_{t} \ph_{t} (a) u_{t}^* - \ps_{t} (a) \| = 0
\]
for all $a \in A.$ We say that two \hm s
$\ph, \ps : A \to D$ are {\em \wyue}\  if the corresponding
constant \amm s with $\ph_t = \ph$ and $\ps_t = \ps$ are \wyue.

\vspace{0.6\baselineskip}

{\bf  1.3.2 Lemma.} Asymptotic
unitary equivalence is the equivalence
relation on \amm s generated  by asymptotic equality and unitary
equivalence in the exact sense (that is,
$u_{t} \ph_{t} (a) u_{t}^* = \ps_{t} (a)$ for all $a \in A$).

\vspace{0.6\baselineskip}

{\em Proof:} The only point needing any work at all is
transitivity of \wue, and this is easy. \QED

\vspace{0.6\baselineskip}

{\bf  1.3.3 Lemma.} Let $A$ and $D$ be \ca s, with $A$ separable.

(1) Let $\ph, \ps : A \to \Kt  D$ be \wyue\  \amm s. Then
$\ph$ is homotopic to $\ps.$

(2) Let $\ph, \ps : A \to \Kt  D$ be \wyue\  \hm s. Then
$\ph$ is homotopic to $\ps$ via a path of \hm s.

\vspace{0.6\baselineskip}

{\em Proof:} (1) Let $t \mapsto u_t \in (\Kt D)^+$ be a \wue.
Modulo the usual isomorphism $M_2 (K) \cong K,$ the \amm s
$\ph$ and $\ps$ are homotopic to the \amm s
$\ph \oplus 0$ and $\ps \oplus 0$ from $A$ to $M_2 (\Kt D).$ Choose a
\ct\  function $(\af, t) \mapsto v_{\af, t}$ from
$[0, 1] \times [0, \infty)$ to $U(M_2 ((\Kt D)^+))$ such that
$v_{0, t} = 1$ and $v_{1, t} = u_t \oplus u_t^*$ for all $t.$
Define a homotopy of \amm s by
$\rh^{(\af)}_t (a) = v_{\af, t} (\ph_t (a)  \oplus 0) v_{\af, t}^*.$
Then $\rh^{(0)} = \ph \oplus 0$ and $\rh^{(1)}$ is
asymptotically equal to $\ps \oplus 0.$ So $\ph$ is homotopic to $\ps.$

(2) Apply the proof of part (1) the the constant paths
$t \mapsto \ph$ and $t \mapsto \ps.$ Putting $t = 0$ gives
homotopies of \hm s from $\ph$ to $\rh^{(1)}_0$ and from $\ps$ to
$\ps \oplus 0.$ The remaining piece of our homotopy is taken to be
defined for $t \in [0, \infty],$ and is given by
$t \mapsto \rh^{(1)}_t$ for $t \in [0, \infty)$ and
$\infty \mapsto \ps \oplus 0.$ \QED

\vspace{0.6\baselineskip}

{\bf  1.3.4 Corollary.}
Two \wyue\  \amm s define the same class in $E$-theory.

\vspace{0.6\baselineskip}

If the domain is nuclear, this corollary shows that
\wyue\  \amm s define the same class in $KK$-theory.
Asymptotic unitary equivalence thus rectifies the most important
disadvantage of \aue\  for \hm s.
Asymptotic unitary equivalence, however, also
has its problems, connected with the extension to \amm s.
The construction of the product of \amm s requires reparametrization
of \amm s, as in the following definition.

\vspace{0.6\baselineskip}

{\bf   1.3.5 Definition.}
Let $A$ and $D$ be \ca s, and let $\ph : A \to D$ be an \amm.
A {\em reparametrization} of $\ph$ is an \amm\  from
$A$ to $D$ of the form $t \mapsto \ph_{f (t)}$ for some
\ct\  nondecreasing function $f : [0, \infty) \to [0, \infty)$
such that $\limi{t} f (t) = \infty.$

\vspace{0.6\baselineskip}

Other versions are possible: one could replace ``nondecreasing''
by ``strictly increasing'', or omit this condition entirely.
The version we give is the most convenient for our purposes.

It is not in general true that an \amm\  is
\wyue\  to its reparametrizations.
(Consider, for example, the \amm\  $\ph : C(S^1) \to {\bf C}$ given
by $\ph_t (f) = f(\exp (it)).$)
The product is thus not
defined on \wue\  classes of \amm s. (The product is defined on
\wue\  classes when one factor is a \hm. We don't prove this fact
because we don't need it, but see the last part of the proof of
Lemma  2.3.5.)
In fact, if an \amm\  is \wyue\  to its reparametrizations, then it
is \wyue\  to a \hm, and this will play an important role in our
proof. The observation that this is true is due to Kirchberg.
It replaces a more complicated argument in the earlier version of this
paper, which involved the use
throughout of ``\wam s'', a generalization of \amm s in which there is
another parameter. We start the proof with a lemma.

\vspace{0.6\baselineskip}

{\bf   1.3.6 Lemma.}
Let $A$ and $D$ be \ca s, and let $\ph : A \to D$ be an \amm.
Suppose $\ph$ is \wyue\  to all its reparametrizations. Then for
any $\ep > 0$ and any finite set $F \subset A$ there is
$M \in [0, \infty)$ such that for any compact interval $I \subset \R$
and any \ct\  nondecreasing functions $f, \, g : I \to [M, \infty),$
there is a \ct\  unitary path $t \mapsto v_t$ in $\tilde{D}$
satisfying $\| v_t \ph_{f (t)} (a) v_t^* - \ph_{g (t)} (a) \| < \ep$
for all $t \in I$ and $a \in F.$

\vspace{0.6\baselineskip}

{\em Proof:}
Suppose the lemma is false. We can obviously change $I$
at will by reparametrizing, so there are $\ep > 0$ and $F \subset A$
finite such that for all $M \in [0, \infty)$ and all
compact intervals $I \subset \R$ there are
\ct\  nondecreasing functions $f, \, g : I \to [M, \infty)$ for which
no \ct\  unitary path $t \mapsto v_t$ in $\tilde{D}$ gives
$\| v_t \ph_{f (t)} (a) v_t^* - \ph_{g (t)} (a) \| < \ep$
for $t \in I$ and $a \in F.$
Choose $f_1$ and $g_1$ for $M = M_1 = 1$ and $I = I_1 = [1, 1 + \half].$
Given $f_n$ and $g_n,$ choose $f_{n + 1}$ and $g_{n + 1}$ as above
for $M = M_{n + 1} = 1 + \max (f_n (n + \half), g_n (n + \half))$
and $I = I_{n + 1} = [n + 1, n + 1 + \half].$
By induction, we have $M_n \geq n.$
Let $f, \, g : [0, \infty) \to  [0, \infty)$ be the unique continuous
functions which are linear on the intervals $[n + \half, n + 1]$
and
satisfy $f |_{[n, n + \half]} = f_n$ and $g |_{[n, n + \half]} = g_n.$
Since  $f$ and $g$ are nondecreasing and satisfy
$f (t), \, g (t) \geq n$ for $t \geq n,$ the functions
$t \mapsto \ph_{f (t)}$ and $t \mapsto \ph_{g (t)}$ are
\amm s which are reparametrizations of $\ph.$
By hypothesis, both are \wyue\  to $\ph,$ and are therefore
also \wyue\  to each other. Let $t \mapsto v_t$ be a unitary
path in $\tilde{D}$ which implements this \wue. Choose $T$ such that
for $a \in F$ and $t > T$ we have
$\| v_t \ph_{f (t)} (a) v_t^* - \ph_{g (t)} (a) \| < \ep/2.$
Restricting to $[n, n + \half]$ for some $n > T$ gives a contradiction
to the choice of $M$ and $\ep.$ This proves the lemma.
\QED

\vspace{0.6\baselineskip}

{\bf  1.3.7 Proposition.} Let $A$ be a separable \ca,
and let $\ph : A \to D$ be a bounded \amm. Suppose that
$\ph$ is \wyue\  to all its reparametrizations.
Then $\ph$ is \wyue\  to a \hm. That is, there exist a
\hm\  $\om : A \to D$ and a \ct\  path $t \mapsto v_t$ of
unitaries in $\tilde{D}$ such that for every $a \in A$, we have
$\lim_{t \to \infty} v_t \ph_t (a) v_t^* = \om (a).$

\vspace{0.6\baselineskip}

{\em Proof:}
Choose finite sets $F_0 \subset F_1 \subset \cdots \subset A$
whose union is dense in $A.$ Choose a sequence
$t_0 < t_1 < \cdots ,$ with $t_n \to \infty,$ such that
\[
\| \ph_t (ab) - \ph_t (a) \ph_t (b) \|, \,\,
   \| \ph_t (a^*)  - \ph_t (a)^* \| < 1 / 2^n
\]
for $a, b \in F_n$ and $t \geq t_n,$
and also such that, as in the previous lemma, for
any compact interval $I \subset \R$
and any \ct\  nondecreasing functions $f, \, g : I \to [t_n, \infty),$
there is a \ct\  unitary path $t \mapsto v_t$ in $\tilde{D}$
satisfying
$\| v_t \ph_{f (t)} (a) v_t^* - \ph_{g (t)} (a) \| < 2^{-n-1}$
for all $t \in I$ and $a \in F_n.$
For $n \geq 0$ let $t \mapsto u_t^{(n)}$ be the unitary path associated
with the particular choices $I = [t_n, t_{n + 1}],$ $f (t) = t,$ and
$g (t) = t_n.$ Set $\tilde{u}_t^{(n)} = (u_{t_n}^{(n)})^* u_t^{(n)}.$
We have
$\| ({u}_{t_n}^{(n)})^* \ph_{t_n} (a) {u}_{t_n}^{(n)} -
                                    \ph_{t_n} (a) \| < 2^{-n-1}$
for $a \in F_n,$ so
$\| \tilde{u}_t^{(n)} \ph_t (a) (\tilde{u}_t^{(n)})^* -
                                    \ph_{t_n} (a) \| < 2^{-n}$
for $t \in [t_n, t_{n + 1}]$ and $a \in F_n.$
Also note that $\tilde{u}_{t_n}^{(n)} = 1.$
Now define a \ct\  unitary function $[0, \infty) \to \tilde{D}$ by
\[
v_t = \tilde{u}_{t_1}^{(0)} \cdot \tilde{u}_{t_2}^{(1)} \cdots
     \tilde{u}_{t_n}^{(n -1)} \cdot \tilde{u}_t^{(n)}
\]
for $t_n \leq t \leq t_{n + 1}$.

We claim that $\om (a) = \lim_{t \to \infty} v_t \ph_t (a) v_t^*$
exists for all $a \in A.$ Since
$\sup_{t \in [0, \infty)} \| \ph_t \| < \infty,$ it suffices to check
this on the dense subset $\bigcup_{k = 0}^{\infty} F_k.$
So let $a \in F_k.$ We prove that the net
$t \mapsto v_t \ph_t (a) v_t^*$ is Cauchy. Let $m \geq k,$ and let
$t \geq t_m.$
Choose $n$ such that $t_n  \leq t \leq t_{n + 1}.$ Then
\beqr
\lefteqn{ \|v_t \ph_t (a) v_t^* - v_{t_m} \ph_{t_m} (a) v_{t_m}^* \|} \\
 & = &
 \left\| \left[ \rule{0em}{2.2ex}
   \tilde{u}_{t_{m + 1}}^{(m)} \cdot \tilde{u}_{t_{m + 2}}^{({m + 1})}
    \cdots
     \tilde{u}_{t_n}^{(n -1)} \cdot \tilde{u}_t^{(n)}
                                                    \right]
        \ph_t (a)
 \left[ \rule{0em}{2.2ex}
   \tilde{u}_{t_{m + 1}}^{(m)} \cdot \tilde{u}_{t_{m + 2}}^{({m + 1})}
    \cdots
     \tilde{u}_{t_n}^{(n -1)} \cdot \tilde{u}_t^{(n)}
                                                    \right]^*
 - \ph_{t_m} (a) \right\|                   \\
 & \leq &
 \left\| \left(\rule{0em}{2.2ex} \tilde{u}_t^{(n)} \right)
        \ph_t (a) \left( \rule{0em}{2.2ex} \tilde{u}_t^{(n)} \right)^*
                          - \ph_{t_n} (a)  \right\|
 + \sum_{j = m}^{n - 1}
   \left\| \left(\rule{0em}{2.2ex} \tilde{u}_{t_{j + 1}}^{(j)} \right)
            \ph_{t_{j + 1}} (a)
         \left(\rule{0em}{2.2ex} \tilde{u}_{t_{j + 1}}^{(j)} \right)^*
                               - \ph_{t_j} (a)  \right\|
  \leq  \sum_{j = m}^{n} \frac{1}{2^j} < \frac{1}{2^{m - 1}}.
\eeqr
Therefore, if $r, t \geq t_m,$ we obtain
\[
\| v_r \ph_r (a) v_r^* - v_t \ph_t (a) v_t^* \| <
                  1 / 2^{m - 2}.
\]
So we have a Cauchy net, which must converge. The claim is now proved.

Since $\ph_t$ is multiplicative and *-preserving to within $2^n$
on $F_n$ for $t \geq t_n,$ it follows that $\om$ is
exactly multiplicative and *-preserving on each $F_n.$
Since $\| \om \| \leq \sup_{t \in [0, \infty)} \| \ph_t \| < \infty,$ it
follows that $\om$ is a \hm.  \QED

\vspace{0.6\baselineskip}

In the rest of this section, we prove some useful facts about \wue.

\vspace{0.6\baselineskip}

{\bf  1.3.8 Lemma.}
Let $\ph : A \to D$ be an \amm, with  $A$ unital.
Then there is a \pj\  $p \in D$ and an
\amm\  $\ps : A \to D$ which is \wyue\  to $\ph$ and
satisfies $\ps_{t} (1) = p$ and
$\ps_{t} (a) \in p D p$ for all $t \in [0, \infty)$ and
and $a \in A.$

\vspace{0.6\baselineskip}

{\em Proof:} Let $t \mapsto q_{t}$ be a tail \pj\  for
$\ph (1),$ as in Definition  1.2.1. Standard results yield a
\ct\  family of unitaries $t \mapsto u_{t}$ in $\tilde{D}$ such that
$u_{0}  = 1$ and $u_{t} q_{t} u_{t}^* =  q_{0}$ for all
$t \in [0, \infty).$
Define $ p  = q_{0}$ and define
$\rh_{t} (a) = u_{t} q_{t} \ph_{t} (a) q_{t} u_{t}^*$
for $t \in [0, \infty)$ and $a \in A.$
Note that the definition of an \amm\  implies that
$(t, a) \mapsto q_{t} \ph_{t} (a) q_{t}$ is
asymptotically equal to $\ph,$ and hence is an \amm. Thus $\rh$ is an
\amm\  which is \wyue\  to $\ph.$

The only problem is that $\rh_{t} (1)$ might not be equal to
$p$. We do know that $\rh_{t} (1) \to p$ as
$t \to \infty.$ Choose a closed subspace $A_0$ of
$A$ which is complementary to ${\bf C} \cdot 1,$ and for $a \in A_0$
and $\ld \in {\bf C}$ define
$\ps_{t} (a + \ld \cdot 1) = \rh_{t} (a) + \ld p.$
\QED

\vspace{0.6\baselineskip}

{\bf  1.3.9 Lemma.}
Let $\ph, \ps : A \to \Kt D$ be \amm s, with  $A$ and $D$
unital.  Suppose that there is a \ct\  family of unitaries
$t \mapsto u_{t}$ in the multiplier algebra $M( \Kt D)$ such
that
$\lim_{t \to \infty} \|u_{t} \ph_{t} (a) u_{t}^* - \ps_{t} (a) \| = 0$
for all $a \in A.$ Then $\ph$ is \wyue\  to $\ps.$

\vspace{0.6\baselineskip}

{\em Proof:}
We have to show that $u_{t}$ can be replaced by
$v_{t} \in (\Kt D)^+.$

Applying the previous lemma twice, and making the corresponding
modifications to the given $u_{t}$, we may assume that
$\ph_{t} (1)$ and $\ps_{t} (1)$ are projections
$p$ and $q$ not depending on $t,$ and that we always have
$\ph_{t} (a) \in  p Dp$ and $\ps_{t} (a) \in  q D q$.

We now want to reduce to the case $p = q.$ The hypothesis implies
that there is $t_0$ such that $\|u_{t_0} p u_{t_0}^* - q \| < 1/2.$
Therefore there is a unitary $w$ in $(\Kt D)^+$
such that $w u_{t_0} p u_{t_0}^* w^* = q.$
Now if $p,$ $q \in \Kt D$ are projections which are unitarily
equivalent in $M (\Kt D),$ then standard arguments show they are
unitarily equivalent in $(\Kt D)^+.$
Therefore
conjugating $\ph$ by $w u_{t_0}$ changes neither its \wue\  class nor
the validity of the hypotheses. We may thus assume \wolog\  that
$p = q.$

Now choose $t_1$ such that $t \geq t_1$ implies
$\|u_{t} p u_{t}^* - p \| < 1.$
Define a \ct\  family of unitaries by
\[
c_{t} = 1 - p +
   p u_{t} p (p u_{t}^* p u_{t} p)^{-1/2}
                          \in (\Kt D)^+
\]
for $t \geq t_1$. (Functional calculus is
evaluated in $p (\Kt D) p.$) For any $d \in p (\Kt D) p,$ we have
\beqr
\lefteqn{
\| c_{t} d c_{t}^* - u_{t} d u_{t}^* \| =
 \| p d p - c_{t}^* u_{t} p d p u_{t}^* c_{t} \|} \\
 & \leq &
 2 \| d \| \| p  - c_{t}^* u_{t} p \|
  \leq 2 \| d \| (\| u_{t} p - p u_{t} \| +
                  \| p  - c_{t}^* p u_{t} p \|).
\eeqr
The first summand in the last factor goes to 0 as $t \to \infty.$
Substituting definitions, the second summand
becomes $\| p  - (p u_{t}^* p u_{t} p)^{1/2} \|,$
which does the same. Since $\ph_{t} (a) \in p (\Kt D) p$
for all $a \in A,$ and since (using Lemma 1.2 of \cite{Sm}
for the first)
\[
\limsup_{t \to \infty}  \| \ph_{t} (a) \|
                         \leq \| a \|
  \andeqn
\lim_{t \to \infty}
          \|u_{t} \ph_{t} (a) u_{t}^* - \ps_{t} (a) \| = 0,
\]
it follows that
$\lim_{t \to \infty}
          \|c_{t} \ph_{t} (a) c_{t}^* - \ps_{t} (a) \| = 0$
as well. This is the desired \wue. \QED

\vspace{\baselineskip}
\section{Asymptotic morphisms to tensor products with $\OI.$}

\vspace{\baselineskip}

The purpose of this section is to prove two things about \amm s from
a separable nuclear unital simple \ca\  $A$ to a \ca\  of the form
$\KOI{D}$ with $D$ unital: homotopy implies \wue, and each such
\amm\  is \wyue\  to a \hm. The basic method is the absorption
technique used in \cite{LP2} and \cite{Ph2}, and in fact this section
is really just the generalization of \cite{Ph2} from \hm s and
\aue\  to \amm s and \wue.

There are three subsections. In the first, we
collect for reference various known results involving Cuntz algebras
(including in particular Kirchberg's theorems on tensor products) and
derive some easy consequences. In the second subsection, we replace
\aue\  by \wue\  in the results of \cite{Rr1} and \cite{LP2}. In the
third, we carry out the absorption argument and derive its consequences.

The arguments involving \wue\  instead of \aue\  are sometimes
somewhat technical.
However, the essential outline of the proof is the same as in the
much easier to read paper \cite{Ph2}.

\vspace{0.6\baselineskip}
\subsection{Preliminaries: Cuntz algebras and Kirchberg's stability
theorems}
\vspace{0.6\baselineskip}

In this subsection, we collect for convenient reference various
results related to Cuntz algebras. Besides R\o rdam's results on
\aue\  and Kirchberg's basic results on tensor products, we need
material on unstable $K$-theory and hereditary subalgebras of tensor
products with $\OI$ and on exact stability of generating relations of
Cuntz algebras.

We start with R\o rdam's work \cite{Rr1}; we also use this
opportunity to establish our notation. The first definition is
used implicitly by  R\o rdam, and appears explicitly in the
work of Ringrose.

\vspace{0.6\baselineskip}

{\bf 2.1.1 Definition.} (\cite{Rn}, \cite{PR})
Let $A$ be a unital \ca. Then its {\em ($C^*$) \exl}\  $\cel (A)$ is
\[
\sup_{u \in U_0 (A)} \inf \left\{\sum_{k = 1}^n \| h_k\|: n \in {\bf N},
          h_1, \dots, h_n \in A \,\, {\rm selfadjoint},
 \,\,  u = \prod_{k = 1}^n \exp(i h_k) \right\}.
\]

\vspace{0.6\baselineskip}

In preparation for the following theorem, and to establish notation,
we make the following remark, most of which is in \cite{Rr1}, 3.3.

\vspace{0.6\baselineskip}

{\bf 2.1.2    Remark.}
Let $B$ be a unital \ca, and let $m \geq 2.$

(1) If $\ph, \ps : \OA{m} \to B$ are unital \hm s, then the element
$u = \sum_{j = 1}^m \ps (s_j) \ph (s_j)^*$ is a unitary in $B$ such
that $u \ph (s_j) = \ps (s_j)$ for $1 \leq j \leq m.$

(2) If $\ph : \OA{m} \to B$ is a unital \hm, then the formula
\[
\ld_{\ph} (a) = \sum_{j = 1}^m \ph (s_j) a \ph (s_j)^*
\]
defines a unital endomorphism  $\ld_{\ph}$ (or just $\ld$ when
$\ph$ is understood) of $B$.

(3) If $\ph$ and $\ld$ are as in (2), and if $u \in B$ has the form
$u = v \ld (v^*)$ for some unitary $v \in B$, then
$v \ph (s_j) v^* = u \ph (s_j)$ for $1 \leq j \leq m.$

\vspace{0.6\baselineskip}

{\bf    2.1.3 Theorem.}
Let $B$ be a unital \ca\  such that $\cel (B)$ is finite and such that
the canonical map $U (B) / U_0 (B) \to K_1 (B)$ is an isomorphism.
Let $m \geq 2,$ and let $\ph, \ps : \OA{m} \to B$, $\ld : B \to B,$
and $u \in U(B)$ be as in Remark 2.1.2    (1) and (2).
Then \tfae:

(1) $[u] \in (m - 1) K_1 (B).$

(2) For every $\ep > 0$ there is $v \in U(B)$ such that
$\| u - v \ld (v^*) \| < \ep.$

(3) $[\ph] = [\ps]$ in $KK^0 (\OA{m}, B).$

(4) The maps $\ph$ and $\ps$ are \ayue.

\vspace{0.6\baselineskip}

{\em Proof:} For $m$ even, this is Theorem 3.6 of \cite{Rr1}.
In Section 3 of \cite{Rr1}, it is also proved that (1) is equivalent
to (3) and (2) is equivalent to (4) for arbitrary $m,$ and Theorem
4.2 of \cite{Ph2} implies that (3) is equivalent to (4) for arbitrary
$m.$ \QED

\vspace{0.6\baselineskip}

We will not actually need to use the equivalence of (3) and (4) for
odd $m.$

That $\cel (D)$ is finite for \pisca s $D$ was first proved in
\cite{Ph1}. We will, however, apply this theorem to
algebras $D$ of the form $\OIA{B}$ with $B$ an arbitrary unital \ca.
Such algebras are shown in Lemma   2.1.7 (2) below to have
finite \exl. Actually, to prove the classification theorem, it suffices
to know that  there is a universal upper bound on
$\cel (C(X) \otimes B)$ for $B$ purely infinite and simple. This follows
from Theorem 1.2 of \cite{Zh1}.

We now state the fundamental results of Kirchberg on which our
work depends.
These were stated in \cite{Kr1}; proofs appear in \cite{KP}.

\vspace{0.6\baselineskip}

{\bf    2.1.4 Theorem.} (\cite{Kr1}; \cite{KP}, Theorem 3.7)
Let $A$ be a separable nuclear unital simple
\ca. Then $\OT{A} \cong \OA{2}.$

\vspace{0.6\baselineskip}

{\bf    2.1.5 Theorem.} (\cite{Kr1}; \cite{KP}, Theorem 3.14)
Let $A$ be a \kfalg.  Then $\OIA{A} \cong A.$

\vspace{0.6\baselineskip}

We now derive some consequences of Kirchberg's results.

\vspace{0.6\baselineskip}

{\bf  2.1.6 Corollary.}
Every \kfalg\  is approximately divisible in the sense of \cite{BKR}.

\vspace{0.6\baselineskip}

{\em Proof:} It suffices to show that $\OI$ is approximately divisible.
Let $\ph : \OIA{\OI} \to \OI$ be an isomorphism, as in the previous
theorem.
Define $\ps : \OI \to \OI$ by $\ps (a) = \ph (1 \otimes a).$
Then $\ps$ is \ayue\  to $\id_{\OI}$ by Theorem 3.3 of \cite{LP2}.
That is, there are unitaries $u_n \in \OI$ such that
$u_n \ph (1 \otimes a) u_n^* \to a$ for all $a \in \OI.$
Let $B \subset \OI$ be a unital copy of $M_2 \oplus M_3.$
Then for large enough $n,$ the subalgebra $u_n \ph (B \otimes 1) u_n^*$
of $\OI$ commutes arbitrarily well with any finite subset of $\OI.$
\QED

\vspace{0.6\baselineskip}

{\bf   2.1.7 Lemma.}
Let $D$ be any unital \ca. Then:

(1) The canonical map $U(\OIA{D})/ U_0(\OIA{D}) \to K_1 (\OIA{D})$ is
an isomorphism.

(2) $\cel (\OIA{D}) \leq 3 \pi.$

\vspace{0.6\baselineskip}

{\em Proof:}
We first prove surjectivity in (1). Let $\et \in K_1 (\OIA{D}).$
Choose $n$ and $u \in U (M_n (\OIA{D}))$ such that $[u] = \et.$
It is easy to find a (nonunital) \hm\  $\ph : M_n (\OI) \to \OI$
which sends $\diag(1, 0, \dots, 0)$ to a \pj\  $p \in \OI$ with
$[p] = [1]$ in $K_0 (\OI).$ Then $\ph$ is an isomorphism on
$K$-theory, so the Five Lemma and the
K\"unneth formula \cite{Sc1} show that
$\ph \otimes \id_D$ is too. Therefore $\ph (u) + 1 - \ph (1)$ is a
unitary in $\OIA{D}$ whose class is $\et.$

Now let $u \in U (\OIA{D})$ satisfy $[u] = 0$ in $K_1 (\OIA{D}).$
We prove that $u$ can be connected to the identity by a path of length
at most $3 \pi + \ep.$ This will simultaneously prove (2) and
injectivity in (1).

Using approximate divisibility of $\OI$ and approximating
$u$ by finite sums of elementary tensors, we can find nontrivial
\pj s $e \in \OI$ with $\|u (e \otimes 1) - (e \otimes 1) u \|$
arbitrarily small. If this norm is small enough, we can find a
unitary $v \in \KOI{D}$ which commutes with $e \otimes 1$ and is
connected to $u$ by a unitary path of length less that $\ep/2.$
Write $v = v_1 + v_2$ with
\[
v_1 \in U( e \OI e  \otimes D) \andeqn
   v_2 \in  U((1 - e) \OI (1 - e) \otimes D).
\]
Choose a partial
isometry $s \in \OI$ with $s^* s = 1 - e$ and $s s^* \leq e.$
The proof of Corollary 5
of \cite{Ph1} shows that $v$ can be connected to the unitary
\beqr
w & = &
 v \left[ \rule{0em}{3ex}
      (e - s s^*) \otimes 1 + (s \otimes 1) v_2 (s \otimes 1)^*
                                                      + v_2^* \right]
                                   \\
 & =  & 1 - e \otimes 1 + v_1 \left[ \rule{0em}{3ex}
   (s \otimes 1) v_2 (s \otimes 1)^* + (e - s s^*) \otimes 1 \right].
\eeqr
by a path of length $\pi$.

Since $\OI$ is purely infinite, there is an embedding of $\Kt e \OI e$
in $\OI$ which extends the obvious identification of
$e_{11} \otimes e \OI e$ with $e \OI e.$
It extends to a
unital \hm\  $\ph : (\Kt e \OI e \otimes D)^+ \to \OIA{D}$ whose range
contains $w,$ and such that $[\ph^{-1} (w)] = 0$ in
$K_1 (\Kt e \OI e \otimes D).$
Thus $\ph^{-1} (w) \in U_0 ((\Kt e \OI e \otimes D)^+).$ Theorem 3.8
of \cite{Ph1.5} shows that the $C^*$ exponential rank of any stable
\ca\  is at most $2 + \ep.$ An examination of the proof, and of the
length of the path used in the proof of Corollary 5 of \cite{Ph1},
shows that in fact any stable \ca\  has \exl\  at most $2 \pi.$
Thus, in particular, $\ph^{-1} (w)$ can be connected to $1$ by a unitary
path of length $2 \pi + \ep / 2.$ It follows that $u$ can be connected
to $1$ by a unitary path of length at most $3 \pi + \ep.$  \QED

\vspace{0.6\baselineskip}

A somewhat more complicated argument shows that in fact
$\cel (\OIA{D}) \leq 2 \pi.$ Details will appear elsewhere \cite{Ph3}.

\vspace{0.6\baselineskip}

{\bf  2.1.8 Lemma.}
Let $D$ be a unital \ca. Then:

(1) Given two full \pj s in $\KOI{D},$ each is \mvn\  to a
subprojection of the other.

(2) If two full \pj s in $\KOI{D}$ have the same $K_0$-class, then
they are homotopic.

\vspace{0.6\baselineskip}

{\em Proof:} Taking direct limits, we reduce to the case that $D$ is
separable. Then $\OIA{D}$ is approximately divisible by
Corollary  2.1.6. It follows from Proposition 3.10 of
\cite{BKR} that two full \pj s in $\KOI{D}$ with the same $K_0$-class
are \mvn. Now (2) follows from the fact that Murray-von Neumann
equivalence implies homotopy in the stabilization of a unital \ca.

Part (1) requires slightly more work. Let $\cal P$ be the set of all
\pj s $p \in \OIA{D}$ such that there are two orthogonal projections
$q_1, q_2 \leq p,$ both \mvn\  to $1.$ One readily verifies that
$\cal P$ is nonempty and satisfies the conditions ($\Pi_1$)-($\Pi_4$)
on page 184 of \cite{Cu2}. Therefore, by \cite{Cu2}, the group
$K_0 (\OIA{D})$ is exactly the set of Murray-von Neumann equivalence
classes of projections in $\cal P$. Since projections in  $\cal P$
are full, Proposition 3.10 of \cite{BKR} now implies that every full
\pj\  is in  $\cal P$. Clearly (1) holds for projections in $\cal P$.
We obtain (1) in general by using the pure infiniteness of $\OI$ to
show that every full \pj\  in $\KOI{D}$ is \mvn\  to a (necessarily
full) \pj\  in $\OIA{D}.$ \QED

\vspace{0.6\baselineskip}

Next, we turn to exact stability. For $\OA{m},$ we need only the
following standard result:

\vspace{0.6\baselineskip}

{\bf  2.1.9 Proposition.} (\cite{LP2}, Lemma 1.3 (1))
For any integer $m,$ the defining relations for $\OA{m}$,
namely
$s_j^* s_j = 1$ and $\sum_{k = 1}^m s_k s_k^* = 1$
for $1 \leq j \leq m,$ are exactly stable.

\vspace{0.6\baselineskip}

We will also need to know about the standard extension
$E_m$ of $\OA{m}$ by the compact operators. Recall from
\cite{Cu1} that $E_m$ is the universal \ca\  on generators
$t_1, \dots, t_m$ with relations
$t_j^* t_j = 1$ and $(t_j t_j^*) (t_k t_k^*) = 0$
for $1 \leq j, k \leq m, \,\, j \neq k.$
Its properties are summarized
in \cite{LP2}, 1.1. In particular, we have $\dirlim E_m \cong \OI$
using the standard inclusions.

Exact stability of the generating relations for $E_m$ is known, but
we need the following stronger result, which can be thought of as a
finite version of exact stability for $\OI.$ Essentially, it says that
if elements approximately satisfy the defining relations for $E_m,$
then they can be perturbed in a functorial way to exactly
satisfy these relations, and that the way the first $k$ elements
are perturbed does not depend on the remaining $m - k$ elements.

Recently, Blackadar has proved that in fact $\OI$ is
semiprojective in the usual sense \cite{Bl2}.

\vspace{0.6\baselineskip}

{\bf   2.1.10 Proposition.}
For each $\dt \geq 0$ and $m \geq 2,$
let $E_m (\dt)$ be the universal unital
\ca\  on generators $t^{(m)}_{j, \dt}$ for $1 \leq j \leq m$
and relations
\[
\| (t^{(m)}_{j, \dt})^* t^{(m)}_{j, \dt} - 1 \| \leq \dt
                      \andeqn
  \left\|
     \left(\rule{0em}{3ex} t^{(m)}_{j, \dt} (t^{(m)}_{j, \dt})^* \right)
    \left(\rule{0em}{3ex}  t^{(m)}_{k, \dt} (t^{(m)}_{k, \dt})^* \right)
    \right\|
                                 \leq \dt
\]
for $j \neq k,$
and let $\kp^{(m)}_{\dt} : E_m (\dt) \to E_m$ be the
\hm\  given by sending $t^{(m)}_{j, \dt}$ to the corresponding standard
generator $t^{(m)}_j$ of $E_m.$
Then there are $\dt(2) \geq \dt(3) \geq \cdots > 0$,
nondecreasing functions $f_m : [0, \dt (m)] \to [0, \infty)$
with $\lim_{\dt \to 0} f_m (\dt) = 0$ for each $m,$ and
\hm s $\ph^{(m)}_{\dt} : E_m \to E_m (\dt)$ for
$0 \leq \dt \leq \dt (m)$,
satisfying the following properties:

(1) $\kp^{(m)}_{\dt} \circ \ph^{(m)}_{\dt} = \id_{E_m}.$

(2) $\| \ph^{(m)}_{\dt} (t^{(m)}_{j})
                            - t^{(m)}_{j, \dt} \| \leq f_m (\dt).$

(3) If $0 \leq \dt \leq \dt' \leq \dt (m),$ then the composite of
$\ph^{(m)}_{\dt'}$ with the canonical map from
$E_m (\dt')$ to $E_m (\dt)$ is $\ph^{(m)}_{\dt}$.

(4) Let $\io^{(m)}_{\dt} : E_m (\dt) \to E_{m + 1} (\dt)$
be the map given by
$\io^{(m)}_{\dt} (t^{(m)}_{j, \dt}) = t^{(m + 1)}_{j, \dt}$
for $1 \leq j \leq m.$ Then for $0 \leq \dt \leq \dt (m + 1)$
and $1 \leq j \leq m,$ we have
$\io^{(m)}_{\dt} (\ph^{(m)}_{\dt} (t^{(m)}_{j}))
                    = \ph^{(m + 1)}_{\dt} (t^{(m + 1)}_{j}).$

\vspace{0.6\baselineskip}

{\em Proof:} The proof of exact stability of $E_m,$ as sketched in the
proof of Lemma 1.3 (2) of \cite{LP2}, is easily seen to yield
\hm s satisfying the conditions demanded here. \QED

\vspace{0.6\baselineskip}

{\bf 2.1.11 Proposition.} Let $D$ be a unital \pisca.
Then any two unital \hm s from $\OI$ to $D$ are homotopic.
Moreover, if $\ph, \ps : \OI \to D$ are unital \hm s such that
$\ph (s_j) = \ps (s_j)$ for $1 \leq j \leq m,$ then there is a
homotopy $t \mapsto \rh_t$ such that $\rh_t (s_j) = \ph (s_j)$ for
$1 \leq j \leq m$ and all $t.$

\vspace{0.6\baselineskip}

{\em Proof:}
We prove the second statement; the first is the special case
$m = 0.$

We construct, by induction on $n \geq m,$ continuous paths
$t \mapsto \rh_t^{(n)}$ of unital \hm s from $\OI$ to $D,$
defined for $t \in [n, n+1]$ and satisfying the following
conditions:

(1) $\rh_n^{(n)} = \rh_n^{(n - 1)}.$

(2) $\rh_t^{(n)} (s_j) = \ps (s_j)$ for $t \in [n, n+1]$ and
$1 \leq j \leq n.$

(3) $\rh_{n + 1}^{(n)} (s_{n + 1}) = \ps (s_{n + 1}).$

(4) $\rh_m^{(m)} = \ph.$

\noindent
Then  we define $\rh_t =  \rh_t^{(n)}$  for $t \in [n, n+1].$
This gives a \ct\  path $t \mapsto \rh_t$ for $t \in [m, \infty),$
with $\rh_m = \ph.$ Furthermore, $\rh_t (s_j) \to  \ps (s_j)$ for
all $j;$ since the $s_j$ generate $\OI$ as a \ca,
standard arguments show that $\rh_t (a) \to  \ps (a)$ for
all $a \in \OI.$ We have therefore constructed the required
homotopy.

It remains to carry out the inductive construction. If we define
$\rh_m^{(m - 1)} = \ph,$ then we only have to worry about (1), (2),
and (3).

Suppose $\rh^{(n - 1)}$ is given. Let
$p =
  \sum_{j = 1}^{n - 1} \rh_n^{(n - 1)} (s_j) \rh_n^{(n - 1)} (s_j)^*,$
which is a \pj\  in $D.$
Then define
\[
e_0 = \rh_n^{(n - 1)} (s_n) \rh_n^{(n - 1)} (s_n)^*
\andeqn
e_1 = \ps (s_n) \ps (s_n)^*.
\]
Both $e_0$ and $e_1$ are proper \pj s in the
\pisca\  $(1 - p) D (1 - p)$
with $K_0$-class equal to $[1_D],$ so they are
homotopic. It follows that there is a unitary path $s \mapsto u_s$
in $(1 - p) D (1 - p)$ such that $u_0 = 1$ and $u_1 e_0 u_1^* = e_1.$
For $s \in [0, 1/3],$
define $\rh_{n + s}^{(n)} (s_j) = \rh_n^{(n - 1)} (s_j)$ for
$1 \leq j \leq n - 1$ and
$\rh_{n + s}^{(n)} (s_j) = u_{3s} \rh_n^{(n - 1)} (s_j)$ for
$j \geq n.$
This yields a homotopy of
\hm s  $\rh_{n + s}^{(n)} : \OI \to D$ such that the isometries
$\rh_{n + 1/3}^{(n)} (s_n)$ and $\ps (s_n)$ have the same range
\pj, namely $e_1,$ although they themselves are probably not equal.

By a similar argument, we extend the homotopy over
$[n + 1/3, n + 2/3]$ in such a way that $\rh_{n + s}^{(n)} (s_j)$
is constant for $s \in [n + 1/3, n + 2/3]$ and $1 \leq j \leq n,$
and so that $\rh_{n + 2/3}^{(n)} (s_{n + 1})$ and
$\ps (s_{n + 1})$ also have the same range \pj,
say $f.$

Now $e_1$ and $f$ are \mvn, so we can
identify $(e_1 + f) D (e_1 + f)$ with $M_2 (e_1 D e_1).$
Since
\[
w_1 =
\left( \begin{array}{cc}
\ps (s_n) \rh_{n + 2/3}^{(n)} (s_{n})^* &               0             \\
              0            & [\ps (s_n) \rh_{n + 2/3}^{(n)} (s_{n})^*]^*
\end{array} \right)
      \in U_0 (M_2 (e_1 D e_1)),
\]
there is a \ct\  path of unitaries $s \mapsto w_s$ in
$M_2 (e_1 D e_1),$ with $w_0 = 1$ and $w_1$ as given.
For $s \in [2/3, 1],$ we now define
$\rh_{n + s}^{(n)} (s_j) = \rh_n^{(n - 1)} (s_j)$ for
$j \neq n, n + 1,$ and
$\rh_{n + s}^{(n)} (s_j) = w_{3s - 2} \rh_{n + 2/3}^{(n)} (s_j)$ for
$j = n, n + 1.$ This is again a homotopy, and gives
$\rh_{n + 1}^{(n)} (s_j) = \ps (s_j)$ for $1 \leq j \leq n,$ as
desired. The induction step is complete.  \QED

\vspace{0.6\baselineskip}

{\bf 2.1.12 Corollary.} Let $D$ be any unital \ca, and let
$p \in \KOI{D}$ be a \pj. Then
$\OIA{ p (\KOI{D}) p} \cong p (\KOI{D}) p.$

\vspace{0.6\baselineskip}

{\em Proof:}
We may replace $p$ by any \mvn\  \pj. So
\wolog\  $p \leq e \otimes 1 \otimes 1$ for some \pj\  $e \in K.$
Using the pure infiniteness of $\OI,$ we can in fact require that
$e$ be a rank one \pj. That is, we may assume $p \in \OIA{D}.$

By Theorem  2.1.5, there is an isomorphism
$\dt : \OIA{\OI} \to \OI.$ Using it, we need only consider \pj s
$p \in \OIA{\OIA{D}}.$ By the previous proposition and Theorem    2.1.5,
$a \mapsto 1 \otimes \dt (a)$ is homotopic to $\id_{\OIA{\OI}}.$
Therefore such a \pj\  $p$ is homotopic to
$q = 1 \otimes (\dt \otimes \id_D) (p),$ and hence also \mvn\  to
$q.$ Now
\[
q (\OIA{\OIA{D}}) q \cong
   \OIA{ [(\dt \otimes \id_D) (p) ][\OIA{D}][ (\dt \otimes \id_D) (p)]},
\]
which is unchanged by tensoring with $\OI$ by Theorem    2.1.5. \QED

\vspace{0.6\baselineskip}

{\bf 2.1.13 Corollary.} Let $D$ be a unital \ca.
Then the hypotheses on $B$ in Theorem    2.1.3 are satisfied for
any unital corner of $\KOI{D}.$

\vspace{0.6\baselineskip}

{\em Proof:} Combine the previous corollary
and Lemma  2.1.7. \QED

\vspace{0.6\baselineskip}
\subsection{Asymptotic unitary equivalence of \hm s from Cuntz algebras}
\vspace{0.6\baselineskip}

In this subsection, we strengthen the main technical
theorems of \cite{Rr1}
(restated here as Theorem  2.1.3) and of \cite{LP2}, replacing
\aue\  by \wue\  in the conclusions. We use the strong versions to
obtain variants of several other known results in which we replace
sequences of \hm s by \ct\  paths.

The first lemma contains the
essential point in the strengthening of Theorem  2.1.3. Its proof
uses the original theorem in a sort of bootstrap argument. The
remaining results lead up to the strengthening of the main theorem
of \cite{LP2}. They are proved by modifying the proofs there.

\vspace{0.6\baselineskip}

{\bf    2.2.1 Lemma.}  (Compare with Theorem    2.1.3.)
Let $D_0$ be a unital \ca, and let $D = \OIA{D_0}.$
Let $m \geq 2,$ and let $t \mapsto \ph_t$ and
$t \mapsto \ps_t$,
for $t \in [0, \infty),$ be two \ct\  paths of unital \hm s from
$\OA{m}$ to $D$. Suppose that the unitary
$u_0 = \sum_{j = 1}^m \ps_0 (s_j) \ph_0 (s_j)^*$ satisfies
$[u_0] \in (m - 1) K_1 (D).$ Then $\ph$ and $\ps$ are \wyue.

\vspace{0.6\baselineskip}

{\bf    2.2.2 Lemma.} (Compare with Proposition 1.7 of \cite{LP2}.)
Let $D$ be a unital \pisca, with $[1] = 0$
in $K_0 (D).$ Let $t \mapsto \ph_t$ and $t \mapsto \ps_t$,
for $t \in [0, \infty),$ be two \ct\  paths of unital \hm s from
$\OI$ to $D$. Then $t \mapsto \ph_t$ and $t \mapsto \ps_t$ are \wyue.

\vspace{0.6\baselineskip}

We will actually only need Lemma    2.2.1 for $m = 2.$

The proofs of the two lemmas are messy. We do the first (which is
easier) in detail, and then describe the modifications needed for
the second.

\vspace{0.6\baselineskip}

{\em Proof of Lemma    2.2.1:}
Corollary 2.1.13 shows that both $D$ and $C([0, 1], D)$
satisfy the hypotheses of Theorem    2.1.3.

By transitivity of \wue, it suffices
to show that $t \mapsto \ph_t$ and $t \mapsto \ps_t$
are both \wyue\  to some constant path. Thus, \wolog\  $t \mapsto \ph_t$
is a constant path $\ph_t = \ph$ for all $t.$
Let $\ld : D \to D$ be $\ld_{\ph}$ as in Remark 2.1.2    (2).

Let $f : [0, \dt] \to [0, \infty)$ be a function associated with
the exact stability of $\OA{m}$ (Proposition  2.1.9) in the
same way the functions $f_m$ of Proposition   2.1.10 are associated with
the exact stability of $E_m.$

Choose $\ep_0' > 0$ with $f (\ep_0') < 1.$ Choose $\ep_0 > 0$ with
$\ep_0 < 1/2,$ and also so small that if $\om : \OA{m} \to A$ is a
unital \hm, and $a_1, \dots, a_m \in A$
satisfy $\| a_j - \om (s_j) \| < \ep_0,$
then the $a_j$ satisfy the relations for $\OA{m}$ to within $\ep_0',$
that is,
\[
\| a_j^* a_j - 1 \| < \ep_0' \andeqn
     \left\| \sum_{k = 1}^m a_k a_k^* - 1 \right\| < \ep_0'
\]
for $1 \leq j \leq m.$
Set
$u_0 = \sum_{j = 1}^m \ps_0 (s_j) \ph (s_j)^*$; this is the same
as the $u_0$ in the statement of the lemma, so its $K_1$-class is in
$(m - 1) K_1 (D).$ Theorem 2.1.3    therefore yields a unitary
$v_0^{(0)} \in D$ such that
$\| u_0 - v_0^{(0)} \ld (v_0^{(0)})^* \| < \ep_0.$
Define $v_t^{(0)} = v_0^{(0)}$ for all $t,$ and define
$\gm_t^{(0)} : \OA{m} \to D$
by $\gm_t^{(0)} (a) = (v_t^{(0)})^* \ps_t (a) v_t^{(0)}.$
Using Remark 2.1.2, we calculate:
\[
\| \ph (s_j) - \gm_0^{(0)} (s_j) \|
   =  \| v_0^{(0)} \ph (s_j) (v_0^{(0)})^* - \ps_0 (s_j) \|
   =  \| v_0^{(0)} \ld (v_0^{(0)})^* \ph (s_j) - u_0 \ph (s_j) \|
                                              < \ep_0
\]
for $1 \leq j \leq m.$

We now construct, by induction on $n,$ numbers
$\ep_n, \ep_n' > 0$ and \ct\  paths $t \mapsto v_t^{(n)}$
of unitaries in $D$ and $t \mapsto \gm_t^{(n)}$
of unital \hm s from $\OA{m} \to D$, for $t \in [0, \infty),$
such that $\ep_0,$ $\ep_0',$ $v_t^{(0)},$ and $\gm_t^{(0)}$ are as
already chosen, and:

(1) $\gm_t^{(n)} (a) = (v_t^{(n)})^* \gm_t^{(n - 1)} (a) v_t^{(n)}$
for $a \in \OA{m}$ and $t \in [0, \infty).$

(2) If $n \geq 1,$ then $v_t^{(n)} = 1$ for $t \leq n.$

(3) If $n \geq 1,$ then
$\| \ph (s_j) - \gm_t^{(n)} (s_j) \| < 2^{-n + 1}$
for $t \in [n - 1, n],$ and if $n \geq 0$ then
$\| \ph (s_j) - \gm_t^{(n)} (s_j) \| < \ep_n$ for $t = n.$

(4) $f (\ep_n') < 2^{-n}.$

(5) Whenever $\om : \OA{m} \to A$ is a
unital \hm, and $a_1, \dots, a_m \in A$
satisfy $\| a_j - \om (s_j) \| < \ep_n,$
then the $a_j$ satisfy the relations for $\OA{m}$ to within $\ep_n'.$

(6) $\ep_n < 2^{-(n + 1)}.$

Suppose that $\ep_n,$ $\ep_n',$ $v_t^{(n)},$ and $\gm_t^{(n)}$ have
been chosen.
Choose $\ep_{n + 1}'$ and then $\ep_{n + 1}$ as in (4), (5), and (6).

For $\af \in [0, 1],$ define
\[
a_j (\af) = ( 1 - \af) (\ph (s_j) - \gm_n^{(n)} (s_j) )
                              + \gm_{n + \af}^{(n)} (s_j).
\]
Then $\| a_j (\af) - \gm_{n + \af}^{(n)} (s_j) \| < \ep_n$
for $1 \leq j \leq m$ and
$\af \in [0, 1]$. Conditions (4) and (5), and the choice of $f,$ provide
a unital \hm\ $\sm : \OA{m} \to C([0, 1], D)$ such that
$\| \sm (s_j) - a_j \| <  2^{-n}$ for $1 \leq j \leq m$.
Define $\sm_{\af} : \OA{m} \to D$ by $\sm_{\af} (a) = \sm (a) (\af)$
for $\af \in [0, 1]$ and $a \in  \OA{m}.$ Then
\[
\| \sm_{\af} (s_j) - \gm_{n + \af}^{(n)} (s_j) \| < \ep_n + 2^{-n}.
\]
Functoriality of the approximating \hm s (the analog of (3) of
Proposition   2.1.10) guarantees that $\sm_0 = \ph$ and
$\sm_1 =\gm_{n + 1}^{(n)}$.

Define a unitary $z \in C([0, 1], D)$ by
$z_{\af} = \sum_{j = 1}^m \sm_{\af} (s_j) \ph (s_j)^*$
for $\af \in [0, 1]$.
Note that $z_0 = 1,$ so $z \in U_0 (C([0, 1], D)).$ Theorem
   2.1.3 provides a unitary $\af \mapsto y_{\af}$ in $C([0, 1], D)$
such that $\| z_{\af} - y_{\af} \ld (y_{\af})^* \| < \ep_{n + 1} / 2$
for $\af \in [0, 1]$.
Putting $\af = 0,$ using $z_0 = 1,$ and rearranging terms,
we obtain $\| y_0^* \ld (y_0) - 1 \|  < \ep_{n + 1} / 2.$
Now define
\[
v_t^{(n + 1)} = \left\{ \begin{array}{ll}
     1                         & t \leq n                \\
     y_{t - n} y_0^*           & n \leq t \leq n + 1     \\
     y_1 y_0^*                 & n + 1 \leq t
\end{array} \right.
\]
and define
$\gm_t^{(n + 1)} (a) = (v_t^{(n + 1)})^* \gm_t^{(n)} (a) v_t^{(n + 1)}.$

It remains only to verify condition (3) in the induction hypothesis.
For $\af \in [0, 1],$
\beqr
  \|z_{\af} - v_{n + \af}^{(n + 1)} \ld ( v_{n + \af}^{(n + 1)})^* \|
 & \leq & \| z_{\af} - y_{\af} \ld (y_{\af})^* \| +
           \| y_{\af} \| \|1 - y_0^* \ld (y_0) \| \| \ld (y_{\af})^* \|
            \\
 &  <  &  \ep_{n + 1} / 2 +  \ep_{n + 1} / 2 =  \ep_{n + 1}.
\eeqr
Therefore, for $t \in [n, n + 1],$ Remark 2.1.2    yields
\beqr
\lefteqn{
   \| \ph (s_j) - \gm_t^{(n + 1)} (s_j) \|
   = \| v_t^{(n + 1)} \ph (s_j) (v_t^{(n + 1)})^* - \gm_t^{(n)} (s_j) \|
      }       \\
  & \leq & \| v_t^{(n + 1)} \ld (v_t^{(n + 1)})^* \ph (s_j)
                                             - z_{t - n} \ph (s_j) \|
                   + \| \sm_{t - n} (s_j) - \gm_t^{(n)} (s_j) \|
            \\
  &  <  &  \ep_{n + 1} + \ep_n + 2^{-n}
            < 2^{-(n + 2)} + 2^{-(n + 1)} + 2^{-n} < 2^{-n + 1}.
\eeqr
Furthermore, if $t = n + 1,$ then actually
$\sm_{t - n} (s_j) =  \gm_t^{(n)} (s_j),$ and we obtain
$\| \ph (s_j) - \gm_t^{(n + 1)} (s_j) \| < \ep_{n + 1}.$
This completes the induction.

To complete the proof, we now define
$v_t = \prod_{n = 0}^{\infty} v_t^{(n)}$
for $t \in [0, \infty).$ Note that the product defines a
\ct\  unitary path $t \mapsto v_t,$ since all but the first $n + 1$
factors are $1$ on $[0, n).$ Furthermore, for $t \in [n, n + 1],$
we have
\[
\| \ph (s_j) - v_t^* \ps_t (s_j) v_t \|
            = \| \ph (s_j) - \gm_t^{(n + 1)} (s_j) \| < 2^{-n + 1}.
\]
This implies that $\ph$ and $t \mapsto \ps_t$ are \wyue. \QED

\vspace{0.6\baselineskip}

For the proof of Lemma    2.2.2, we need the following lemma.

\vspace{0.6\baselineskip}

{\bf    2.2.3 Lemma.} Let $D$ be a unital \pisca\ with $[1] = 0$
in $K_0 (D).$ Let $m < n,$ and identify $E_m$ with the
subalgebra of $\OA{n}$ generated by $s_1, \dots, s_m.$ Let
$\ph : E_m \to D$ be a unital \hm. Then there exists a unital
\hm\  $\tilde{\ph} : \OA{n} \to D$ such that
$\tilde{\ph} |_{E_m} = \ph.$ Moreover, if we are already given a
unital \hm\  $\ps : \OA{n} \to D$, then $\tilde{\ph}$ can be chosen
to satisfy $[\tilde{\ph}] = [\ps]$ in $KK^0 (\OA{n}, D).$

\vspace{0.6\baselineskip}

{\em Proof:} This is essentially contained in the proof of
Proposition 1.7 of \cite{LP2}, using the equivalence of conditions
(1) and (3) in Theorem    2.1.3. \QED

\vspace{0.6\baselineskip}

{\em Proof of Lemma    2.2.2:} We describe how to modify the proof
of Lemma    2.2.1 to obtain this result.

First, note that $U(D) / U_0 (D) \to K_1 (D)$ is an isomorphism
because $D$ is purely infinite simple. Furthermore,
$\cel (C([0, 1], D)) \leq 5 \pi / 2 < \infty$ by Theorem 1.2
of \cite{Zh1}.
(It turns out that we only need this result for $D = \OI,$ so we could
use Corollary 2.1.13 here instead.)
Thus, the conditions on $D$ in Lemma    2.2.1
are satisfied.

As in the proof of Lemma    2.2.1, we may assume that
$t \mapsto \ph_t$ is a constant path $\ph_t = \ph$ for all $t.$

Let the functions $f_m$ be the ones associated with
the exact stability of $E_m$ as in Proposition   2.1.10.

The proof uses an induction argument similar to that of the proof
of Lemma    2.2.1, except that at the $n$-th stage we work with
extensions to $\OA{2n}$ of $\ph |_{E_n}$ and $\ps_t |_{E_n}$.
To avoid confusion, we let $s_1, s_2, \dots$ be the standard
generators of $\OI,$ with the first $n$ of them generating $E_n,$
and we let $s_1^{(2n)}, \dots, s_{2n}^{(2n)}$ be the standard
generators of $\OA{2n},$ with $E_k$, for $k < 2n,$ being identified
with the subalgebra generated by the first $k$ of them.

We start the construction at $n = 2$ so as not to have to worry about
$E_0$ and $E_1.$

In the preliminary step, we choose $\ep_2 > 0$ and $\ep_2' > 0$ so that
$\ep_2 < 1/8,$ $f_4 (\ep_2') < 1/4,$ and whenever $\om : E_2 \to A$
is a unital \hm, and
$a_1, a_2 \in A$ satisfy $\| a_j - \om (s_j) \| < \ep_0,$
then the $a_j$ satisfy the relations for $E_{2}$ to within $\ep_0'.$
Use Lemma    2.2.3 to choose unital \hm s
$\tilde{\ph}^{(2)}, \tilde{\ps}^{(2)}_2 : \OA{4} \to D$
such that $\tilde{\ph}^{(2)} |_{E_2} = \ph |_{E_2},$
$\tilde{\ps}^{(2)}_2 |_{E_2} = \ps_2 |_{E_2},$
and $[\tilde{\ph}^{(2)}] = [\tilde{\ps}^{(2)}_2]$ in $KK^0 (\OA{4}, D).$
Set
\[
u = \sum_{j = 1}^4
    \tilde{\ps}^{(2)}_2 (s_j^{(4)}) \tilde{\ph}^{(2)} (s_j^{(4)})^*.
\]
{}From (2) implies (3) in Theorem    2.1.3, we obtain a unitary
$v_2^{(2)} \in D$ such that
\[
\| u - v_2^{(2)} \ld_{\ph^{(2)}} (v_2^{(2)})^* \| < \ep_2.
\]
Define $v_t^{(2)} = v_2^{(2)}$ for $t \in [2, \infty),$ and define
$\gm_t^{(2)} : \OI \to D$ by
$\gm_t^{(2)} (a) = (v_t^{(2)})^* \ps_t (a) v_t^{(2)}.$
As in the proof of Lemma    2.2.1, a calculation shows that
\[
\| \tilde{\ph}^{(2)} (s_j^{(4)})
    - (v_2^{(2)})^* \tilde{\ps}_2^{(2)} (s_j^{(4)}) v_2^{(2)} \| < \ep_2
\]
for $1 \leq j \leq 4.$ It follows that
\[
\| \ph (s_j) - \gm_2^{(2)} (s_j) \| < \ep_2
\]
for $1 \leq j \leq 2.$

In the induction step, we now require that $t \in [2, \infty),$
that $\ep_2,$ $\ep_2',$ $\gm_t^{(2)},$ and $v_t^{(2)}$ be as already
given, that $\gm_t^{(n)} : \OI \to D$, and that:

(1) $\gm_t^{(n)} (a) =
   (v_t^{(n)})^* \gm_t^{(n - 1)} (a) v_t^{(n)}$
for $a \in \OI$ and $t \in [2, \infty).$

(2) If $n \geq 3,$ then $v_t^{(n)} = 1$ for $t \leq n.$

(3) If $n \geq 3,$ then
$\| \ph (s_j) - \gm_t^{(n)} (s_j) \| < 2^{-n + 1}$
for $t \in [n - 1, n]$ and $1 \leq j \leq n- 1,$ and if $n \geq 2$
then $\| \ph (s_j) - \gm_t^{(n)} (s_j) \| < \ep_n$
for $t = n$ and $1 \leq j \leq n.$

(4) $f_{n} (\ep_n') < 2^{-n}.$

(5) Whenever $\om : E_{n} \to A$ is a
unital \hm, and $a_1, \dots, a_n \in A$ satisfy
$\| a_j - \om (s_j) \| < \ep_n,$
then the $a_j$ satisfy the relations for $E_{n}$ to within $\ep_n'.$

(6) $\ep_n < 2^{-(n + 1)}.$

For the proof of the inductive step, we first
choose $\ep_{n + 1}'$ and $\ep_{n + 1}$ to satisfy (4), (5), and (6).
Then construct, as in the proof of Lemma    2.2.1, a \ct\  path of
\hm s $\sm_{\af} : E_n \to D$ such that $\sm_0 = \ph |_{E_n},$
$\sm_1 = \gm_{n + 1}^{(n)} |_{E_n},$ and
\[
\| \sm_{\af} (s_j ) - \gm_{n + \af}^{(n)} (s_j) \| < \ep_n + 2^{-n}
\]
for $1 \leq j \leq n.$

We now claim that there is a unitary path $\af \mapsto w_{\af}$
in $D$ such that $w_0 = 1,$
$w_{\af} \sm_{\af} (s_j) = \sm_0 (s_j)$
for $\af \in [0, 1]$ and $1 \leq j \leq n,$ and
$w_1 \gm_{n + 1}^{(n)} (s_{n + 1}) = \ph   (s_{n + 1}).$
To prove this, start by defining
$q_{\af} = \sum_{j = 1}^{n} \sm_{\af} (s_j) \sm_{\af} (s_j)^*.$
Then set
$w'_{\af} = \sum_{j = 1}^{n} \sm_0 (s_j) \sm_{\af} (s_j)^*,$
which is a partial isometry from $q_{\af}$ to $q_0$ such that
$w'_{\af} \sm_{\af} (s_j) = \sm_0 (s_j)$
for $1 \leq j \leq n.$ Next, define
\[
p_1 = \gm_{n + 1}^{(n)} (s_{n + 1}) \gm_{n + 1}^{(n)} (s_{n + 1})^*
   \andeqn    p_0 = \ph (s_{n + 1}) \ph (s_{n + 1})^*.
\]
Since
$\gm_{n + 1}^{(n)} |_{E_n} = \sm_1$ and
$\ph |_{E_n} = \sm_0,$ we see that $p_1$ and $p_0$ are proper
sub\pj s of $1 - q_1$ and $1- q_0$ respectively, both with the
same class (namely $[1] = 0$) in $K_0 (D).$ Standard methods therefore
yield a unitary path $\af \mapsto c_{\af}$ in $D$ such that
$c_0 = 1,$ $c_{\af} q_{\af} c_{\af}^* = q_0,$ and $c_1 p_1 c_1^* = p_0.$
Then $\ph (s_{n + 1}) \gm_{n + 1}^{(n)} (s_{n + 1})^* c_1^*$
is a unitary in $p_0 D p_0,$ so there is a unitary
$d \in (1 - q_0 - p_0 ) D (1 - q_0 - p_0 )$ such that
\[
\ph (s_{n + 1}) \gm_{n + 1}^{(n)} (s_{n + 1})^* c_1^*
                               + d \in U_0 ((1 - q_0 ) D (1 - q_0)),
\]
and a unitary path $\af \mapsto w_{\af}''$ in $(1 - q_0 ) D (1 - q_0)$
such that
\[
w_0'' = 1 \andeqn
     w_1'' = \ph (s_{n + 1}) \gm_{n + 1}^{(n)} (s_{n + 1})^* c_1^* + d.
\]
Set $w_{\af} = w_{\af}' + w_{\af}'' c_{\af};$ this is the path that
proves the claim.

Use Lemma  2.2.3 to choose a unital
\hm\  $\tilde{\ph}^{(n + 1)} : \OA{2n + 2} \to D$ such that
$\tilde{\ph}^{(n + 1)} |_{E_{n + 1}} = \ph |_{E_{n + 1}}.$ Define
unital \hm s $\tilde{\sm}_{\af} : \OA{2n + 2} \to D$ by
\[
\tilde{\sm}_{\af} (s_j^{(2n + 2)})
                      = w_{\af}^* \tilde{\ph}^{(n + 1)} (s_j^{(2n + 2)})
\]
for $1 \leq j \leq 2n + 2.$ Then
\[
\tilde{\sm}_0 = \tilde{\ph}^{(n + 1)},  \,\,\,\,\,\,
   \tilde{\sm}_{\af} |_{E_n} =\sm_{\af}, \andeqn
     \tilde{\sm}_1 |_{E_{n + 1}} = \gm_{n + 1}^{(n)} |_{E_{n + 1}}.
\]

Define $z$ and choose $y$ as in the proof of Lemma    2.2.1,
using $\OA{2 n + 2}$ in place of $\OA{m},$
$\tilde{\sm}$ in place of $\sm,$ $\tilde{\ph}^{(n + 1)}$ in
place of $\ph,$ and $\ld = \ld_{\tilde{\ph}^{(n + 1)}}.$
Define $v_t^{(n + 1)}$ and $\gm_t^{(n + 1)}$ as there. The same
computations as there show that
\[
\| \ph(s_j) - \gm_t^{(n + 1)} (s_j) \| =
  \|\tilde{\ph}^{(n + 1)} (s_j^{(2n + 2)})
          - (v_t^{(n + 1)})^* \tilde{\sm}_{t - n} (s_j^{(2n + 2)})
                                                       v_t^{(n + 1)} \|
           < 2^{-n + 1}
\]
for $1 \leq j \leq n$ and $t \in [n, n + 1],$ and
\[
\| \ph(s_j) - \gm_{n + 1}^{(n + 1)} (s_j) \| =
  \|\tilde{\ph}^{(n + 1)} (s_j^{(2n + 2)})
          - (v_t^{(n + 1)})^* \tilde{\sm}_{1} (s_j^{(2n + 2)})
                                                       v_t^{(n + 1)} \|
< \ep_{n + 1}
\]
for $1 \leq j \leq n + 1.$ This completes the induction step.

Define $v_t = \prod_{n = 2}^{\infty} v_t^{(n)}.$ Calculations
analogous to those in the proof of Lemma    2.2.1 show that
$t \mapsto v_t$ is a \ct\  unitary path in $D,$ and that for
$n \geq 2$ we have
\[
\| \ph(s_j) - v_t^* \ps_t (s_j) v_t \| < 2^{-n + 1}
\]
for $t \in [n, n + 1]$ and $1 \leq j \leq n.$ This implies that
\[
\lim_{t \to \infty}
         \left( \rule{0em}{3ex} \ph(a) - v_t^* \ps_t (a) v_t \right) = 0
\]
for all $a \in \OI.$ \QED

\vspace{0.6\baselineskip}

{\bf    2.2.4 Lemma.} There exists a \ct\  family
$t \mapsto \ph_t$ of unital endomorphisms of $\OI,$ for
$t \in [0, \infty),$ which is asymptotically central in the sense that
\[
\lim_{t \to \infty}
        \left( \rule{0em}{3ex} \ph_t (b) a - a \ph_t (b) \right) = 0
\]
for all $a, b \in \OI.$

\vspace{0.6\baselineskip}

{\em Proof:} Let $A_n$ be the tensor product of $n$ copies of $\OI,$
and define $\mu_n : A_n \to A_{n + 1}$ by $\mu_n (a) = a \otimes 1.$
Set $A = \dirlim A_n,$ which is just $\bigotimes_1^{\infty} \OI.$
Theorem    2.1.5 implies that $A_n \cong \OI,$ so Theorem 3.5 of
\cite{LP2} implies that $A \cong \OI.$
(Actually, that $A \cong \OI$ is shown in the course of the proof of
Theorem 2.1.5. See \cite{KP}.) It therefore suffices to
construct a \ct\  asymptotically central inclusion of $\OI$ in $A$
rather than in $\OI.$

Let $\nu_n : A_n \to A$ be the inclusion.
Proposition   2.1.11 provides a homotopy $\af \mapsto \ps_{\af}$ of
unital \hm s $\ps_{\af} : \OI \to \OIA{\OI}$ such that
$\ps_0 (a) = a \otimes 1$ and $\ps_1 (a) = 1 \otimes a.$ For $n \geq 1$
and $t \in [n, n + 1],$ we write $t = n + \af$ and define
\[
\ph_t (a) = \nu_{n + 2} (1  \otimes 1 \otimes \cdots
                      \otimes 1 \otimes \ps_{\af} (a)),
\]
where the factor $1$ appears $n$ times in the tensor product.
The two definitions of $\ph_n (a)$ agree, so $t \mapsto \ph_t$ is
\ct. We clearly have
$\lim_{t \to \infty} \left( \ph_t (b) a - a \ph_t (b) \right) = 0$
for $b \in\OI$ and $a \in \bigcup_{n = 1}^{\infty} \nu_n (A_n),$
and a standard argument then shows this is true for all $a \in A.$
\QED

\vspace{0.6\baselineskip}

The notation introduced in the following definition is the same as
in \cite{LP1}, \cite{LP2}, and \cite{Ph2}.

\vspace{0.6\baselineskip}

{\bf  2.2.5 Definition.} Let $A$ be any unital \ca ,
and let $D$ be a \pisca .
Let $\varphi, \psi: A \to D$ be two homomorphisms, and assume that
$\varphi(1) \neq 0$ and
$[\psi(1)] = 0$ in $K_0 (D)$. We define a homomorphism
$\varphi \tdsum \psi : A \to D$, well defined up to unitary equivalence,
by the following construction. Choose a projection $q \in D$ such that
$0 < q < \varphi(1)$ and $[q] = 0$ in $K_0 (D).$
Since $D$ is purely infinite and
simple, there are partial isometries $v$ and $w$ such that
\[
vv^* = \varphi(1) -q, \,\,\,\,\,\, v^*v = \varphi(1), \,\,\,\,\,\,
ww^* = q, \andeqn w^*w = \psi(1).
\]
Now define $(\varphi \tdsum \psi)(a) = v\varphi(a)v^* + w\psi(a)w^*$
for $ a \in A$.

\vspace{0.6\baselineskip}

{\bf  2.2.6 Lemma.}
(Compare with Proposition 2.3 of \cite{LP2}.)
Let $D$ be a unital \pisca, and let $q \in D$ be a \pj\  with
$[q] = 0$ in $K_0 (D).$ Let $\ph : \OI \to D$ and
$\ps : \OI \to qDq$ be unital \hm s. Then $\ph$ is \wyue\  to
$\ph \tdsum \ps.$

\vspace{0.6\baselineskip}

{\em Proof:} Let $t \mapsto \gm_t$ be a \ct ly parametrized
asymptotically central inclusion of $\OI$ in $\OI$, as in
Lemma    2.2.4. Let $e \in \OI$ be a nonzero \pj\  with
$[e] = 0$ in $K_0 (\OI),$ and set $e_t = \gm_t (e).$ Choose a
\ct\   unitary path $t \mapsto u_t$ such that $u_t e_t u_t^* = e_0.$

Let the functions $f_m : [0, \dt(m)] \to [0, \infty)$ be as in
Proposition   2.1.10. Choose numbers $\ep_2 > \ep_3 > \cdots 0$
and  $\ep'_2 > \ep'_3 > \cdots 0$ such that:

(1) $\ep_m' < \dt (m)$ and $f_{m} (\ep_m') < 1/m.$

(2) Whenever $\om : E_{m} \to A$ is a
unital \hm, and $a_1, \dots, a_m \in A$ satisfy
$\| a_j - \om (s_j) \| < \ep_m,$
then the $a_j$ satisfy the relations for $E_{m}$ to within $\ep_m'.$

(3) $\ep_m < 1/m.$

Next, use the asymptotic centrality of $t \mapsto e_t$ to choose
$t_2 < t_3 < \cdots,$ with $t_m \to \infty$ as $m \to \infty$,
such that
\[
\left\|s_j - \left[  \rule{0em}{2.2ex}
   e_t s_j e_t + (1 - e_t) s_j (1 - e_t)\right] \right\| < \ep_m
\]
for $1 \leq j \leq m$ and $t \geq t_m.$ Define
\[
a_j (t) = u_t  \left[ \rule{0em}{2.2ex}
                    e_t s_j e_t + (1 - e_t) s_j (1 - e_t)\right] u_t^*
        \in e_0 \OI e_0 \oplus (1 - e_0) \OI (1 - e_0).
\]
Conditions (1) and (2), and
Proposition   2.1.10, then yield \ct\  paths $t \mapsto \sm_t^{(m)}$
of \hm s from $E_m$ to $e_0 \OI e_0 \oplus (1 - e_0) \OI (1 - e_0),$
defined for $t \geq t_m,$ such that
$\| \sm_t^{(m)} (s_j) - a_j (t) \| < 1/m$
for $1 \leq j \leq m,$ and $\sm_t^{(m + 1)} |_{E_m} = \sm_t^{(m)}$
for $t \geq t_{m + 1}.$

Define
\[
\af_t^{(m)} : E_m \to e_0 \OI e_0 \andeqn
\bt_t^{(m)} : E_m \to (1 - e_0) \OI (1 - e_0)
\]
by
\[
\af_t^{(m)} (a) = e_0 \sm_t^{(m)} (a) e_0 \andeqn
     \bt_t^{(m)} (a) = (1 - e_0) \sm_t^{(m)} (a) (1 - e_0).
\]

Note that $\af_{t_m}^{(m)}$ is homotopic to
$\af_{t_{m + 1}}^{(m+ 1)} |_{E_m};$ since
$\af_{t_{m + 1}}^{(m+ 1)} |_{E_m}$ is injective, it follows that
$\af_{t_m}^{(m)}$ is injective. Since $e_0 \OI e_0$ is purely
infinite simple, it is easy to extend  $\af_{t_m}^{(m)}$ to a
\hm\  $\af_{t_m} : \OI \to e_0 \OI e_0$. Proposition   2.1.11
provides homotopies $t \mapsto \af_t$ of \hm s from $\OI$ to
$e_0 \OI e_0$, defined for $t \in [t_m, t_{m + 1}],$ such that
$\af_t |_{E_m} = \af_{t}^{(m)}$ and such that $\af_{t_m}$ and
$\af_{t_{m + 1}}$ are as already given. Putting these homotopies
together, and defining $\af_t = \af_{t_2}$ for $t \in [0, t_2],$
we obtain a \ct\  path $t \mapsto \af_t$ of unital \hm s from $\OI$ to
$e_0 \OI e_0$, defined for $t \in [0, \infty),$ such that
$\af_t |_{E_m} = \af_{t}^{(m)}$ whenever $t \geq t_m.$
Similarly, there is a \ct\  path $t \mapsto \bt_t$ of unital \hm s
from $\OI$ to $(1 - e_0) \OI (1 - e_0)$, defined for
$t \in [0, \infty),$ such that $\bt_t |_{E_m} = \bt_{t}^{(m)}$
whenever $t \geq t_m.$ Define $\sm_t : \OI \to \OI$ by
$\sm_t (a) = \af_t (a) + \bt_t (a).$

For $t \geq t_m$ and $1 \leq j \leq m,$ we have
$u_t^* \sm_t (s_j) u_t = u_t^* \sm_t^{(m)} (s_j) u_t,$ and
\[
\| u_t^* \sm_t^{(m)} (s_j) u_t - s_j \|  \leq
       \|\sm_t^{(m)} (s_j) - a_j (t) \| + \|u_t^* a_j (t) u_t - s_j \|
 < 1/m + \ep_m < 2/m.
\]
Therefore
$\lim_{t \to \infty} \| u_t^* \sm_t^{(m)} (s_j) u_t - s_j \| = 0$
for all $j.$ Thus $t \mapsto \sm_t$ is \wyue\  to $\id_{\OI}.$
So $\ph$ is \wyue\  to $t \mapsto \ph \circ \sm_t.$

Let $f < \ph (e_0)$ be a nonzero \pj\  with $[f] = 0$ in $K_0 (D).$
Let $w_1, w_2 \in D$ be partial isometries satisfying
\[
w_1^* w_1 = 1, \,\,\, w_1 w_1^* = 1 - f, \andeqn
         w_1 (1 - \ph (e_0)) = (1 - \ph (e_0)) w_1 =1 - \ph (e_0)
\]
and
\[
w_2^* w_2 = q \andeqn  w_2 w_2^* = f.
\]
The \hm\  $\ph \tdsum \ps$ is only defined  up to unitary equivalence,
and we can take it to be
\[
(\ph \tdsum \ps) (x) = w_1 \ph (x) w_1^* + w_2 \ps (x) w_2^*.
\]
We make the same choices when defining $(\ph \circ \sm_t) \tdsum \ps.$
Writing $\ph \circ \sm_t = \ph \circ \af_t + \ph \circ \bt_t,$
with
\[
\ph \circ \af_t : \OI \to \ph (e_0) D \ph (e_0) \andeqn
    \ph \circ \bt_t : \OI \to \ph (1 - e_0) D \ph( 1 -e_0),
\]
this choice gives
\[
(\ph \circ \sm_t) \tdsum \ps = [(\ph \circ \af_t) \tdsum \ps]
                                   + \ph \circ \bt_t.
\]
By Lemma    2.2.2, $t \mapsto (\ph \circ \af_t) \tdsum \ps$ is
\wyue\  to $\ph \circ \af_t.$ Therefore, with $\sim$ denoting
\wue, we have
\[
\ph \tdsum \ps  \sim
   (\ph \circ \af_t) \tdsum \ps + \ph \circ \bt_t \sim
   \ph \circ \af_t + \ph \circ \bt_t \sim \ph.
\]
This is the desired result. \QED

\vspace{0.6\baselineskip}

{\bf    2.2.7 Proposition.}
(Compare with Theorem 3.3 of \cite{LP2}.)
Let $D$ be a unital \pisca, and let
$\ph, \ps : \OI \to D$ be two unital \hm s. Then $\ph$ is \wyue\  to
$\ps.$

\vspace{0.6\baselineskip}

{\em Proof:} Let $e = 1 - s_1 s_1^* - s_2 s_2^* \in \OI,$ and let
$f = \ph (e) \in D.$ Define $\overline{\ph} : \OI \to f D f$ by
$\overline{\ph} (s_j) = \ph (s_{j + 2})f.$ Let $w \in M_2 (D)$ be a
partial isometry with
$w^* w = 1 \oplus f$ and $w w^* = q \oplus 0$ for some $q \in D.$
We regard $w (\ph \oplus \overline{\ph})(-) w^*$ and
$w (\ps \oplus \overline{\ph})(-) w^*$ as \hm s from $\OI$ to $qDq.$
Furthermore, $[q] = 0$ in $K_0 (D),$ so
\[
\ph \tdsum w (\ps \oplus \overline{\ph})(-) w^*   \andeqn
   \ps \tdsum w (\ph \oplus \overline{\ph})(-) w^*
\]
are defined; they are easily seen to be unitarily equivalent.
Using Lemma    2.2.6 for the other two steps, we therefore obtain \wue s
\[
\ph \sim \ph \tdsum w (\ps \oplus \overline{\ph})(-) w^* \sim
       \ps \tdsum w (\ph \oplus \overline{\ph})(-) w^*  \sim \ps.
\]
\QED

\vspace{0.6\baselineskip}

{\bf    2.2.8 Corollary.} Let $A$ be any unital \ca\  such that
$\OIA{A} \cong A.$ Then there exists an isomorphism
$\bt : \OIA{A} \to A$ such that the
\hm\  $a \mapsto \bt (1 \otimes a)$ is \wyue\  to $\id_A.$

\vspace{0.6\baselineskip}

{\em Proof:} We first prove this for $A = \OI.$ Theorem    2.1.5 implies
that $\OIA{\OI} \cong \OI;$ let
$\bt_0 : \OIA{\OI} \to \OI$ be an isomorphism.
Then $a \mapsto \bt_0 (1 \otimes a)$
and $\id_{\OI}$ are two unital \hm s from $\OI$ to $\OI,$ so they
are \wyue\  by Proposition    2.2.7. Let $t \mapsto u_t$ be a unitary
path such that
$\lim_{t \to \infty} (\bt_0 (1 \otimes a) - u_t a u_t^*) = 0$ for all
$a \in \OI.$

Now let $A$ be as in the hypotheses. We may as well prove the result for
$\OIA{A}$ instead of $A$.  Take $\bt = \bt_0 \otimes \id_A;$ then
$a \mapsto \bt (1 \otimes a)$ is \wyue\  to $\id_{\OIA{A}}$ via the
unitary path $t \mapsto u_t \otimes 1.$ \QED

\vspace{0.6\baselineskip}
\subsection{When homotopy implies \wue}
\vspace{0.6\baselineskip}

In this subsection, we will prove that if $A$ is a separable
nuclear unital simple \ca\  and
$D_0$ is unital, then two homotopic \amm s
from $A$ to $\Kt \OIA{D_0}$ are \wyue. We will furthermore prove that
an  \amm\  from $A$ to $\Kt \OIA{D_0}$ is \wyue\  to a \hm . The method
of proof of the first statement will generalize the methods of
\cite{Ph2}. We will obtain the second via a trick.

The following two definitions will be convenient. The first is used,
both here and in Section 3, to simplify terminology, and the second
is the analog of Definition 2.1 of \cite{Ph2}.

\vspace{0.6\baselineskip}

{\bf  2.3.1 Definition.} Let $A,$ $D,$  and $Q$ be \ca s, with
$A$ and $Q$ separable and with $Q$ also unital and nuclear.
Let $\ph : A \to D$ be an \amm.
A {\em standard factorization} of $\ph$ through $Q \otimes A$ is an
\amm\  $\ps : Q \otimes A \to D$ such that
$\ph_t (a) = \ps_t (1 \otimes a)$ for all $t$ and all $a \in A.$
An {\em asymptotic standard factorization}
of $\ph$ through $Q \otimes A$ is an
\amm\  $\ps : Q \otimes A \to D$ such that $\ph$ is \wyue\  to
the \amm\  $(t, a) \mapsto \ps_t (1 \otimes a)$.

\vspace{0.6\baselineskip}

{\bf  2.3.2 Definition.} Let $A,$ $D,$ and $\ph$ be as in the
previous definition. An
{\em (asymptotically) \tfn}\  of $\ph$ is a (asymptotic)
standard factorization with $Q = \OA{2}.$ In this case, we say that
$\ph$ is {\em (asymptotically) trivially factorizable.}

\vspace{0.6\baselineskip}

{\bf  2.3.3 Lemma.} (Compare \cite{Ph2}, Lemma 2.2.)
Let $A$ be \snus, let $D_0$ be a unital \ca, and let $D = \OIA{D_0}.$
Then any two full   \wtf\  \amm s from $A$ to $\Kt D$ are
\wyue.

\vspace{0.6\baselineskip}

{\em Proof:}
It suffices to prove this for full \amm s $\ph, \ps : A \to \Kt D$
with \tfn s $\ph', \ps' : \OT{A} \to \Kt D.$
Note that $\ph'$ and $\ps'$
are again full, and it suffices to prove that $\ph'$ is \wyue\  to
$\ps'.$ By Theorem    2.1.4, we have
$\OT{A} \cong \OA{2},$ and Proposition
1.1.7 then implies that $\ph'$ and $\ps'$ are
asymptotically equal to \ct\  families $\ph''$ and $\ps''$ of \hm s.

We now have two \ct\  families of \dpj s
$t \mapsto \ph''_{t} (1)$
and $t \mapsto \ps''_{t} (1)$ in $\Kt D,$ parametrized by
$[0, \infty).$ Standard methods show that each family is unitarily
equivalent to a constant projection.
Moreover, the \pj s $\ph''_0 (1)$ and $\ps''_0 (1)$
have trivial $K_0$ classes, so are homotopic by Lemma  2.1.8 (2).
Therefore they are unitarily equivalent. Combining the unitaries
involved and conjugating by the result,
we can assume $\ph''_{t} (1)$ and $\ps''_{t} (1)$
are both equal to the constant family $t \mapsto p$ for a
suitable \dpj\  $p$. Now replace $\Kt D$ by $p (\Kt D) p,$ and apply
Lemma  2.2.1; its hypotheses are satisfied by Corollary 2.1.12. \QED

\vspace{0.6\baselineskip}

{\bf  2.3.4 Corollary.} (Compare \cite{Ph2}, Lemma 2.3.) Under the
hypotheses of Lemma    2.3.3, the direct sum of two
full \wtf\  \amm s $\ph, \ps : A \to \Kt D$ is again
full and \wtf.

\vspace{0.6\baselineskip}

{\em Proof:} Since \wue\  respects direct sums, the previous lemma
implies we may assume $\ph = \ps.$  We may further assume that $\ph$
actually has a \tfn\  $\ph' : \OT{A} \to \Kt D$.
Then $\ph \oplus \ps$ has the standard factorization
$\id_{M_2} \otimes \ph'$ through $({M_2} \otimes \OA{2}) \otimes A,$
and this is a \tfn\  because $M_2 \otimes \OA{2} \cong \OA{2}.$

Fullness follows from Lemma  1.2.6 (2). \QED

\vspace{0.6\baselineskip}

We also need asymptotically standard factorizations through
$\OIA{A}.$ The special properties required in the following lemma will
be used in the proof of Theorem  2.3.7.

\vspace{0.6\baselineskip}

{\bf  2.3.5 Lemma.}
Let $A$ be a separable unital nuclear \ca,
let $D_0$ be unital, and let $D = \OIA{D_0}.$
Let  $\ph : A \to \Kt {D}$ be an \amm. Then $\ph$ has an
asymptotic standard factorization through $\OIA{A}.$ In fact, $\ph$ is
\wyue\  to an \amm\  of the form
$\ps_{t} (a) =
       \dt \circ (\id_{\OI} \otimes \tilde{\ph}_{t}) (1 \otimes a),$
in which $\dt : \OIA{\Kt {D}} \to \Kt {D}$ is an isomorphism,
$\tilde{\ph}$ is completely positive contractive and
asymptotically equal
to $\ph,$ and $\id_{\OI} \otimes \tilde{\ph}_{t}$
is defined to be the tensor product of completely positive maps and
is again completely positive contractive.

\vspace{0.6\baselineskip}

{\em Proof:}
Lemma  1.1.5 provides a completely positive contractive
\amm\  $\tilde{\ph}$ which is \sye\  to $\ph.$ Then
$\id_{\OI} \otimes \tilde{\ph}_{t}$ is the minimal tensor
product of two completely positive contractive linear maps, and is
therefore bounded and completely positive by Proposition
IV.4.23 (i) of \cite{Tk}. Looking at the proof of that proposition
and of Theorem IV.3.6 of \cite{Tk}, we see that such a tensor
product is in fact contractive. Thus,
$\| \id_{\OI} \otimes \tilde{\ph}_{t} \| \leq 1$
for all $t.$ One checks that
$t \mapsto ( \id_{\OI} \otimes \tilde{\ph}_{t}) (b)$
is continuous for $b$ in the algebraic tensor product of $\OI$ and
$A.$ It follows that continuity holds for all $b \in \OIA{A}.$
Similarly, one checks that
$t \mapsto \id_{\OI} \otimes \tilde{\ph}_{t}$ is
asymptotically multiplicative, so is an \amm.

Use Corollary    2.2.8 to find an isomorphism
$\dt_0 : \OIA{D} \to D$ such that $d \mapsto \dt_0 (1 \otimes d)$ is
\wyue\  to $\id_D.$ This induces an isomorphism
$\dt : \OIA{ \Kt D} \to \Kt D,$ and a unitary path
$t \mapsto u^{(0)}_t \in M (\Kt D)$ such that
$\| u^{(0)}_t \dt (1 \otimes d) (u^{(0)}_t)^* - d \| \to 0$
for all $d \in \Kt D.$
By Lemma 1.3.9, there is a unitary path $t \mapsto u_t \in (\Kt D)^+$
such that $\| u_t \dt (1 \otimes d) u_t^* - d \| \to 0$
for all $d \in \Kt D.$

We prove that the $\ps$ that results from
these choices is in fact \wyue\  to $\tilde{\ph};$ this will prove the
lemma. Choose finite subsets $F_1 \subset F_2 \subset \cdots$ whose
union is dense in $A.$
For each $n,$ note that the set
$S_n = \{ \tilde{\ph}_{t} (a) : a \in F_n,\, t \in [0, n] \}$
is compact in $D,$ so that there is $r_n$ with
$\| u_t \dt (1 \otimes d) u_t^* - d \| < 2^{-n}$ for all $d \in S_n$
and $t \geq r_n.$
For $\af \in [0, 1]$ define
$f (n + \af) = (1 - \af) r_n + \af r_{n + 1}.$ Then define
unitaries $v_{t} \in (\Kt D)^+$ by $v_{t} = u_{f (t)}.$
For $t \in [n, n + 1]$ and $a\in F_n,$ this gives
(using $f (t) \geq r_n$)
\[
\|v_{t} \ps_{t} (a) v_{t}^* - \tilde{\ph}_{t} (a) \| =
   \|u_{f (t)} \dt (1 \otimes \tilde{\ph}_{t} (a)) u_{f (t)}^* -
                                 \tilde{\ph}_{t} (a) \| < 2^{-n}.
\]
Thus $\ps$ is in fact \wyue\  to $\tilde{\ph}.$ \QED

\vspace{0.6\baselineskip}

{\bf  2.3.6 Lemma.} (Compare \cite{Ph2}, Proposition 3.3.) Assume the
hypotheses of Lemma    2.3.3. Let
$\ph, \ps : A \to \Kt D$ be full   \amm s with
$\ps$ \wtf. Then $\ph \oplus \ps$ is \wyue\  to $\ph.$

\vspace{0.6\baselineskip}

{\em Proof:}
By Lemma  2.3.5, we may assume that $\ph$ has a standard
factorization through $\OIA{A},$ say
$\ph' : \OIA{A} \to \Kt D.$
Using Lemma  1.3.8 on $\ph'$ and on an
asymptotically trivializing factorization for $\ps,$ we may assume
without loss of generality that there are \pj s $p, q \in \Kt D$ such
that $\ph'$ is a unital \amm\  from $A$ to $p(\Kt D)p$
and $\ps$ is an \wtf\  unital \amm\  from $A$ to $q(\Kt D)q.$

Choose a nonzero
\pj\  $e \in \OI$ with trivial $K_0$ class.
Let $t \mapsto f_t$ be a tail \pj\  for $\ph' (e \otimes 1).$
Choose a \ct\   unitary family $t \mapsto u_{t}$ in
$p(\Kt D)p$ such that $u_{t} f_{t} u_{t}^* = f_{0}$ for all $t.$
Define bounded \amm s
\[
\sm : (1 - e) \OI (1 - e) \otimes A \to
              (p - f_{0}) (\Kt D) (p - f_{0})
\andeqn
\ta : e \OI e  \otimes A \to f_{0} (\Kt D) f_{0}
\]
by
\[
\sm_{t} (x) = u_{t} (p - f_{t}) \ph'_{t} (x)
             (p - f_{t}) u_{t}^*
\andeqn
\ta_{t} (x) = u_{t} f_{t} \ph'_{t} (x) f_{t} u_{t}^*.
\]
These are in fact \amm s, because
$\lim_{t \to \infty} \| f_{t} - \ph'_{t} (e \otimes 1) \| = 0.$
Then define \amm s
\[
\tilde{\sm} : A \to (p - f_{0}) (\Kt D) (p - f_{0})
\andeqn
\tilde{\ta} : A \to f_{0} (\Kt D) f_{0}
\]
by
\[
\tilde{\sm}_{t} (a) = \sm_{t} ((1 - e) \otimes a))
       \andeqn
\tilde{\ta}_{t} (a) = \ta_{t} (e \otimes a).
\]

It follows that
\[
\lim_{t \to \infty}
           \| u_{t} \ph'_{t} (1 \otimes a) u_{t}^*
                 - \tilde{\sm}_{t} (a) - \tilde{\ta}_{t} (a)\| = 0
\]
for all $a \in A,$ so $\ph$ is \wyue\  to
$\tilde{\sm} \oplus \tilde{\ta}.$ Since $[e] = 0$ in $K_0 (\OI),$
there is a unital \hm\  $\nu : \OA{2} \to  e \OI e,$ and the
formula
$\tilde{\ta}_{t} (a)
        = (\ta_{t} \circ (\nu \otimes \id_A)) (1 \otimes a)$
shows that $\tilde{\ta}$ has a \tfn.
Furthermore, $\tilde{\ta}$ is full   because $\ph'$ is.
So $\tilde{\ta} \oplus \ps$ is full   and \wtf\  by Corollary
 2.3.4, and therefore \wyue\  to $\tilde{\ta}$ by Lemma    2.3.3.
The \wue\  of $\ph$ and $\tilde{\sm} \oplus \tilde{\ta}$ now implies
that $\ph \oplus \ps$ is \wyue\  to $\ph.$ \QED

\vspace{0.6\baselineskip}

We now come to the main technical theorem of this section.

\vspace{0.6\baselineskip}

{\bf  2.3.7 Theorem.} Let $A$ be  \snus.  Let $D_0$ be a unital
\ca, and let $D = \OIA{D_0}.$ Then two full \amm s from $A$ to $\Kt D$
are \wyue\  if and only if they are homotopic.

\vspace{0.6\baselineskip}

This result is a continuous analog of Theorem 3.4 of \cite{Ph2},
which gives a similar result for \aue. In the proof of that theorem,
to get \aue\  to within $\ep$ on a finite set $F,$ it was
necessary to approximately absorb a large direct sum of
\wtf\  \hm s---a direct sum which had to be larger for smaller $\ep$ and
larger $F.$ In the proof of the theorem stated here, we must
continuously interpolate between approximate absorption of
ever larger numbers of \amm s. The resulting argument is rather messy.
We try to make it easier to follow by isolating two pieces as lemmas.

\vspace{0.6\baselineskip}

{\bf  2.3.8 Lemma.} Let $A$ and $D$ be \ca s, with $A$ separable.
Let
$\af \mapsto \ph^{(\af)}$ be a bounded homotopy of \amm s from $A$ to
$D.$ Then there exists a \ct\  function
$f : [0, \infty) \to (0, \infty)$ such that for every
$a \in A,$ we have
\[
\lim_{t \to \infty} \left( \sup_{ | \af_1 - \af_2 | \leq 1/f(t) }
   \| \ph_{t}^{(\af_1)} (a) - \ph_{t}^{(\af_2)} (a) \| \right)
                         = 0.
\]

\vspace{0.6\baselineskip}

{\em Proof:} Choose finite sets
$F_0 \subset F_1 \subset \cdots \subset A$ whose union is dense in $A.$

For each $n$ and each fixed $a \in A,$ the map
$(t, \af) \mapsto \ph^{(\af)}_{t} (a)$ is uniformly \ct\  on
$[0, n] \times [0, 1].$ So there is $\dt_n > 0$ such that
\[
\sup \{ \| \ph_{t}^{(\af_1)} (a) - \ph_{t}^{(\af_2)} (a) \| :
  t \in [0, n], \, | \af_1 - \af_2 | \leq \dt_n, \, a \in F_n \}
  < 2^{-n}.
\]
We may clearly assume $\dt_1 \geq \dt_2 \geq \cdots .$
Let $t \mapsto \dt (t)$ be a \ct\  function such that
$0 < \dt (t) \leq \dt_n$ for $t \in [n - 1, n].$

We claim that if $a \in \bigcup_{n = 0}^{\infty} F_n,$ then
\[
\lim_{t \to \infty} \left( \sup_{ | \af_1 - \af_2 | \leq \dt(t) }
   \| \ph_{t}^{(\af_1)} (a) - \ph_{t}^{(\af_2)} (a) \| \right)
                                   = 0.
\]
To see this, let $a \in F_m.$ For $n \geq m + 1,$ $t \in [n - 1, n],$
and $| \af_1 - \af_2 | \leq \dt(t),$ we have in particular
$| \af_1 - \af_2 | \leq \dt_n,$ so that
$\| \ph_{t}^{(\af_1)} (a) - \ph_{t}^{(\af_2)} (a) \| \leq 2^{-n}.$

The statement of the lemma, using $f (t) = 1 / \dt (t),$
follows from the claim by a standard
argument, since $\ph$ is bounded and $\bigcup_{n = 0}^{\infty} F_n$ is
dense in $A.$ \QED

\vspace{0.6\baselineskip}

{\bf  2.3.9 Lemma.} Let $A$ and $Q$ be \ca s, with $Q$ unital and
nuclear. Let $N \geq 2,$ let $e_0, e_1, \dots, e_N \in Q$ be
\mops\  which sum to $1,$ and let $w \in Q$ be a unitary such that
$w e_0 w^* \leq e_1,$ $w e_j w^* \leq e_j + e_{j + 1}$ for
$1 \leq j \leq N - 1,$ and $w e_N w^* \leq e_N + e_0.$ Let
$a_0, \dots, a_N, b_0, \dots, b_N \in A.$ Then in $Q \otimes A$ we
have
\beqr
\lefteqn{ \left\| (w \otimes 1)
        \left( \sum_{j = 0}^N e_j \otimes a_j \right)
                     (w \otimes 1)^* -
            \sum_{j = 0}^N e_j \otimes b_j \right\|} \\
 & \leq &
 \max \{ \| a_N - b_0 \|, \|a_0 - b_1 \| + \| a_1 - b_1 \|,
   \|a_1 - b_2 \| + \| a_2 - b_2 \|,                \\
 & & \hfill  \dots,
              \|a_{N - 1} - b_N \| + \| a_N - b_N \| \}.
\eeqr

\vspace{0.6\baselineskip}

{\em Proof:} Let
\[
x = (w \otimes 1) \left( \sum_{j = 0}^N e_j \otimes a_j \right)
             (w \otimes 1)^* = \sum_{j = 0}^N w e_j w^* \otimes a_j
      \andeqn
y = \sum_{j = 0}^N e_j \otimes b_j.
\]
Observe that if we take the indices mod $N + 1,$ then
$e_k$ is orthogonal to $w e_j w^*$
whenever $k \neq j, j + 1,$ and also if $j = k = 0.$ Therefore we
can calculate
\beqr
x - y & = & \left( \sum_{i = 0}^N e_i \otimes 1 \right) (x - y)
               \left( \sum_{k= 0}^N e_k \otimes 1 \right) \\
      & = & \sum_{j = 0}^N \left[ \rule{0em}{3ex}
                       (e_j \otimes 1)(w e_j w^* \otimes a_j
                    + w e_{j - 1} w^* \otimes a_{j - 1})(e_j \otimes 1)
                                   - e_j \otimes b_j \right]    \\
      & & \mbox{} + \sum_{j = 0}^N \left[ \rule{0em}{3ex}
               e_j (w e_j w^*) e_{j + 1}
              + e_{j + 1} (w e_j w^*) e_j\right] \otimes b_j.
\eeqr

We now claim that the second term in the last expression is zero.
The \pj s $w e_k w^*$ are orthogonal and add up to $1,$ and $e_{j + 1}$
is orthogonal to all of them except for $k = j$ and $k = j + 1.$
Therefore $e_{j + 1} \leq w e_j w^* + w e_{j + 1} w^*.$ Also,
$ e_j w e_{j + 1} w^* = 0,$  so we obtain
\[
e_j (w e_j w^*) e_{j + 1} =
  e_j (w e_j w^* + w e_{j + 1} w^*) e_{j + 1} = e_j e_{j + 1} = 0.
\]
Similarly, $e_{j + 1} (w e_j w^*) e_j = 0.$ So the claim is proved.

It remains to estimate the first term. Since the summands are
orthogonal, the norm of this term is bounded by the maximum of the norms
of the summands. Using again $e_j \leq w e_{j - 1} w^* + w e_j w^*,$
we obtain
\beqr
\lefteqn{\| (e_j \otimes 1)(w e_j w^* \otimes a_j
                    + w e_{j - 1} w^* \otimes a_{j - 1})(e_j \otimes 1)
                                   - e_j \otimes b_j \| }   \\
 & \leq & \| a_{j - 1} - a_j \| +
   \| (e_j \otimes 1)(w e_j w^* \otimes a_j
                    + w e_{j - 1} w^* \otimes a_{j})(e_j \otimes 1)
                                   - e_j \otimes b_j \|    \\
 & = & \| a_{j - 1} - a_j \| + \| e_j \otimes (a_j - b_j) \|
  \leq \| a_{j - 1} - a_j \| + \| a_j - b_j \|.
\eeqr
If $j = 0,$ then $j - 1 = N.$ We then have also $e_0 w e_0 w^* = 0,$
so $e_0 \leq w e_N w^*,$ whence
\[
\| (e_0 \otimes 1)(w e_0 w^* \otimes a_0
    + w e_N w^* \otimes a_N)(e_0 \otimes 1) - e_0 \otimes b_0 \|
  =  \| e_0 \otimes (a_N - b_0) \| \leq \| a_N - b_0 \|.
\]
This proves the lemma. \QED

\vspace{0.6\baselineskip}

{\em Proof of Theorem  2.3.7:}
That \wue\  implies homotopy is Lemma    1.3.3 (1). We therefore prove
the reverse implication.

Using Lemma  2.3.5, we may \wolog\  assume our homotopy has the form
$\tilde{\ph}^{(\af)}_{t} (a) =
                \dt (1_{\OI} \otimes \ph^{(\af)}_{t} (a)),$
where $\dt : \OIA{\Kt D} \to \Kt D$ is a \hm\  and $\ph$ is
a completely positive contractive \amm\  from
$A$ to $C ([0, 1], \Kt D).$
It then suffices to prove the theorem for the homotopy of \amm s
from $A$ to $\OIA{\Kt D}$ given by
$\overline{\ph}^{(\af)}_{t} (a) = 1 \otimes \ph^{(\af)}_{t} (a).$
(We get an \wue\  of $\tilde{\ph}^{(0)}$ and
$\tilde{\ph}^{(1)}$ by applying $\dt.$)

The next step is to do some constructions in $\OI$ and $\OA{2}.$
Choose a \pj\  $e \in \OI$ with $e \neq 1$ and $[e] = [1]$ in
$K_0 (\OI).$  Choose a unital
\hm\  $\gm : \OA{2} \to (1 - e) \OI (1 - e).$
Define isometries $\tilde{s_j} \in \OI$ by
$\tilde{s}_j = \gm (s_j).$
Let $\ld : \OA{2} \to \OA{2}$ be the standard shift
$\ld (c) = s_1 c s_1^* + s_2 c s_2^*.$
Since any two unital endomorphisms of $\OA{2}$ are
homotopic (by Remark 2.1.2    (1) and the connectedness of the
unitary group of $\OA{2}$),
there is a homotopy $\af \mapsto \om_{\af}$ of
endomorphisms of $\OA{2}$ with $\om_0 = \id_{\OA{2}}$ and
$\om_1 = \ld.$

We will now suppose that we are given \ct\  functions
$\af_{n} : [n - 1, \infty) \to [0, 1]$ for $n \geq 1$
such that
\beq
\af_{n + 1} (n) = \af_{n} (n) \label{1}
\eeq
for all $n,$ and a \ct\  function $F : [0, \infty) \to (0, \infty).$
(These will be chosen below.) Then we define
$\ps_t : \OT{A} \to  \OIA{\Kt D}$ by
\begin{eqnarray}
\ps_{t} (c \otimes a) & = &
  \sum_{k = 1}^n
        \tilde{s}_2^{k - 1} \tilde{s}_1 \gm (c) (\tilde{s}_2^{k - 1}
             \tilde{s}_1)^* \otimes \ph_{t}^{(\af_{k} \circ F(t))} (a)
\nonumber \\
& & \mbox{}
       + \tilde{s}_2^n \gm (\om_{F (t) - n} (c)) ( \tilde{s}_2^n)^*
                       \otimes \ph_{t}^{(\af_{n + 1} \circ F(t))} (a)
\nonumber \\
  & &       \,\,\,\,\,\,\,\,\,\,\,\,\,\,\,\,\,\,\,\,\,\,\,\,
             {\rm for} \,\,\, F(t) \in [n, n + 1].
        \label{2}
\end{eqnarray}
(This is an orthogonal sum since the \pj s
\[
\tilde{s}_1 \tilde{s}_1^*, \,\,
\tilde{s}_2 \tilde{s}_1 (\tilde{s}_2 \tilde{s}_1)^*, \dots , \,\,
\tilde{s}_2^{n - 1} \tilde{s}_1 (\tilde{s}_2^{n - 1} \tilde{s}_1)^*,\,\,
\tilde{s}_{2}^n (\tilde{s}_2^n )^*
\]
are mutually orthogonal.)
As in the proof of Lemma  2.3.5, each $\ps_{t}$ is well
defined,  linear, and contractive, and $t \mapsto \ps_{t} (b)$
is continuous for $b$ in the algebraic tensor product of $\OA{2}$ and
$A$ (using (\ref{1}) when $F (t) \in {\bf N}$), and so for all
$b \in \OT{A}.$

We now claim that $\ps,$ as defined by (\ref{2}), is actually an
\amm\  from $\OT{A}$ to $\OIA{\Kt D}.$ It only
remains to prove asymptotic multiplicativity. By linearity and
finiteness of $\sup_{t} \| \ps_{t} \|,$ it suffices to do this
on elementary tensors. Since $\gm,$ $\om_{\af},$ and the maps
$c \mapsto
\tilde{s}_2^{k - 1} \tilde{s}_1 c (\tilde{s}_2^{k - 1} \tilde{s}_1)^*$
and $c \mapsto \tilde{s}_{2}^n c (\tilde{s}_2^n )^*$
are \hm s (and so contractive), a calculation gives,
for $F(t) \in [n, n + 1],$
\beqr
\lim_{t \to \infty}
\lefteqn{ \| \ps_{t} ((c_1 \otimes a_1) (c_2 \otimes a_2)) -
   \ps_{t} (c_1 \otimes a_1) \ps_{t} (c_2 \otimes a_2) \| }  \\
 & \leq & \lim_{t \to \infty} \| c_1 c_2 \|
      \left( \sup_{\af \in [0, 1]}
        \left\| \ph_{t}^{(\af)} (a_1 a_2) -
           \ph_{t}^{(\af)} (a_1) \ph_{t}^{(\af)} (a_2) \right\|
                                     \right) = 0.
\eeqr

Define $\io : A \to \OT{A}$ by $\io (a) = 1 \otimes a.$ Then
$\ps \circ \io$ is an \amm\  from $A$ to
$\OIA{\Kt D}.$ By definition, it has a \tfn, so
Lemma    2.3.6 implies that
$\overline{\ph}^{(\af)} \oplus (\ps \circ \io)$ is \wyue\  to
$\overline{\ph}^{(\af)}.$ The theorem will therefore be proved if we
can choose the functions $F$ and $\af_{n}$ in such a way that
$\overline{\ph}^{(0)} \oplus (\ps \circ \io)$ is \wyue\  to
$\overline{\ph}^{(1)} \oplus (\ps \circ \io).$

Before actually choosing $F$ and the $\af_{n}$, we construct, in
terms of $F,$ the unitary path we will use for the desired \wue.
Let $\ta$ be an automorphism of $M_2 (\OI)$ which sends
$1 \oplus 0$ to $e \oplus 0$ and $0 \oplus c$
to $c \oplus 0$ for all $c \in (1 - e) \OI (1 - e).$
Let $\tilde{\ta}$ be the obvious induced automorphism of
$M_2 (\OIA{\Kt D}).$ It suffices to prove \wue\  of
$\tilde{\ta} \circ (\overline{\ph}^{(0)} \oplus (\ps \circ \io))$ and
$\tilde{\ta} \circ (\overline{\ph}^{(1)} \oplus (\ps \circ \io)).$
Furthermore, these two \amm s take values in
$\OIA{\Kt D},$ embedded as the upper left corner, so we
only work there. This results in the identification
\[
\tilde{\ta} \circ (\overline{\ph}^{(\af)} \oplus (\ps \circ \io))
 = e \otimes \ph^{(\af)} (-) + (\ps \circ \io).
\]
We further note that, by
Lemma  1.3.9, it suffices to construct a \ct\  family of
unitaries in the multiplier algebra $M( \OIA{\Kt D}).$

With these identifications and reductions, our unitary path will take
the form $u_{t} = v (F(t)) \otimes 1$ for a suitable unitary path
$r \mapsto v(r)$ in $\OI,$ defined for $r \in [0, \infty).$
The construction of $v$ requires further notation.

Define projections in $\OI$ by
$p_k =
\tilde{s}_2^{k - 1} \tilde{s}_1 (\tilde{s}_2^{k - 1} \tilde{s}_1)^*$
and $q_n = \tilde{s}_{2}^n (\tilde{s}_2^n )^*.$
Then
\[
p_1 + p_2 + \cdots + p_n + q_n + e = 1 \andeqn
          p_{n + 1} + q_{n + 1} = q_{n}.
\]
Choose projections $f_k < p_k$ with $[f_k] = 1$ in $K_0 (\OI).$
Note that $f_{n + 1} < q_n.$
Then there are partial isometries $v_k$ with
\[
v_0^* v_0 = e, \,\,\,\,\,\, v_0 v_0^* = f_1, \,\,\,\,\,\,
v_k^* v_k = f_k, \andeqn      v_k v_k^* = f_{k + 1},
\]
and $w_n$ with
\[
w_n^* w_n = f_{n + 1} \andeqn  w_n w_n^* = e.
\]
Using the connectedness of the unitary group of
$(f_{n + 1} + f_{n + 2}) \OI (f_{n + 1} + f_{n + 2}),$ choose a
\ct\  path $\af \mapsto y_n (\af)$ of partial isometries from
$f_{n + 1} + f_{n + 2}$ to $f_{n + 2} + e$ such that
$y_n (0) = w_n + f_{n + 2}$ and $y_n (1) = v_{n + 1} + w_{n + 1}.$
Then define
\[
v (n + \af) = (p_1 - f_1) + \cdots + (p_{n + 1} - f_{n + 1}) +
   (q_{n + 1} - f_{n + 2}) + v_0 + v_1 + \cdots + v_n + y_n (\af)
\]
for $n \in {\bf N}$ and $\af \in [0, 1].$ There are two definitions at
each integer, but they agree, so $v$ is a \ct\  path of unitaries.
Furthermore, one immediately verifies that for fixed $r \in [n, n + 1],$
the unitary $w = v(r)$ and sequence of \pj s
\beq
e_0 = e, \,\, e_1 = p_1, \,\, e_2 = p_2, \dots, e_n = p_n, \,\,
   e_{n + 1} = q_n      \label{3}
\eeq
satisfy the hypotheses in Lemma  2.3.9.

Now take $f$ to be as in Lemma    2.3.8, and set $F (t) = f (t) + 2.$
Define $\af_{0} : [0, \infty ) \to [0, 1]$ by $\af_0 (r) = 0$ for all
$r,$ and choose the functions $\af_{n} : [n - 1, \infty ) \to [0, 1]$
to be continuous, to satisfy (\ref{1}), and such that
$\af_{n + 1} (r) = 1$ for $r \in [n, n + 1]$ and
\[
| \af_{k + 1} (r) - \af_{k} (r)| \leq 1 / (n - 1) \,\,\,\,\,\,
{\rm for}\,\,\, r \in [n, n + 1] \,\,\, {\rm and}\,\,\, 0 \leq k \leq n.
\]

Take $\ps$ and $u$ to be defined using these choices of $F$ and the
$\af_{n}.$
Let $t \in [0, \infty)$. Set $r = F(t)$
and choose $n \in {\bf N}$ such that $r \in [n, n + 1].$ Let $w = v(r)$
and let $e_0, \dots, e_{n + 1}$ be as in (\ref{3}).
For $a \in A,$ we then have
\beqr
\lefteqn{
   \left\| u_{t} \left[ e \otimes \ph^{(0)} (a)
                                 + \ps (1 \otimes a) \right] u_{t}^* -
   \left[ e \otimes \ph^{(1)} (a) + \ps (1 \otimes a) \right] \right\| }
                          \\
 & = & \left\| (w \otimes 1) \left[
        \sum_{k = 0}^{n + 1} e_k \otimes \ph_{t}^{(\af_{k} (r))} (a)
                               \right] (w \otimes 1)^*
   - \left[ e_0 \otimes \ph_{t}^{(1)} (a)
      + \sum_{k = 1}^{n + 1} e_k \otimes \ph_{t}^{(\af_{k} (r))} (a)
                                                   \right] \right\|.
\eeqr
Apply Lemma    2.3.9 with
$a_k = b_k = \ph_{t}^{(\af_{k} (r))} (a)$
for $1 \leq k \leq n + 1,$ and with $a_0 = \ph_{t}^{(\af_0 (r))} (a)$
and $b_0 = \ph_{t}^{(1)} (a) = \ph_{t}^{(\af_{n + 1} (r))} (a)
                               = a_{n + 1}.$
It follows that the expression above is at most
\beqr
\max (0, \| a_0 - a_1 \|, \dots, \| a_n - a_{n + 1} \| )
   & = & \max \{ \|  \ph_{t}^{(\af_{k} (r))} (a) -
           \ph_{t}^{(\af_{k + 1} (r))} (a) \| : 0 \leq k \leq n \} \\
   & \leq & \sup \{ \|  \ph_{t}^{(\af_1)} - \ph_{t}^{(\af_2)} \| :
                               | \af_1 - \af_2 | \leq 1/(n - 1) \}.
\eeqr
Since $n - 1 \geq r - 2 = f (t),$ we have
$1/(n - 1) \leq 1/f(t),$ and this last expression converges to $0$
as $t \to \infty.$ Thus we have shown that
\[
e \otimes \ph^{(0)} (-) + (\ps \circ \io)   \andeqn
       e \otimes \ph^{(1)} (-) + (\ps \circ \io)
\]
are \wyue. This completes the proof. \QED

\vspace{0.6\baselineskip}

{\bf  2.3.10 Corollary.}
Let $A$ be \snus, let $D_0$ be unital, and
let $D = \OIA{D_0}.$  Then any full
\amm\  $\ph : A \to \Kt D$ is \wyue\  to a \hm.

\vspace{0.6\baselineskip}

{\em Proof:}
It is obvious that an \amm\  is homotopic to all of its
reparametrizations.
The result therefore follows from Theorem 2.3.7 and Proposition 1.3.7.
\QED

\vspace{0.6\baselineskip}

{\bf 2.3.11 Remark.} The hypothesis of fullness can be removed in
Theorem  2.3.7 (and in Corollary 2.3.10) in the following way.
Let $\af \mapsto \ph^{(\af)}$  be a homotopy of
\amm s from $A$ to $\Kt D,$ with $D = \OIA{D_0}.$ Applying
Lemma  1.3.8, we can assume $\af \mapsto \ph^{(\af)}$ is a homotopy of
unital (hence full) \amm s from $A$ to $D' = p (\Kt D) p$ for
a suitable \pj\  $p.$ The algebra $D'$ is stable under tensoring with
$\OI$ by Corollary 2.1.12. So we can apply the result already
proved to \amm s from
$A$ to $\Kt D'$. Then embed  $\Kt D'$ in  $\Kt D.$

\vspace{\baselineskip}
\section{Unsuspended $E$-theory for simple nuclear \ca s}
\vspace{\baselineskip}

In \cite{DL},  D\v{a}d\v{a}rlat and Loring proved that for certain
\ca s $A,$ one can obtain the groups $KK^0 (A, B)$ via
``unsuspended $E$-theory'': $KK^0 (A, B) \cong [[\Kt A, \Kt B]]$
(notation from Definition  1.1.2) for all separable $B.$
The terminology comes from the
omission of the suspension that is normally required. The conditions
on $A$ are quite restrictive, and in particular fail for trivial reasons
as soon as $A$ has even one nonzero \pj.

In this section, we want to take $A$ to be \snus. To make enough room,
we assume $B$ is a tensor product $\OIA{D}$ with $D$ unital. We then
discard the the class of the zero \amm\  (the source of the difficulty
with projections). We are able to prove, with the help of Kirchberg's
results as stated in Section 2.1 and also using Theorem  2.3.7, that
we do in fact get $KK^0 (A, B)$ as a set of suitable homotopy classes
of \amm s from $\Kt A$ to $\Kt B.$ (Corollary 2.3.10 implies that
we can even use \wue\  classes of \hm s. See Section 4.1.)

In the first subsection,
we construct for fixed $A$ a middle exact homotopy invariant
functor from separable \ca s to abelian groups in a manner
analogous to the definition of $K_0 (D),$ but using \amm s from
$A$ to $\KOI{D^+}$ in place  of projections in $\Kt D^+.$
The fact that the target algebra is infinite means that,
as for $K_0$ of a \pisca, we do not need to take formal differences
of classes. We do, however, need to introduce the unitization
of the target algebra for essentially the same reason that it is
necessary in the definition of $K_0.$
In the second subsection, we then show that this functor
is naturally isomorphic to $KK^0(A, -).$

\vspace{0.6\baselineskip}
\subsection{The groups $[[A, \KOI{D}]]_+$ and $\ET_A (D)$}
\vspace{0.6\baselineskip}

Let $A$ be \snus. In this subsection we construct a functor
$[[A, \KOI{-}]]_+$ on unital \ca s and the corresponding functor
$\ET_A (-)$ on general \ca s (obtained via the unitization). We then
prove that $\ET_A$ is a cohomology theory on separable \ca s in the
usual sense. This information is needed in order to apply the
uniqueness theorems for $KK$-theory in the next subsection.

\vspace{0.6\baselineskip}

{\bf  3.1.1 Definition.}
Let $A$ be separable and unital, and assume each ideal of $A$ is
generated by its projections.  Let $B$ have an approximate identity
of \pj s.  Then $[[A, B]]_+$ denotes the
set of homotopy classes of full \amm s from $A$ to $B.$

\vspace{0.6\baselineskip}

{\bf  3.1.2 Proposition.}
Let $A$ be simple, separable, unital, and nuclear.
For any unital \ca\  $D,$ give
$[[A, \KOI{D}]]_+$ the addition operation that it receives from being
a subset of $[[A, \KOI{D}]].$ Then  $[[A, \KOI{-}]]_+$
is a functor from separable unital \ca s and
homotopy classes of unital \amm s to
abelian groups. The zero element is the class of any full   \amm\  from
$A$ to $\KOI{D}$ with a standard factorization (see Definition 2.3.1)
through $\OT{A}.$

\vspace{0.6\baselineskip}

{\em Proof:} Lemma  1.2.6 (2)  shows that $[[A, \KOI{D}]]_+$ is closed
under the addition in $[[A, \KOI{D}]].$ Therefore  $[[A, \KOI{D}]]_+$
is an abelian semigroup, provided it is not empty.

According to Theorem 2.3.7,
homotopy is the same relation as \wue\  in this
set. So we can use them interchangeably.

For functoriality, let $E$ be another unital \ca, and let
$\ph : D \to E$ be a unital \amm. Let
$\overline{\ph} = \id_{\Kt \OI} \otimes \ph$ (see Proposition    1.1.8)
be the induced \amm\  from $\KOI{D}$ to $\KOI{E}.$ It is
full because if $e \in K$ is any nonzero projection, then
$e \otimes 1 \otimes 1$ is a \dpj\  in $\KOI{D}$ which is sent to the
\dpj\  $e \otimes 1 \otimes 1$  in $\KOI{E}.$ Lemmas  1.2.6 (2) and
 2.1.8 (1) now imply that $\et \mapsto [[\ph]] \cdot \et$ sends
full   \amm s to full   \amm s.

We now construct an identity element.
Theorem    2.1.4 provides an isomorphism $\nu : \OT{A} \to \OA{2}.$
Let $\ta : \OA{2} \to \OI$ be an injective \hm\  (sending $1$ to a
nonzero \pj\  in $\OI$ with trivial $K_0$-class), and define a full
\hm\  $\zt : A \to \OI$ by $\zt (a) = (\ta \circ \nu) ( 1 \otimes a).$
Composing it with the
full   \hm\  $x \mapsto e \otimes x  \otimes 1$ from $\OI$ to
$\KOI{D},$ where $e \in K$ is any nonzero projection, we obtain
a full   \amm\  from $A$ to $\KOI{D}$ which has a
standard factorization through $\OT{A}.$

Lemma    2.3.3 implies that any other full    \amm\  with a \tfn\  is
\wyue\  to $\zt.$ This class acts as the identity by Lemma    2.3.6.

Finally, we must construct additive inverses.
Let $\et \in [[A, \KOI{D}]]_+$. By Lemma  2.3.5, we can take
$\et = [[\ph]],$ where $\ph$ has a standard factorization through
$\OIA{A},$ say $\ph_t (a) = \ps_t (1 \otimes a)$ for some
\amm\ $\ps : \OIA{A} \to \KOI{D}.$
Choose a \pj\  $f \in \OI$ with $[f] = -1$
in $K_0 (\OI).$  Define $\overline{\ps}_t =\ps_t |_{f \OI f \otimes A},$
and define $\overline{\ph} : A \to \KOI{D}$ by
$\overline{\ph}_t (a) = \ph_t (f  \otimes a).$ Choose
a unital \hm\  
\[
\nu : \OA{2} \to
\left( \begin{array}{cc} 1 & 0 \\ 0 & f \end{array} \right)
 M_2 (\OI)
\left( \begin{array}{cc} 1 & 0 \\ 0 & f \end{array} \right).
\]
Then $(\id_{M_2} \otimes \ps) \circ \nu$
provides a standard factorization of  $\ph \oplus \overline{\ph}$
through $\OT{A}.$ Note that $\ph \oplus \overline{\ph}$ is
full   because $\ph$ is, so it is \wyue\  to $\zt$ by Lemma    2.3.3.
This shows that $[[\overline{\ph}]]$ is the inverse
of $\et.$ \QED

\vspace{0.6\baselineskip}

{\bf  3.1.3 Definition.} If $D$ is any \ca, then we denote by
$D^{\#}$ the \ca\  $\KOI{D^+}.$ We use the analogous notation for \hm s.
If $D$ is separable, we define
$\ET_A (D)$ to be the kernel of the map
$[[A, D^{\#}]]_+ \to [[A, \Kt \OI]]_+$ induced by the unitization
map $D^+ \to {\bf C}.$

\vspace{0.6\baselineskip}

{\bf  3.1.4 Proposition.}
Let $A$ be \snus.  Then  $\ET_A$
is a functor from separable \ca s and
homotopy classes of \amm s to
abelian groups.

\vspace{0.6\baselineskip}

{\em Proof:} This follows from Proposition  3.1.2 and the fact that
unitizations and tensor products
of \amm s are well defined (Lemma  1.1.6 and Proposition   1.1.8). \QED

\vspace{0.6\baselineskip}

{\bf  3.1.5 Remark.} It is obvious that if $D_1$ and $D_2$ are
unital, then there is a natural isomorphism
\[
[[A, \KOI{ (D_1 \oplus D_2) } ]]_+ \cong
  [[A, \KOI{D_1} ]]_+ \oplus [[A, \KOI{ D_2) } ]]_+.
\]
It follows that for unital $D,$ there is a natural isomorphism
\[
\ET_A (D) \cong [[A, \KOI{D}]]_+.
\]

\vspace{0.6\baselineskip}

We will sometimes denote by $\ph_*$ the map
$[[A, D_1]]_+ \to [[A, D_2]]_+$ or the map $\ET_A (D_1) \to \ET_A (D_2)$
induced by a (full) \hm\ $\ph : D_1 \to D_2.$

\vspace{0.6\baselineskip}

{\bf  3.1.6 Lemma.} Let $A$ \snus. Let
\[
0 \longrightarrow J  \stackrel{\mu}{\longrightarrow} D
         \stackrel{\pi}{\longrightarrow} D/J \longrightarrow 0
\]
be a short exact sequence of separable \ca s.
Then the sequence
\[
\ET_A (J)  \stackrel{\mu_*}{\longrightarrow} \ET_A ( D)
         \stackrel{\pi_*}{\longrightarrow}  \ET_A (D/J)
\]
is exact in the middle.

\vspace{0.6\baselineskip}

{\em Proof:} It is immediate that $\pi_* \circ \mu_* = 0.$

For the other half, we introduce the maps $\ch_D : D^{\#} \to \Kt \OI$
and $\io_D : \Kt \OI \to D^{\#}$ associated with the unitization maps
$D^+ \to {\bf C}$ and ${\bf C} \to D^+.$ Define $\ch_{D/J},$
$\io_{D/J},$ etc.\  similarly. Now let $\et \in \ker (\pi_*),$ and
choose a full    \amm\  $\ph : A \to D^{\#}$ whose class is $\et.$
By definition, we have $[[\pi^{\#} \circ \ph]] = 0$ in
$[[A, (D/J)^{\#}]]_+.$ Choose a full    \hm\  $\zt : A \to \Kt \OI$
with a standard factorization through $\OT{A},$ as in the proof of
Proposition    3.1.2.    Theorem 2.3.7 then implies that
$\pi^{\#} \circ \ph$ is \wyue\  to $\io_{D/J} \circ \zt,$
so there is a unitary path $t \to u_t$ in $((D/J)^{\#})^+$ such that
$u_t (\pi^{\#} \circ \ph_t) (a) u_t^*  \to (\io_{D/J} \circ \zt) (a)$
for all $a \in A.$

Without changing homotopy classes, we may replace $\ph$ by
$\ph \oplus 0$
and $\zt$ by $\zt \oplus 0.$ This also replaces $\pi^{\#} \circ \ph$
and $\io_{D/J} \circ \zt$ by their direct sums with the zero \amm.
We then replace $u_t$ by $u_t \oplus u_t^*$. We may thus assume
\wolog\  that $u$ is in the identity component of the unitary group of
$\Cb ([0, \infty), ((D/J)^{\#})^+).$
Therefore there is $v \in U_0 (\Cb ([0, \infty), (D^{\#})^+))$
whose image is $u.$ Then $\pi^{\#} (v_t) = u_t$ for all
$t,$ whence
\[
\lim_{t \to 0} \pi^{\#} (v_t \ph_t (a) v_t^* -
        (\io_{D} \circ \zt) (a)) = 0
\]
for all $a \in A.$

Let $\sm : (D/J)^{\#} \to D^{\#}$ be a continuous (nonlinear) cross
section for $\pi^{\#}$ satisfying $\sm (0) = 0.$ (See \cite{BG}.) Define
$\ps_t : A \to D^{\#}$ by
\[
\ps_t (a) = v_t \ph_t (a) v_t^* -
 (\sm \circ \pi^{\#})
    \left( \rule{0em}{3ex}
                 v_t \ph_t (a) v_t^* - (\io_{D} \circ \zt) (a) \right).
\]
This yields an \amm\  asymptotically equal to
$t \mapsto v_t \ph_t (-) v_t^*,$
and hence \wyue\  to $\ph.$ Furthermore,
$\pi^{\#} (\ps_t (a) - (\io_{D} \circ \zt) (a)) = 0$ for all
$t$ and $a.$ It follows that $\ps_t (a) \in J^{\#}$ and that
$\ch_J (\ps_t (a)) = \zt (a).$ So $\ps$ is in fact
a full    \amm\  from $A$ to $J^{\#}$ such that
$[[\ch_J \circ \ps]] = 0,$ from which it follows that $\ps$ defines a
class $[[\ps]] \in \ET_A (J).$ Clearly $\mu_* ([[\ps]]) = \et.$
This shows that $\ker(\pi_*) \i {\rm Im}(\mu_*).$  \QED

\vspace{0.6\baselineskip}

{\bf  3.1.7 Proposition.} Let $A$ and
\[
0 \longrightarrow J  \stackrel{\mu}{\longrightarrow} D
         \stackrel{\pi}{\longrightarrow} D/J \longrightarrow 0
\]
be as in Lemma  3.1.6. Then there is a natural exact sequence
\[
 \cdots  \stackrel{(S\mu)_*}{\longrightarrow}
 \ET_A ( SD  )  \stackrel{(S\pi)_*}{\longrightarrow}
 \ET_A (S(D/J)) {\longrightarrow}
\ET_A (J)  \stackrel{\mu_*}{\longrightarrow} \ET_A ( D)
         \stackrel{\pi_*}{\longrightarrow}  \ET_A (D/J).
\]

\vspace{0.6\baselineskip}

{\em Proof:} This follows from middle exactness (the previous lemma)
and homotopy invariance by standard methods. See, for example,
Section 7 of \cite{Ksp}. \QED

\vspace{0.6\baselineskip}

{\bf  3.1.8 Remark.} It should be pointed out that we need much less
than the full strength of    Theorem 2.3.7 here. Only knowing that
homotopy implies \wue\  for full \amm s from $A$ to
$\Kt \OIA{C([0, 1])},$ it is possible to prove middle exactness in the
first stage of the Puppe sequence, namely
\[
\ET_A (C\pi) \longrightarrow \ET_A (D) \longrightarrow  \ET_A (D/J).
\]
This sequence can be extended to the left as in the proof of Proposition
2.6 of \cite{Sc2}. Proposition 3.2 of \cite{DL} can then be used to show
that $\ET_A$ is split exact; this is the property we actually use in the
next section.

\vspace{0.6\baselineskip}

We now prove stability of $\ET_A$ under formation of tensor products
with both $K$ and $\OI.$

\vspace{0.6\baselineskip}

{\bf  3.1.9 Lemma.}
Let $A$ be \snus, and let $D$ be a separable \ca. Then the
map $d \mapsto 1 \otimes d,$ from $D$ to $\OIA{D},$ induces
an isomorphism $\ET_A (D) \to \ET_A (\OIA{D}).$

\vspace{0.6\baselineskip}

{\em Proof:} By naturality, Proposition  3.1.7, and the Five Lemma,
it suffices to prove this for unital $D.$ By Remark  3.1.5, we have
to prove that $d \mapsto 1 \otimes d$ induces an isomorphism
$[[A, \KOI{D}]]_+ \to [[A, \KOI{\OIA{D}}]]_+.$ This follows from
Theorem  2.1.5 and Proposition 2.1.11, since these results imply that
the map $x \mapsto x \otimes 1,$ from $\OI$ to $\OIA{\OI},$ is
homotopic to an isomorphism. \QED

\vspace{0.6\baselineskip}

The other stability result requires the
following lemma. We really want an increasing \ct ly parametrized
approximate identity of \pj s, but of course such a thing does not
exist. The quasiincreasing version in the lemma is good enough.

\vspace{0.6\baselineskip}

{\bf   3.1.10 Lemma.}
Let $D$ be a unital \pisca, and let $e_0 \in \Kt D$ be a nonzero \pj.
Then there exists a \ct\  family
$t \mapsto e_t$ of \pj s in $\Kt D$ such that, for every $b \in \Kt D$,
we have
\[
\lim_{t \to \infty} (e_t b - b) = \lim_{t \to \infty} (b e_t - b)
     = \lim_{t \to \infty} (e_t b e_t - b) = 0,
\]
such that $e_0$ is the given \pj, and such that $e_s \geq e_t$ for
$s \geq t + 1.$

\vspace{0.6\baselineskip}

{\em Proof:} Choose a nonzero \pj\  $p  \in \Kt D$ such that
$[p] = 0$ in $K_0 (D).$ We start by constructing a family
$t \mapsto f_t$ in $\Kt pDp.$ Note that
\[
[\diag (1_{pDp}, 0, 0)] = [\diag (1_{pDp}, 1_{pDp}, 0)] = 0
\]
in $K_0 (M_3 (pDp)).$ Therefore there is a homotopy $t \mapsto q_t$
of \pj s in $M_3 (pDp)$ such that
\[
q_0 = \diag (1, 0, 0) \andeqn q_1 = \diag (1, 1, 0).
\]
Now define
\[
f_{n + s} = 1_{M_{n + 1} (pDp)} \oplus q_s \oplus 0 \in \Kt pDp
\]
for $n = 0, 1, \dots$ and $s \in [0, 1].$ The family $f_t$ is clearly
\ct. It satisfies $f_0 = p \oplus p.$ We have
$f_t \geq 1_{M_{n + 1} (pDp)}$ for
$t \geq n,$ so $t \mapsto f_t$ really is an approximate identity.
Finally, $f_t \leq 1_{M_{n + 3} (pDp)}$ for $t \leq n,$
so $f_s \geq f_t$ for $s \geq t + 4.$ We can replace $4$ by $1$ in
this last statement by a reparametrization.

To get the general case, choose a \pj\  $r \in pDp$ with
$[r] = - [e_0]$ in $K_0 (D).$ Then $f_t \geq p \geq r$ for
all $t,$ so $t \mapsto f_t - r$ is
a \ct ly parametrized approximate identity of \pj s for
$(1 - r) (\Kt pDp) (1 - r).$ (Here $1$ is the identity
of $(\Kt pDp)^+.$)
There is an isomorphism $\ph : \Kt D \to (1 - r) (\Kt pDp) (1 - r),$
and since $[f_0 - r] = [e_0]$ in $K_0 (D),$ we can require
that $\ph (e_0) =f_0 - r.$ Now set $e_t = \ph^{-1} (f_t - r).$
Then clearly $e_t b - b, \, b e_t - b \to 0$ as $t \to \infty.$
It follows that
\[
\| e_t b e_t - b \| \leq \|e_t b - b \| \| e_t \| + \| b e_t - b \|
         \to 0
\]
as well. \QED

\vspace{0.6\baselineskip}

{\bf 3.1.11 Lemma.} Let $A$ be \snus, let $D$ be
separable, and let $e \in K$ be a rank one \pj. Then the
map $d \mapsto e \otimes d,$ from $D$ to $\Kt D,$ induces
an isomorphism $\ET_A (D) \to \ET_A (\Kt D).$

\vspace{0.6\baselineskip}

{\em Proof:} By Lemma 3.1.9, we may use $\OIA{D}$ in place of $D,$ and
as in its proof we may assume $D$ is unital.

Let $s \in \OI$ be a proper isometry, and define
$\gm : \OIA{D} \to \OIA{D}$ by
$\gm (a) = (s \otimes 1) a (s \otimes 1)^*.$ We claim that
$\gm_* : \ET_A (\OIA{D}) \to \ET_A (\OIA{D})$ is an isomorphism.
It follows from Remark  3.1.5 and Definition  3.1.3 that this map
can be thought of as composition with $\id_{\Kt \OI} \otimes \gm$
from $[[A, \KOI{\OIA{D}}]]_+$ to itself, even though $\gm$ is not
unital. (The discrepancy is an orthogonal sum with an \amm\  which
up to homotopy has a \tfn. Note that the composition with $\gm$
is still full.) Now $\Kt s s^* \OI s s^*$ and
$\Kt (1 - s s^*) \OI (1 - s s^*)$ are both isomorphic to $\Kt \OI,$
so we may as well consider the map from $[[A, \KOI{\OIA{D}}]]_+$ to
$[[A, M_2 ( \KOI{\OIA{D}})]]_+$ induced by inclusion in the upper right
corner. Let $\ta : M_2 (K) \to K$ be an isomorphism. Then
$a \mapsto \ta (a \oplus 0)$ is homotopic to $\id_K$ and
$b \mapsto \ta (b) \oplus 0$ is homotopic to $\id_{M_2 (K)}.$
So our map has an inverse given by composition with
$\ta \otimes \id_{\OIA{\OIA{D}}}.$

We next require a construction involving $\OI$ and $\Kt \OI.$
Define $\ph : \OI \to \Kt \OI$ by $\ph (x) =  e \otimes x.$ Let
$t \mapsto e_t$ be a \ct ly parametrized approximate identity for
$\Kt \OI$ which satisfies the properties of the previous lemma and has
$e_0 = e \otimes 1.$ Let $t \mapsto u_t$ be a \ct\  family of unitaries
in $(\Kt \OI)^+$ such that $u_0 = 1$ and $u_t e_t u_t^* = e_0$ for all
$t.$ Define $\ps^{(0)}_t : \Kt \OI \to \Kt \OI$ by
$\ps^{(0)}_t (a) = u_t e_t a e_t u_t^*.$ One immediately checks that
$\ps^{(0)}$ is an \amm\ whose values are in
$(e \otimes 1) (\Kt \OI) (e \otimes 1),$ so that there is an
\amm\  $t \mapsto \ps_t$ from $\Kt \OI$ to $\OI$ such that
$\ph \circ \ps_t = \ps^{(0)}_t$ for all $t.$

The composite \amm s $\ph \circ \ps$ and $\ps \circ \ph$ can be
computed without reparametrization, because $\ph$ is a \hm.
Now $\ph \circ \ps = \ps^{(0)},$
which is \wyue\  to $(t, a) \mapsto e_t a e_t,$ which in turn is
asymptotically equal to $\id_{\Kt \OI}.$ So
$\ph \circ \ps$ is homotopic to $\id_{\Kt \OI}.$ Also, $\ps \circ \ph$
is clearly homotopic to a map of the form
$x \mapsto s x s^*$ for a proper isometry $s \in \OI.$

We now observe that $\id_{\Kt \OI} \otimes (\ph \otimes \id_D)^+$ and
$\id_{\Kt \OI} \otimes (\ps \otimes \id_D)^+$ define full \amm s from
$(\OIA{D})^{\#}$ to $(\KOI{D})^{\#}$ and back. The composite from
$(\KOI{D})^{\#}$ to itself is homotopic to the identity, and
therefore induces the identity map on $\ET_A (\KOI{D}).$
Composition on the right with the
composite from $(\OIA{D})^{\#}$ to itself is a map of the form $\gm_*$
as considered at the beginning of the proof, and is thus an
isomorphism from $\ET_A (\OIA{D})$ to itself.
It follows that $\ph_*$ is an isomorphism. \QED

\vspace{0.6\baselineskip}
\subsection{The isomorphism with $KK$-theory}
\vspace{0.6\baselineskip}

In this subsection, we prove that if $A$ is \snus,  and
$D$ is separable, then the natural map from
$\ET_A (D)$ to $KK^0 (A, D)$ is an isomorphism. Combined with
Remark  3.1.5, this gives for unital $D$ a form of
``unsuspended $E$-theory'' as in \cite{DL}, in which
we need only discard the zero \amm.

We will use the universal property of $KK$-theory with respect to
split exact, stable, and homotopy invariant functors on separable
\ca s \cite{Hg}. (We use this instead of the related property of
$E$-theory because it is more convenient for the proof of Lemma
3.2.4 below.)

\vspace{0.6\baselineskip}

{\bf  3.2.1 Lemma.} Let $A$  be \snus. Then $\ET_A$ sends split exact
sequences to split exact sequences.

\vspace{0.6\baselineskip}

{\em Proof:} Let
\[
0 \longrightarrow J  \stackrel{\mu}{\longrightarrow} D
         \stackrel{\pi}{\longrightarrow} D/J \longrightarrow 0
\]
be a split short exact sequence of separable \ca s, with
splitting map $\sm : D/J \to D.$
{}From Proposition  3.1.7, we obtain the exact sequence
\[
0 \longrightarrow
\ET_A (J)  \stackrel{\mu_*}{\longrightarrow} \ET_A ( D)
         \stackrel{\pi_*}{\longrightarrow}  \ET_A (D/J).
\]
Using $\sm_*,$ we obtain a splitting; this also shows that the last map
is surjective.
\QED

\vspace{0.6\baselineskip}

{\bf  3.2.2 Notation.} In this subsection,
we denote by $\cal S$ the category of separable \ca s and \hm s
and by ${\cal KK}$ the category  of separable \ca s with
morphisms $KK^0 (A, B)$ for \ca s $A$ and $B.$ If
$\et \in KK^0 (A, B)$ and $\ld \in KK^0 (B, C),$ we denote their
product by $\ld \times \et \in KK^0 (A, C).$ We further denote
by $k$ the functor from $\cal S$ to ${\cal KK}$ which sends a
\hm\  to the class it defines in $KK$-theory.

\vspace{0.6\baselineskip}

{\bf  3.2.3 Corollary.} Let $A$ be \snus.
Then there is a functor $\EH_A$ from ${\cal KK}$ to the
category of abelian groups such that $\EH_A \circ k = \ET_A.$

\vspace{0.6\baselineskip}

This simply means that one can make sense of
$\ET_A (\et) : \ET_A (D) \to \ET_A (F)$  not only when
$\et$ is an \amm\  from $D$ to $F$, but also when
$\et$ is merely an element of $KK^0 (D, F).$

\vspace{0.6\baselineskip}

{\em Proof of Corollary  3.2.3:}
The result is immediate from Theorem 4.5 of \cite{Hg}, since
$\ET_A$ is a stable (Lemma 3.1.11), split exact (Lemma  3.2.1),
and homotopy invariant (Proposition    3.1.4) functor from
separable \ca s to abelian groups. \QED

\vspace{0.6\baselineskip}

We want to show that $\ET_A (D)$ is naturally isomorphic to
$KK^0 ({A}, {D}).$
Our argument is based on an alternate proof of the main theorem
of \cite{DL} suggested by the referee of that paper; we are
grateful to Marius D\v{a}d\v{a}rlat for telling us about it.
The argument requires the construction of certain
natural transformations. (The argument used in Section 4 of \cite{DL}
presumably also works.)

Before starting the construction, we prove a lemma on the functors
$\hat{F}$ of Higson \cite{Hg} (as used in the previous corollary).

\vspace{0.6\baselineskip}

{\bf  3.2.4 Lemma.} Let $F$ and $G$ be stable, split exact,
and homotopy invariant functors from $\cal S$ to the category of
abelian groups, and let $\hat{F}$ and $\hat{G}$ be the unique extensions
to functors from ${\cal KK}$ of Theorem 4.5 of \cite{Hg}. If
$\af$ is a natural transformation from $F$ to $G,$ then
$\af$ is also a natural transformation from $\hat{F}$ to $\hat{G}.$

\vspace{0.6\baselineskip}

{\em Proof:}
Let $\mu \in KK^0 (A, B).$ By Lemma 3.6 of \cite{Hg}, we can choose
a representative cycle (in the sense of Definition 2.1  of \cite{Hg})
of the form $\Ph = (\ph_+, \ph_-, 1),$
where $\ph_+, \ph_- : A \to M(\Kt B)$
are \hm s such that $\ph_+ (a) - \ph_- (a) \in \Kt B$ for $a \in A.$
The \hm\  $\hat{F} (\mu)$ is then the composite
\[
F(A) \stackrel{F (\hat{\ph}_+) - F (\hat{\ph}_-)}{\longrightarrow}
    F(A_{\Ph}) \stackrel{F (\pi)}{\longrightarrow} F(\Kt B)
  \stackrel{F (\ep)^{-1}}{\longrightarrow} F(B),
\]
for a certain \ca\  $A_{\Ph}$, certain \hm s  $\pi$, $ \hat{\ph}_+,$
and $\hat{\ph}_-,$ and with $\ep (a) = 1 \otimes a.$
(See Definition 3.4 and the proofs of Theorems 3.7 and 4.5 in
\cite{Hg}.)
{}From this expression, it is obvious that naturality  with respect to
\hm s implies naturality  with respect to classes in $KK$-theory.
\QED

\vspace{0.6\baselineskip}

{\bf  3.2.5 Definition.}
Let $A$  be \snus.
We regard $KK^0 ({A}, {-})$ and $\EH_{{A}}$  as functors
from ${\cal KK}$ to abelian groups. (On morphisms, the first of
these sends $\et \in KK^0 (D_1, D_2)$ to Kasparov product with $\et.$)
We now define natural transformations
\[
\af : KK^0 ({A}, {-}) \to \EH_{{A}} \andeqn
\bt : \EH_{{A}} \to KK^0 ({A}, {-}).
\]

To define $\af_D,$ let $e \in K$ be a rank one \pj, let
$\io_A : A \to \Kt \OIA{A}$ be the map
$\io_A (a) = e \otimes 1 \otimes a,$
and let $[[\io_{{A}}]] \in \ET_{{A}} (A)$ denote its class in
$[[{A}, \Kt \OIA{A}]]_+ \cong \ET_{{A}} ({A}).$
(Recall that $A$ is unital, so that Remark  3.1.5 applies.)
Now let $\et \in KK^0 ({A}, {D}).$
Then $\EH_{{A}} (\et)$ is a \hm\  from
$\ET_{{A}} ({A})$ to $\ET_{{A}} ({D}).$
Define
\[
\af_D (\et) = \EH_{{A}} (\et) ([[\io_{{A}}]])
                            \in \ET_{{A}} (D).
\]

To define $\bt_D,$ let $\ch_D : D^{\#} \to \Kt \OI$
be the standard map (as in the proof of Lemma    3.1.6).
Starting with
$\et \in \ET_{{A}} (D) \subset [[{A}, D^{\#}]],$
choose a full   \amm\  $\ph : {A} \to D^{\#}$ with
$[[\ch_D]] \cdot [[\ph]] = 0$ which represents $\et.$
Then form the second suspension
\[
[[S^2 \ph]] \in [[S^2 A, S^2 D^{\#}]] \cong KK^0 ({A}, \OIA{D^+}).
\]
Since $[[S^2 \ch_D]] \cdot [[S^2 \ph]] = 0,$ split exactness of
$KK^0 ({A}, -)$ implies that $[[S^2 \ph]]$ is actually in
$KK^0 ({A}, \OIA{D}).$ In this last expression, we can use the
$KK$-equivalence of $\OI$ and ${\bf C}$, given by the unital
\hm\  ${\bf C} \to \OI,$ to drop
$\OI$.  We thus obtain an
element $\bt_D (\et) \in KK^0 ({A}, {D}).$

\vspace{0.6\baselineskip}

{\bf  3.2.6 Lemma.} The maps $\af_D$ and $\bt_D$ of the previous
definition are in fact natural transformations.

\vspace{0.6\baselineskip}

{\em Proof:} It is easy to check that both $\af$
and $\bt$ are natural with respect to
\hm s, so  naturality  with respect to classes in $KK$-theory follows
from Lemma  3.2.4. \QED

\vspace{0.6\baselineskip}

{\bf    3.2.7 Theorem.}
Let $A$ be \snus.
Then for every separable  $D,$ the maps $\af_D$ and $\bt_D$ of
Definition    3.2.4 are mutually inverse isomorphisms.

\vspace{0.6\baselineskip}

{\em Proof:}
It is convenient to prove this first under the assumptions that
$\OIA{A} \cong A$ and $\OIA{D} \cong D.$ It then follows
that the map $a \mapsto 1 \otimes a$ from $A$ to
$\OIA{A}$ is homotopic to an isomorphism, and similarly for $D.$
(This is true for $\OI$ by Theorem    2.1.5 and Proposition   2.1.11.
Therefore it is true for $\OIA{A}$ and $\OIA{D},$ hence for $A$ and
$D.$)
Thus, $A$ and $\OIA{A}$ are naturally
homotopy equivalent, and therefore also naturally equivalent in
${\cal KK}$ as well. Similar considerations apply to $D.$
Thus, $\ET_A (D)$ becomes just $[[A, \Kt D]]_+$.
The natural transformations above are then given by
\[
\af_D (\et) = \EH_A (\et) ([[\id_A]])
\]
(with $\id_A$ being the obvious map from $A$ to $\Kt A$), and
\[
\bt_D ([[\ph]]) =
             [[S^2 \ph]] \in [[S^2 A, \Kt S^2 D]] \cong KK^0 (A, D).
\]
Letting $1_A$ denote the class in $KK^0 (A, A)$ of the identity map,
we then immediately verify that
\[
\af_A (1_A) = [[\id_A]] \andeqn \bt_A ([[\id_A]]) = 1_A.
\]

We now show that these two facts imply the theorem for unital $D$.
Let $\et \in KK^0 (A, D).$ Then $\et = 1_A \times \et,$
and naturality implies that
\[
\bt_D (\af_D (1_A \times \et)) = \bt_D (\EH_A (\et)(\af_A (1_A)))
          =  \bt_A (\af_A (1_A)) \times \et = 1_A \times \et.
\]
So $\bt_D \circ \af_D = \id.$ For the other direction,
let $\mu \in \EH_A (D).$ Using Corollary    2.3.10,
represent $\mu$ as the class of
a full   \hm\  $\ph : A \to \Kt D.$
Let $\et = [[S^2 \ph]]$ be the $KK$-class determined by $[[\ph]].$
Then, identifying $\Kt K$ with
$K$ as necessary, we have
\[
\mu = \ph_*  ([[\id_A]]) = \EH_A (\et) ([[\id_A]]).
\]
The same argument as above now shows that
\[
(\af_D  \circ \bt_D)
            \left( \rule{0em}{3ex} \EH_A (\et) ([[\id_A]]) \right)
   = \EH_A (\et) ([[\id_A]]).
\]
So $\af_D \circ \bt_D = \id$ also.

The result for nonunital algebras follows from naturality,
split exactness, and the Five Lemma.

To remove the assumption that $\OIA{D} \cong D,$ use Lemma  3.1.9.

Finally, we remove the assumption that $\OIA{A} \cong A.$ Let
$\dt_0 : \OIA{\OI} \to \OI$ be an isomorphism (from Theorem  2.1.5),
and let $\dt : \OIA{\Kt \OI} \to \Kt \OI$ be the obvious corresponding
map. Define $i_D : \ET_A (D) \to \ET_{\OIA{A}} (D)$ by
\[
i_D ([[\et]]) =
      [[\dt \otimes \id_{D^+}]] \cdot [[\id_{\OI} \otimes \et]].
\]
Let $j_D : KK^0 (A, D) \to KK^0 (\OIA{A}, D)$ be the isomorphism
induced by the $KK$-equivalence of ${\bf C}$ and $\OI.$ Both $i$ and $j$
are natural transformations. Using Theorem  2.1.5 and Proposition
2.1.11, we can rewrite $j_{\OIA{D}} (\mu)$ as
$(\dt_0 \otimes \id_D)_* (1_{\OI} \otimes \mu).$ This formula and
Remark   3.1.5 imply that $i_D \circ \af_D = \af_D \circ j_D$ when $D$
is unital and $\OIA{D} \cong D.$ The previous paragraph and the
definition of $\ET_A (D)$ in terms of $[[A, D^{\#}]]_+$ now imply that
$i_D \circ \af_D = \af_D \circ j_D$ for all $D.$  A related argument
shows that also $j_D \circ \bt_D = \bt_D \circ i_D$ for all $D.$

It now suffices to prove that $i_D$ is an isomorphism for all $D.$
By naturality, split exactness, and the Five Lemma, it suffices to do so
for unital $D.$ In this case, we have
\[
i_D : [[A, \KOI{D}]]_+ \to [[\OIA{A}, \KOI{D}]]_+
\]
given by
$i_D ([[\et]]) = [[\dt \otimes \id_D]] \cdot [[\id_{\OI} \otimes \et]].$
Define a map $k_D$ in the opposite direction by restriction to
$1 \otimes A \subset \OIA{A}.$ We prove that $k_D = i_D^{-1}.$

Let $\tilde{\dt} (x) = \dt (1 \otimes x).$ Proposition 2.1.11 implies
that there is a homotopy $\tilde{\dt} \simeq \id_{\Kt \OI}.$
It is easy to check directly that $k_D \circ i_D$ sends $[[\et]]$ to
$[[(\tilde{\dt} \otimes \id_D) \circ \et]],$ so $k_D \circ i_D$ is the
identity. For the reverse composition, let $\ta_A$ be the inclusion of
$A = 1 \otimes A$ in $\OIA{A},$ and let $\ph : \OIA{\OI} \to \OIA{\OI}$
be the flip $\ph (x \otimes y) = y \otimes x.$ Then
$\ph \simeq \id_{\OIA{\OI}}$ by Proposition 2.1.11 and Theorem  2.1.5.
Therefore, for $[[\et]] \in [[\OIA{A}, \KOI{D}]]_+,$ we have
\[
(\dt \otimes \id_D) \circ (\id_{\OI} \otimes (\et \circ \ta_A))
   \simeq (\dt \otimes \id_D) \circ (\id_{\OI} \otimes \et) \circ
        (\ph \otimes \id_A) \circ  (\id_{\OI} \otimes \ta_A)
  = (\tilde{\dt} \otimes \id_D) \circ \et \simeq \et.
\]
This shows that $i_D \circ k_D$ is the identity. \QED

\vspace{0.6\baselineskip}

{\bf  3.2.8 Remark.} We used Corollary 2.3.10 in this proof because
we had it available. It is, however, not necessary for the argument.
Using methods similar to, but a bit more complicated than, the proof of
Lemma  3.2.4, one can show that if $F$ as there is in fact a functor
on homotopy classes of \amm s, then $F([[\ph]])$ is equal to $\hat{F}$
applied to the $KK$-theory class given by $\ph.$

\vspace{0.6\baselineskip}

{\bf  3.2.9 Theorem.} Let $A$ be a separable unital
nuclear simple \ca. Then for separable unital \ca s $D$, the set of
homotopy classes of full \amm s from
$A$ to $\KOI{D}$ is naturally  isomorphic to $KK^0 (A, {D})$
via the map sending an \amm\  to the $KK$-class it determines.

\vspace{0.6\baselineskip}

{\em Proof:} This follows from Theorem  3.2.7 and Remark  3.1.5. \QED

\vspace{\baselineskip}
\section{Theorems on $KK$-theory and classification}
\vspace{\baselineskip}

In this section, we present our main results. The first subsection
contains the alternate descriptions of $KK$-theory in terms of
homotopy classes and \wue\  classes of \hm s, in case the first
variable is \snus. We also give here a proof that homotopies of
automorphisms of \kfalg s can in fact be chosen to be isotopies.
The second subsection contains the classification
theorem and its corollaries.
The third subsection contains the nonclassification results.

\vspace{0.6\baselineskip}
\subsection{Descriptions of $KK$-theory}
\vspace{0.6\baselineskip}

Probably the most striking of our descriptions of $KK$-theory is the
following:

\vspace{0.6\baselineskip}

{\bf  4.1.1 Theorem.} For a separable unital nuclear simple \ca\  $A$
and a separable unital \ca\  $D,$ the obvious maps define natural
isomorphisms of abelian groups between the following three objects:

(1) The set of \wue\  classes of full \hm s from $A$ to $\KOI{D},$ with
the operation given by direct sum (Definition  1.1.3).

(2) The set of homotopy classes of full \hm s from
$A$ to $\KOI{D},$ with the operation given by direct sum as above.

(3) The group $KK^0 (A, D).$

\vspace{0.6\baselineskip}

{\em Proof:} For the purposes of this proof, denote the set in (1)
by $KU (A, D)$ and the set in (2) by $KH (A, D).$ The map from
$KH (A, D)$ to $KK^0 (A, D)$ is the one from Theorem  3.2.9.  By this
theorem, we can use $[[A, \KOI{D}]]_+$ in place of $KK^0 (A, D).$

Lemma  1.3.3 (2) implies that the map from $KU (A, D)$ to $KH (A, D)$
is well defined, and it is then clearly surjective. Injectivity is
immediate from    Theorem 2.3.7. Thus this map is an isomorphism.
Theorem  3.2.9 implies that the map from $KH (A, D)$ to
$[[A, \KOI{D}]]_+$ is injective, while Corollary    2.3.10 implies
that the map from $KU (A, D)$ to $[[A, \KOI{D}]]_+$ is surjective.
Therefore these maps are in fact both isomorphisms. It now follows that
$KU (A, D)$ and $KH (A, D)$ are both abelian groups. \QED

\vspace{0.6\baselineskip}

We now want to give a stable version of this theorem, in which
the Kasparov product will reduce exactly to composition of \hm s.
We need the following lemma. The hypotheses allow one \ct\  path of
\hm s, and require unitaries in  $U_0 ((\Kt D)^+),$ for use in the next
subsection.

\vspace{0.6\baselineskip}

{\bf   4.1.2 Lemma.} Let $A$ be \snus, let $D_0$ be
separable and unital, and
let $D = \OIA{D_0}.$  Let $t \mapsto \ph_t,$ for $t \in [0, \infty),$
be a \ct\  path of full \hm s from $\Kt A$ to $\Kt D,$ and let
$\ps : \Kt A \to \Kt D$ be a full \hm.  Assume that $[\ph_0] = [\ps]$
in $KK^0 (A, D).$
Then there is a \wue\  from $\ph$ to $\ps$ which consists of unitaries
in $U_0 ((\Kt D)^+).$

\vspace{0.6\baselineskip}

{\em Proof:}
Let $\{e_{ij}\}$ be a system of matrix units for $K.$
Identify $A$ with the subalgebra $e_{11} \otimes A$ of $\Kt A.$
Define $\ph_t^{(0)}$ and $\ps^{(0)}$ to be the restrictions of
$\ph_t$ and $\ps$ to $A.$ Then $[\ph_0^{(0)}] = [\ps^{(0)}]$
in $KK^0 (A, D).$ It follows from Theorem    4.1.1 that
$\ph_0^{(0)}$ is homotopic to $\ps^{(0)}.$ Therefore $\ph^{(0)}$ and
$\ps^{(0)}$ are homotopic as \amm s, and    Theorem 2.3.7 provides
an \wue\  $t \mapsto u_t$ in $U ((\Kt D)^+)$ from $\ph^{(0)}$ to
$\ps^{(0)}.$ Let $c \in U ((\Kt D)^+)$ be a unitary with
$c \ps^{(0)} (1) = \ps^{(0)} (1) c = \ps^{(0)} (1)$ and such that $c$ is
homotopic to $u_0^{-1}.$ Then $c$ commutes with every $\ps^{(0)} (a).$
Replacing $u_t$ by $c u_t,$ we obtain an \wue, which we
again call $t \mapsto u_t,$  from
$\ph^{(0)}$ to $\ps^{(0)}$ which is in $U_0 ((\Kt D)^+).$

Define $\overline{e}_{ij} = e_{ij} \otimes 1.$ Then in particular
$u_t \ph_t (\overline{e}_{11}) u_t^* \to \ps (\overline{e}_{11})$
as $t \to \infty.$ Therefore
there is a continuous path $t\to z_t^{(1)} \in U_0 ((\Kt D)^+)$
such that $z_t^{(1)} \to 1$ and
$z_t^{(1)}u_t \ph_t (\overline{e}_{11}) u_t^* (z_t^{(1)})^* =
             \ps (\overline{e}_{11})$
for all $t.$ We still have
$z_t^{(1)}u_t \ph_t (e_{11} \otimes a) u_t^* (z_t^{(1)})^* \to
             \ps (e_{11} \otimes a)$
for $a \in A.$

For convenience, set
$f_{i j t}= z_t^{(1)}u_t \ph_t (\overline{e}_{ij}) u_t^* (z_t^{(1)})^*,$
for all $t$ and for $1 \leq i, j \leq 2.$
For each fixed $t,$ the $f_{i j t}$ are matrix units, and
$f_{1 1 t} = \ps (\overline{e}_{11}).$
Set
$w_t =
  \ps (\overline{e}_{21}) f_{1 2 t} + 1 - f_{2 2 t} \in U((\Kt D)^+).$
Then one checks that $w_t f_{i j t} w_t^* = \ps (\overline{e}_{ij})$
for all $t$ and for $1 \leq i, j \leq 2.$
Choose $c \in U((\Kt D)^+)$ with
$c \ps (\overline{e}_{11} + \overline{e}_{22}) =
        \ps (\overline{e}_{11} + \overline{e}_{22}) c =
        \ps (\overline{e}_{11} + \overline{e}_{22})$
and $c w_1  \in U_0 ((\Kt D)^+).$ Set $z_t^{(2)} = c w_t$ for
$t \geq 1$ and extend $z_t^{(2)}$ over $[0, 1]$ to be \ct, unitary,
and satisfy $z_0^{(2)} = 1.$
This gives $z_t^{(2)} = 1$ for $t = 0,$
$z_t^{(2)} \ps (\overline{e}_{11}) =\ps (\overline{e}_{11}) z_t^{(2)}
                      = \ps (\overline{e}_{11})$
for all $t,$ and
\[
z_t^{(2)}
 \left[z_t^{(1)}u_t \ph_t (\overline{e}_{ij}) u_t^* (z_t^{(1)})^*\right]
                              (z_t^{(2)})^*
   = \ps (\overline{e}_{ij})
\]
for $t \geq 1$ and $1 \leq i, j \leq 2.$

We continue inductively, obtaining by the same method a sequence of
continuous paths $t \mapsto z_t^{(n)}$ such that
$z_t^{(n + 1)} = 1$ for $t \leq n - 1,$
\[
z_t^{(n + 1)}\left( \sum_{j = 1}^n \ps (\overline{e}_{jj}) \right) =
  \left(\sum_{j = 1}^n \ps (\overline{e}_{jj}) \right) z_t^{(n + 1)} =
  \sum_{j = 1}^n \ps (\overline{e}_{jj})
\]
for all $t,$ and
\[
z_t^{(n + 1)} \left[ \left( \prod_{k = 1}^n z_t^{(k)} \right) u_t
\ph_t (\overline{e}_{ij})
  u_t^*\left( \prod_{k = 1}^n  z_t^{(k)} \right)^* \right]
                           (z_t^{(n + 1)})^*
   = \ps (\overline{e}_{ij})
\]
for $t \geq n$ and $1 \leq i, j \leq n + 1.$

Now define
\[
z_t = \left( \prod_{k = 1}^{\infty} z_t^{(k)} \right) u_t.
\]
In a neighborhood of each $t,$ all but finitely many of the factors
are equal to $1,$ so this product yields a \ct\  path of unitaries
in $U_0 ((\Kt D)^+).$ Moreover,
$z_t \ph_t (\overline{e}_{ij}) z_t^* = \ps (\overline{e}_{ij})$
whenever $t \geq i, j,$ so that
$\lim_{t \to \infty} z_t \ph_t (\overline{e}_{ij}) z_t^* =
            \ps (\overline{e}_{ij})$
for all $i$ and $j,$ while
\[
\lim_{t \to \infty} z_t \ph_t (e_{11} \otimes a) z_t^* =
\lim_{t \to \infty}
        z_t^{(1)} u_t \ph_t (e_{11} \otimes a) u_t^* (z_t^{(1)})^* =
  \ps (e_{11} \otimes a)
\]
for all $a \in A.$ Since the $\overline{e}_{ij}$ and
$e_{11} \otimes a$ generate $\Kt A,$ this shows that $t \mapsto z_t$
is an \wue. \QED

\vspace{0.6\baselineskip}

{\bf    4.1.3 Theorem.} For a separable unital nuclear simple \ca\  $A$
and a separable unital \ca\  $D,$ the obvious maps and the
isomorphism $KK^0 (A, D) \to KK^0 (\KOI{A}, \KOI{D})$ define natural
isomorphisms of abelian groups between the following three objects:

(1) The set of \wue\  classes of full \hm s from
$\KOI{A}$ to $\KOI{D},$ with the operation given by direct sum
(as in Theorem    4.1.1).

(2) The set of homotopy classes of full \hm s from
$\KOI{A}$ to $\KOI{D},$ with the operation given by direct sum as above.

(3) The group $KK^0 (A, D).$

\noindent
Moreover, if $B$ is another a separable unital nuclear simple \ca,
then the Kasparov product
$KK^0 (A, B) \times KK^0 (B, D) \to KK^0 (A, D)$
is given in the groups in (1) and (2) by composition of \hm s.

\vspace{0.6\baselineskip}

{\em Proof:} The last statement will follow immediately from the
rest of the theorem, since if two $KK$-classes are represented by
\hm s, then their product is represented by the composition.

For the rest of the theorem, first note that the map
$KK^0 (A, D) \to KK^0 (\KOI{A}, \KOI{D})$ is a natural isomorphism
because it is induced by the $KK$-equivalence ${\bf C} \to \Kt \OI,$
given by $1 \mapsto e \otimes 1$ for some rank one \pj\  $e \in K,$
in each variable.

Now observe that the previous lemma implies that the map from the
set in (1) to $KK^0 (A, D)$ is injective. Moreover, the map from the
set in (1) to the set in (2) is well defined by Lemma    1.3.3 (2),
and is then obviously surjective. It therefore suffices to prove that
the map from the set in (2) to $KK^0 (A, D)$ is surjective, that is,
that every class in $KK^0 (A, D)$ is represented by a \hm\  from
$\KOI{A}$ to $\KOI{D}.$ It follows from Theorem    4.1.1 that
every such class is represented by a \hm\  from $A$ to $\KOI{D},$
and we obtain a \hm\  from $\KOI{A}$ to $\KOI{D}$ by tensoring
with $\id_{\Kt \OI}$ and composing with the tensor product of
$\id_D$ and an isomorphism $\KOI{\Kt \OI} \to \Kt \OI$ which is the
identity on $K$-theory. \QED

\vspace{0.6\baselineskip}

We finish this section with one other application. Following terminology
from differential topology, we define an {\em isotopy} to be a homotopy
$t \mapsto \ph_t$ in which each $\ph_t$ is an isomorphism.

\vspace{0.6\baselineskip}

{\bf  4.1.4 Theorem.} Let $A$ be a \kfalg.

(1) If $U(A)$ is connected, then two automorphisms of $A$ with the same
class in $KK^0 (A, A)$ are isotopic.

(2) Any two automorphisms of $\Kt A$ with the same
class in $KK^0 (A, A)$ are isotopic.

\vspace{0.6\baselineskip}

{\em Proof:} For (2), take $D = A$ in Lemma  4.1.2, note that
$\OIA{A} \cong A$ (Theorem  2.1.5), and note that an \wue\  with
unitaries in $U_0 ((\Kt A)^+)$ gives an isotopy, not just a homotopy.

For (1), let $\ph$ and $\ps$ be automorphisms of $A$ with the same
class in $KK^0 (A, A).$ Let $e \in K$ be a rank one projection.
Apply (2) to $\id_K \otimes \ph$ and $\id_K \otimes \ps.$
Thus, there is a unitary path $t \mapsto u_t$ in $(\Kt A)^+$ with
$u_t \ph (e \otimes a) u_t^* \to \ps (e \otimes a)$ for $a \in A.$
In particular, $u_t (e \otimes 1) u_t^* \to (e \otimes 1).$
Replacing $u_t$ by $v_t u_t$ for a suitable unitary path
$t \mapsto v_t,$ we may therefore assume
that $u_t (e \otimes 1) u_t^* = e \otimes 1$ for all $t.$ Cut down
by $e \otimes 1,$ and observe that the hypotheses imply that
$(e \otimes 1) u_0 (e \otimes 1)$ is homotopic to 1.  Now finish as in
the proof of (2). \QED

\vspace{0.6\baselineskip}
\subsection{The classification theorem}
\vspace{0.6\baselineskip}

The following theorem is the stable version of the main classification
theorem. Everything else will be an essentially immediate corollary.
In the proof, it is easy to get the existence of the isomorphism, but
we have to do some work to make sure that it has the right class in
$KK$-theory.

\vspace{0.6\baselineskip}

{\bf  4.2.1 Theorem.} Let $A$ and $B$ be \kfalg s, and suppose that
there is an invertible element $\et \in KK^0 (A, B).$
Then there is an isomorphism $\ph : \Kt A \to \Kt B$ such
$[\ph] = \et$ in $KK^0 (A, B).$

\vspace{0.6\baselineskip}

{\em Proof:} Theorems  3.2.9 and  2.1.5 provide a full
\amm\  $\af: A \to {\Kt B}$ whose class in $KK^0 (A, B)$ is
$\et.$ By Corollary    2.3.10, we may in fact take $\af$ to be
a \hm. Let $\mu : \Kt K \to K$ be an isomorphism, and set
$\ph_0 = (\mu \otimes \id_B) \circ (\id_K \otimes \af).$ Then $\ph_0$
is a nonzero (hence full) \hm\  from $\Kt A$ to ${\Kt B}$ whose class
in $KK^0 (A, B)$ is also $\et.$ Similarly, there is a full
\hm\  $\ps_0 : \Kt B \to \Kt A$
whose class in $KK^0 (B, A)$ is $\et^{-1}.$ It follows from
Theorems    4.1.3 and  2.1.5 that $\ps_0 \circ \ph_0$ is homotopic
to $\id_{\Kt A}$ and $\ph_0 \circ \ps_0$ is homotopic  to $\id_{\Kt B}.$

We now construct \hm s $\ph^{(n)} : \Kt A \to \Kt B,$
$\ps^{(n)} : \Kt B \to \Kt A,$
homotopies $\af \mapsto \tilde{\ph}^{(n)}_{\af}$ (for $\af \in [0, 1]$)
of \hm s from $\Kt A$ to $\Kt B,$ and finite subsets
$F_n \subset \Kt A$ and $G_n \subset \Kt B$
such that the following conditions are satisfied:

(1) $\ph^{(0)} = \ph_0.$

(2) Each $\ph^{(n)}$ is of the form $a \mapsto v \ph_0 (a) v^*$ for
some suitable $v \in U_0 ((\Kt B)^+),$ and similarly
each $\ps^{(n)}$ is of the form $b \mapsto u \ph_0 (b) u^*$ for
some suitable $u \in U_0 ((\Kt A)^+).$

(3) $F_0 \subset F_1 \subset \cdots$ and
$\bigcup_{n = 0}^{\infty} F_n$ is dense in $\Kt A,$ and similarly
$G_0 \subset G_1 \subset \cdots$ and
$\bigcup_{n = 0}^{\infty} G_n$ is dense in $\Kt B.$

(4) $\ph^{(n)} (F_n) \subset G_n$ and
$\ps^{(n)} (G_n) \subset F_{n + 1}.$

(5) $\| \ps^{(n)} \circ \ph^{(n)} (a) - a \| < 2^{-n}$ for
$a \in F_n$ and
$\| \ph^{(n + 1)} \circ \ps^{(n)} (b) - b \| < 2^{-n}$ for
$b \in G_n.$

(6) $\| \tilde{\ph}_{\af}^{(n + 1)} (a)
                   - \tilde{\ph}_{\af}^{(n)} (a) \| < 2^{-n}$
for $a \in F_n$ and $\af \in [0, 1].$

(7) $\tilde{\ph}_{\af}^{(n)} = \ph_0$ for $\af \geq 1 - 2^{-n}$
and $\tilde{\ph}_0^{(n)} = \ph^{(n)}.$

\noindent
This will yield the following approximately commutative diagram:

\begin{picture}(435, 160)(-15, -40)

\put( 0,80){\makebox(0,0){$A$}}
\put(60,80){\makebox(0,0){$A$}}
\put(210,80){\makebox(0,0){$A$}}
\put(270,80){\makebox(0,0){$A$}}
\put( 0,0){\makebox(0,0){$B$}}
\put(60,0){\makebox(0,0){$B$}}
\put(210,0){\makebox(0,0){$B$}}
\put(270,0){\makebox(0,0){$B$}}

\put(10,  80){\vector(1,0){40}}
\put(70,  80){\vector(1,0){40}}
\put(135, 80){\makebox(0,0){$\cdots \cdots$}}
\put(160, 80){\vector(1,0){40}}
\put(220, 80){\vector(1,0){40}}
\put(280, 80){\vector(1,0){40}}
\put(345, 80){\makebox(0,0){$\cdots \cdots$}}
\put(10,  0){\vector(1,0){40}}
\put(70,  0){\vector(1,0){40}}
\put(135, 0){\makebox(0,0){$\cdots \cdots$}}
\put(160, 0){\vector(1,0){40}}
\put(220, 0){\vector(1,0){40}}
\put(280, 0){\vector(1,0){40}}
\put(345, 0){\makebox(0,0){$\cdots \cdots$}}

\put( 30, 86){\makebox(0,0)[b]{$\id_A$}}
\put( 90, 86){\makebox(0,0)[b]{$\id_A$}}
\put(180, 86){\makebox(0,0)[b]{$\id_A$}}
\put(240, 86){\makebox(0,0)[b]{$\id_A$}}
\put(300, 86){\makebox(0,0)[b]{$\id_A$}}
\put( 30, -2){\makebox(0,0)[t]{$\id_B$}}
\put( 90, -2){\makebox(0,0)[t]{$\id_B$}}
\put(180, -2){\makebox(0,0)[t]{$\id_B$}}
\put(240, -2){\makebox(0,0)[t]{$\id_B$}}
\put(300, -2){\makebox(0,0)[t]{$\id_B$}}

\put(  0, 72){\vector(0,-1){64}}
\put( 60, 72){\vector(0,-1){64}}
\put(210, 72){\vector(0,-1){64}}
\put(270, 72){\vector(0,-1){64}}

\put( -2, 40){\makebox(0,0)[r]{$\ph^{(0)}$}}
\put( 58, 30){\makebox(0,0)[r]{$\ph^{(1)}$}}
\put(212, 60){\makebox(0,0)[l]{$\ph^{(n-1)}$}}
\put(272, 50){\makebox(0,0)[l]{$\ph^{(n)}$}}

\put( 6,  8){\vector(3,4){48}}
\put(66,  8){\vector(3,4){44}}
\put(159, 12){\vector(3,4){44}}
\put(216,  8){\vector(3,4){48}}
\put(276,  8){\vector(3,4){44}}

\put( 35, 46){\makebox(0,0)[br]{$\ps^{(0)}$}}
\put( 95, 46){\makebox(0,0)[br]{$\ps^{(1)}$}}
\put(185, 46){\makebox(0,0)[br]{$\ps^{(n-2)}$}}
\put(263, 20){\makebox(0,0)[br]{$\ps^{(n-1)}$}}

\put(  0, 95){\makebox(0,0){$\cap$}}
\put(  0,-15){\makebox(0,0){$\cup$}}
\put( 60, 95){\makebox(0,0){$\cap$}}
\put( 60,-15){\makebox(0,0){$\cup$}}
\put(210, 95){\makebox(0,0){$\cap$}}
\put(210,-15){\makebox(0,0){$\cup$}}
\put(270, 95){\makebox(0,0){$\cap$}}
\put(270,-15){\makebox(0,0){$\cup$}}

\put(  0,110){\makebox(0,0){$F_{0}$}}
\put(  0,-30){\makebox(0,0){$G_{0}$}}
\put( 60,110){\makebox(0,0){$F_{1}$}}
\put( 60,-30){\makebox(0,0){$G_{1}$}}
\put(210,110){\makebox(0,0){$F_{n-1}$}}
\put(210,-30){\makebox(0,0){$G_{n-1}$}}
\put(270,110){\makebox(0,0){$F_{n}$}}
\put(270,-30){\makebox(0,0){$G_{n}$}}

\end{picture}

\noindent
The diagram will remain approximately commutative if we replace
each $\ph^{(n)}$ by $\tilde{\ph}_{\af}^{(n)}$ (with $\af \in [0, 1]$
fixed) and delete the diagonal arrows.

The proof is by induction on $n.$ We start by choosing finite
sets $F^{(0)}_0 \subset F^{(0)}_1 \subset \cdots \subset \Kt A$ with
$\bigcup_{n = 0}^{\infty} F^{(0)}_n$ dense in $\Kt A,$ and similarly
$G^{(0)}_0 \subset G^{(0)}_1 \subset \cdots \subset \Kt B$ with
$\bigcup_{n = 0}^{\infty} G^{(0)}_n$ dense in $\Kt B.$
For the initial step of the induction, we take $F_0 = F^{(0)}_0,$
$\ph^{(0)} = \ph^{(0)}_{\af} = \ph_0,$ and
$G_0 = G_0^{(0)} \cup \ph^{(0)} (F_0).$ We then assume we are given
$F_k,$ $\ph^{(k)},$ $G_k,$ and $\ph^{(k)}_{\af}$ for $0 \leq k \leq n$
and $\ps^{(k)}$ for $0 \leq k \leq n - 1,$ and we construct
$\ps^{(n)},$ $F_{n + 1},$ $\ph^{(n + 1)},$ $G_{n + 1},$ and
$\af \mapsto \tilde{\ph}_{\af}^{(n + 1)}.$
That is, we are given the diagram above through the column containing
$F_n$ and $G_n,$ as well as the corresponding homotopies
$\tilde{\ph}^{(k)},$ and we construct the next rectangle (consisting
of two triangles) and the corresponding homotopy
$\tilde{\ph}^{(n + 1)}.$

Define $\sm : \Kt A \to C([0, 1]) \otimes \Kt A$ by
$\sm (a) (\af) = \ps_0 (\tilde{\ph}^{(n)}_{\af} (a)).$
Note that $\sm$ is homotopic to $a \mapsto 1 \otimes \ps_0 (\ph_0 (a)),$
and so has the same class in $KK$-theory as
$a \mapsto 1 \otimes a.$ Lemma   4.1.2 provides a unitary path
$(\af, t) \mapsto u_{\af, t} \in U_0 ((\Kt A)^+)$
such that
\[
\lim_{t \to \infty} \sup_{\af \in [0, 1]}
   \| u_{\af, t} \ps_0 (\tilde{\ph}^{(n)}_{\af} (a)) u_{\af, t}^* - a \|
                    = 0
\]
for all $a \in \Kt A.$
Next, define an \amm\  $\ta$ from $\Kt B$ to $C([0, 1]) \otimes \Kt B$
by $\ta_t (b) (\af) = \ph_0 (u_{\af, t} \ps_0 (b) u_{\af, t}^*).$
Then $\ta$ is homotopic to $b \mapsto 1 \otimes \ph_0 (\ps_0 (b)),$
and so has the same class in $KK$-theory as
$b \mapsto 1 \otimes b.$ Again by Lemma   4.1.2,
there is a unitary path
$(\af, t) \mapsto v_{\af, t} \in U_0 ((\Kt B)^+)$
such that
\[
\lim_{t \to \infty} \sup_{\af \in [0, 1]}
   \| v_{\af, t} \ph_0 (u_{\af, t} \ps_0 (b) u_{\af, t}^*) v_{\af, t}^*
                                                        - b \| = 0
\]
for all $b \in \Kt B.$

Since
$\tilde{G} = G_{n} \cup
                \bigcup_{\af \in [0, 1]} \tilde{\ph}^{(n)}_{\af} (F_n)$
is a compact subset of $\Kt B,$ we can choose $T$ so large that
\[
\| v_{\af, t} \ph_0 (u_{\af, t} \ps_0 (b) u_{\af, t}^*) v_{\af, t}^*
                              - b \| < 2^{-(n + 1)}
\]
for all $b \in \tilde{G}$ and $t \geq T.$ Increasing $T$ if necessary,
we can also require
\[
\| u_{\af, t} \ps_0 (\tilde{\ph}^{(n)}_{\af} (a)) u_{\af, t}^* - a \|
                                           < 2^{-(n + 1)}
\]
for all $a \in F_n$ and $t \geq T.$ Now define
\[
\ps^{(n)} (b) = u_{0, T} \ps_0 (b) u_{0, T}^* \andeqn
\ph^{(n + 1)} (a) = v_{0, T} \ph_0 (a) v_{0, T}^*,
\]
and
\[
F_{n + 1} = F_{n + 1}^{(0)} \cup F_n \cup \ps^{(n)} (G_n) \andeqn
G_{n + 1} = G_{n + 1}^{(0)} \cup G_n \cup \ph^{({n + 1})} (F_{n + 1}).
\]
The relevant parts of conditions (2)--(4) are then certainly satisfied.
For (5), we have in fact
\[
\| \ps^{(n)} \circ \ph^{(n)} (a) - a \|
   = \| u_{0, T} \ps_0 (\tilde{\ph}_0^{(n)} (a)) u_{0, T}^* - a \|
   < 2^{-(n + 1)}
\]
for $a \in F_n$ by the choice of $T,$ and similarly
\[
\| \ph^{(n + 1)} \circ \ps^{(n)} (b) - b \| =
  \| v_{0, T} \ph_0 (u_{0, T} \ps_0 (b) u_{0, T}^* ) v_{0, T}^* - b\|
 < 2^{-(n + 1)}
\]
for $b \in G_n.$

Now choose a \ct\  function $f : [0, 1 - 2^{-(n + 1)}) \to [T, \infty)$
such that $f ( \af) = T$ for $0 \leq \af \leq 1 - 2^{-n}$ and
$f( \af) \to \infty$ as $\af \to 1 - 2^{-(n + 1)}.$
Define $\af \mapsto \tilde{\ph}^{(n + 1)}_{\af}$ by
\[
\tilde{\ph}^{(n + 1)}_{\af} (a) =
\left\{ \begin{array}{lll}
  v_{\af, f(\af)} \ph_0 (a) v_{\af, f(\af)}^* & \hspace{1cm} &
                                       0 \leq \af < 1 - 2^{-(n + 1)} \\
  \ph_0 (a)  & \hspace{1cm} & 1 - 2^{-(n + 1)} \leq \af \leq 1.
\end{array} \right.
\]
We first have to show that the functions
$\af \mapsto \tilde{\ph}^{(n + 1)}_{\af} (a)$ are \ct\  at
$1 - 2^{-(n + 1)}$
for $a \in \Kt A.$ Set $\af_0 = 1 - 2^{-(n + 1)},$ and consider
$\af$ with $1 - 2^{-n} \leq \af < 1 - 2^{-(n + 1)}.$ By the induction
hypothesis, we then have $\tilde{\ph}^{(n)}_{\af} (a) = \ph_0 (a).$ For
$a \in \Kt A,$ set $b = \ph_0 (a);$ then
\beqr
\lefteqn{
\|\tilde{\ph}^{(n + 1)}_{\af} (a) - \tilde{\ph}^{(n + 1)}_{\af_0} (a) \|
} \\
 &  \leq &  \| a -
u_{\af, f(\af)} \ps_0 (\tilde{\ph}^{(n)}_{\af} (a)) u_{\af, f(\af)}^* \|
+
\| v_{\af, f(\af)} \ph_0 ( u_{\af, f(\af)} \ps_0 (b)
                                u_{\af, f(\af)}^* )v_{\af, f(\af)}^*
       - b \|.
\eeqr
The requirement that $f(\af) \to \infty$ as $\af \to 1 - 2^{-(n + 1)},$
together with the condition of uniformity in $\af$ in the limits
used in the choices of $u_{\af, t}$ and $v_{\af, t},$ implies that both
terms on the right converge to $0.$ So the required continuity holds.

The relevant part of condition (7) is satisfied by definition, so it
remains only to check (6). We may assume $\af < 1 - 2^{-(n + 1)}.$
So let $a \in F_n.$ Then
$b = \tilde{\ph}^{(n)}_{\af} \in \tilde{G}.$ So
\beqr
\lefteqn{\| \tilde{\ph}^{(n + 1)}_{\af} (a)
                                - \tilde{\ph}^{(n)}_{\af} (a) \|}  \\
 & \leq & \sup_{\af \in [0, 1], t \geq T}
          \| v_{\af, t} \ph_0 (a) v_{\af, t}^*
                                - \tilde{\ph}^{(n)}_{\af} (a) \|   \\
 & \leq & \sup_{\af \in [0, 1], t \geq T}
   \| a - u_{\af, t} \ps_0 (\tilde{\ph}^{(n)}_{\af} (a)) u_{\af, t}^* \|
+ \sup_{\af \in [0, 1], t \geq T}
\| v_{\af, t} \ph_0 ( u_{\af, t} \ps_0 (b)
                           u_{\af, t}^* )v_{\af, t}^* - b \|  \\
 & < & 2^{-(n + 1)} + 2^{-(n + 1)} = 2^{-n}.
\eeqr
This proves (6), and finishes the inductive construction.  Note that the
set $\bigcup_{n = 0}^{\infty} F_n$ is dense in $\Kt A$ because it
contains the dense subset $\bigcup_{n = 0}^{\infty} F_n^{(0)},$ and
similarly $\bigcup_{n = 0}^{\infty} G_n$ is dense in $\Kt B.$

We now define $\ph : \Kt A \to \Kt B$ by
$\ph (a) = \lim_{n \to \infty} \ph^{(n)} (a),$ and define
$\ps : \Kt B \to \Kt A$ and the homotopy
$\tilde{\ph} : \Kt A \to C([0, 1]) \otimes \Kt B$ analogously. As in
Section 2 of \cite{Ell2}, these limits all exist and define  \hm s;
moreover, $\ps \circ \ph = \id_{\Kt A},$ $\ph \circ \ps = \id_{\Kt B},$
$\tilde{\ph}_0 = \ph,$ and $\tilde{\ph}_1 = \ph_0.$ So $\ph$ is an
isomorphism from $\Kt A$ to $\Kt B$ which is homotopic to $\ph_0$ and
therefore satisfies $[\ph] = \et$ in $KK^0 (A, B).$ \QED

\vspace{0.6\baselineskip}

{\bf    4.2.2 Corollary.} Let $A$ and $B$ be \kfalg s, and suppose that
there is an invertible element $\et \in KK^0 (A, B)$ such that
$[1_A] \times \et = [1_B].$
Then there is an isomorphism $\ph : A \to B$ such
$[\ph] = \et$ in $KK^0 (A, B).$

\vspace{0.6\baselineskip}

{\em Proof:} The previous theorem provides an isomorphism
$\af : \Kt A \to \Kt B$ such that $[\af] = \et$ in $KK^0 (A, B).$
Choose a rank one \pj\  $e \in K.$ Then
$[\af (e \otimes 1_A)] = [1_A] \times \et = [e \otimes 1_B]$ in
$K_0 (B).$ Since $\Kt B$ is purely infinite simple, it follows that
there is a unitary $u \in (\Kt B)^+$ such that
$u \af (e \otimes 1_A) u^* = e \otimes 1_B.$ Define
$\ph (a) = u \af (e \otimes a) u^*,$ regarded as an element of
$(e \otimes 1_B) (\Kt B) (e \otimes 1_B) = B.$ \QED

\vspace{0.6\baselineskip}

The remaining corollaries require some hypotheses on the
Universal Coefficient Theorem. (See \cite{RS}.)
The following terminology is convenient.

\vspace{0.6\baselineskip}

{\bf  4.2.3 Definition.} Let $A$ and $D$ be separable nuclear \ca s.
We say that the pair $(A, D)$ {\em satisfies the Universal Coefficient
Theorem} if the sequence
\[
0 \longrightarrow {\rm Ext}_1^{\bf Z} (K_* (A), K_* (D)) \longrightarrow
    KK^0 (A, D) \longrightarrow {\rm Hom} (K_* (A), K_* (D))
    \longrightarrow 0
\]
of Theorem 1.17 of \cite{RS} is defined and exact. (Note that the second
map is always defined, and the first map is the inverse of a map that is
always defined.) We further say that $A$ satisfies the Universal
Coefficient Theorem if $(A, D)$ does for every separable \ca\  $D.$

\vspace{0.6\baselineskip}

{\bf    4.2.4 Theorem.} Let $A$ and $B$ be separable nuclear \pisca s
which satisfy the Universal Coefficient Theorem.
Assume that $A$ and $B$ are either both unital or both nonunital.
If there is a graded isomorphism
$\af : K_* (A) \to K_* (B)$ which (in the unital case) satisfies
$\af_* ([1_A]) = [1_B],$ then there is an isomorphism
$\ph : A \to B$ such that $\ph_* = \af.$

\vspace{0.6\baselineskip}

{\em Proof:} The proof of Proposition 7.3 of \cite{RS} shows that
there is a $KK$-equivalence $\et \in KK^0 (A, B)$ which induces
$\af.$ Now use Theorem    4.2.1 or Corollary    4.2.2 as appropriate.
\QED

\vspace{0.6\baselineskip}

This theorem gives all the classification results of \cite{Rr1},
\cite{Rr2}, \cite{LP1}, \cite{ER}, \cite{LP2}, \cite{Rr4},
\cite{LP3}, \cite{Ph2}, and \cite{Rr5}. Of course, we have
used the main technical theorem of \cite{Rr1}, as well as
substantial material from \cite{LP2}, in the proof. We do not obtain
anything new about the Rokhlin property of \cite{BKRS}; indeed, our
results show that the \ca s of \cite{Rr4} are classifiable as long as
they are purely infinite and simple, regardless of whether the Rokhlin
property is satisfied. On the other hand, the Rokhlin property has been
verified in many cases; see \cite{Ksh1} and \cite{Ksh2}.

We finish this section by giving some further corollaries.
Let $\Class$ be the ``classifiable class'' given in Definition 5.1
of \cite{ER}, and let $\Boot$ denote the bootstrap category of
\cite{RS}, for which the Universal Coefficient Theorem was shown to
hold (Theorem 1.17 of \cite{RS}).

\vspace{0.6\baselineskip}

{\bf    4.2.5 Theorem.} Let $G_0$ and $G_1$ be countable abelian groups,
and let $g \in G_0$.
Then:

(1) There is a \kfalg\  algebra $A \in \Boot$ such that
\[
\left( K_0 (A), [1_{A}], K_1 (A) \right) \cong (G_0, g, G_1).
\]

(2) There is a separable nuclear nonunital \pisca\  $A \in \Boot$
such that
\[
\left( K_0 (A), K_1 (A) \right) \cong (G_0, G_1).
\]

\vspace{0.6\baselineskip}

{\em Proof:} The construction of Theorem 5.6 of \cite{ER}
gives algebras which are easily seen to be in $\Boot.$ \QED

\vspace{0.6\baselineskip}

{\bf    4.2.6 Corollary.} Every \ca\  in $\Class$ is in $\Boot.$
Every \pisca\  in $\Boot,$ and more generally every
separable nuclear \pisca\  satisfying the Universal Coefficient Theorem,
is in $\Class.$

\vspace{0.6\baselineskip}

{\em Proof:} The first part follows immediately from the previous
theorem, since it follows from the definition of $\Class$ that any
$A \in \Class$ must be isomorphic to the \ca\  of that theorem
with the same
$K$-theory. The second part follows from Theorem
   4.2.4, since Theorem 1.17 of \cite{RS} states that every
\ca\  in $\Boot$ satisfies the Universal Coefficient Theorem. \QED

\vspace{0.6\baselineskip}

{\bf    4.2.7 Corollary.} Let $A \in \Class,$ and let $B$ be
a separable nuclear unital simple \ca\  which satisfies the
Universal Coefficient Theorem. (In particular, $B$ could be a unital
simple \ca\  in $\Boot.$) Then $A \otimes B \in \Class.$

\vspace{0.6\baselineskip}

{\em Proof:} The \ca\  $A \otimes B$ is \snus, and Theorem 7.7
of \cite{RS} (and the remark after this theorem) shows that it
satisfies the Universal Coefficient Theorem. Furthermore, $A$ is
approximately divisible by Corollary  2.1.6, and it follows from
the remark after Theorem 1.4 of \cite{BKR} that $A \otimes B$ is
approximately divisible. Clearly $A \otimes B$ is infinite, so it is
purely infinite by Theorem 1.4 (a) of \cite{BKR}. The result now follows
from the previous corollary. \QED

\vspace{0.6\baselineskip}

{\bf    4.2.8 Corollary.} The class $\Class$ is closed under tensor
products.

\vspace{0.6\baselineskip}

{\bf    4.2.9 Corollary.} For any $m, \, n \geq 2,$ we have
$\OA{m} \otimes \OA{n} \in \Class.$ In particular, if $m - 1$ and
$n - 1$ are relatively prime, then $\OA{m} \otimes \OA{n} \cong \OA{2}.$

\vspace{0.6\baselineskip}

{\bf   4.2.10 Corollary.} Let $A_1$ and $A_2$ be two higher
dimensional noncommutative toruses  of the same dimension, and
let $B$ be any simple Cuntz-Krieger algebra. Then
$A_1 \otimes B \cong A_2 \otimes B.$

\vspace{0.6\baselineskip}

{\em Proof:} The K\"{u}nneth formula \cite{Sc1} shows that
$A_1 \otimes B$ and $A_2 \otimes B$ have the same $K$-theory. \QED

\vspace{0.6\baselineskip}

{\bf   4.2.11 Theorem.} Let $A$ be a \kfalg\  satisfying the
Universal Coefficient Theorem. Let $A^{\rm op}$ be the opposite
algebra, that is, $A$ with the multiplication reversed but all other
operations the same. Then $A \cong A^{\rm op}.$

\vspace{0.6\baselineskip}

{\em Proof:} The identity map from $A$ to  $A^{\rm op}$ is an
antiisomorphism which induces an isomorphism on $K$-theory sending
$[1_A]$ to $[1_{A^{\rm op}}].$ Also, the pair $(A^{\rm op}, B)$
always satisfies the Universal Coefficient Theorem, because
$(A, B^{\rm op})$ does. \QED

\vspace{0.6\baselineskip}

By way of contrast, we note that Connes has shown \cite{Cn} that there
is a type III factor not isomorphic to its opposite algebra. It is also
known (although apparently not published) that there are nonsimple
separable nuclear (even type I) \ca s not isomorphic to their opposite
algebras.

\vspace{0.6\baselineskip}
\subsection{Nonclassification}
\vspace{0.6\baselineskip}

In this subsection, we give some results which show how badly the
classification theorem fails if the algebras are not nuclear. The
results are mostly either proved elsewhere or follow fairly easily
from results proved by other people. There are three main results.
First, nonnuclear separable \pisca s need not be approximately divisible
in the sense of \cite{BKR}, but whenever $A$ is a \pisca, then
$\OIA{A}$ is an approximately divisible \pisca\  with exactly the same
K-theoretic invariants. Second, there are infinitely many
mutually nonisomorphic
approximately divisible separable exact unital \pisca s $A$ satisfying
$K_* (A) = 0.$ Finally, given arbitrary countable abelian groups
$G_0$ and $G_1,$ and $g \in G_0,$ there are uncountably many
mutually nonisomorphic
approximately divisible separable unital \pisca s $A$ satisfying
$K_j (A) \cong G_j$ with $[1] \mapsto g_0.$ Unfortunately these algebras
are not exact, and it remains unknown whether the same is true with
the additional requirement of exactness.

The first result is taken straight from a paper of Dykema and R\o rdam.

\vspace{0.6\baselineskip}

{\bf   4.3.1  Theorem.} (\cite{DR}, Theorem 1.4)
There exists a separable unital \pisca\  which is not
approximately divisible.

\vspace{0.6\baselineskip}

{\bf   4.3.2 Remark.}
In fact, there exists a separable unital \pisca\  $A$ which is not
approximately divisible and such that $K_* (A) = 0.$

One way to see this is to modify the proof of Proposition 1.3 of
\cite{DR} so as to ensure that $K_* (A_n) \to K_* (B)$ is injective
for all $n.$ This is done by enlarging the set $X_{n + 1}$ in the
proof so as to include appropriate partial isometries (implementing
equivalences between \pj s) and paths of unitaries
(implementing triviality of classes of unitaries in $K_1$). See the
proof of Theorem 4.3.11 below for this argument in a related context.

\vspace{0.6\baselineskip}

The second result is a fairly easy consequence of a computation
of Cowling and Haagerup and of unpublished work of Haagerup.
The key invariant is described in the following definition.
I am grateful to Uffe Haagerup for explaining the properties of this
invariant and where to find proofs of them.

\vspace{0.6\baselineskip}

{\bf   4.3.3 Definition.}
(Haagerup \cite{Ha}; also see Section 6 of \cite{CHa}.)
Let $A$ be a \ca. Define $\Ld (A)$ to be the infimum of numbers
$C$ such that there is a net of finite rank operators
$T_{\af} : A \to A$ for which $\| T_{\af} (a) - a \| \to 0$ for all
$a \in A$ and the completely bounded norms satisfy
$\sup_{\af} \|T_{\af} \|_{\rm cb} \leq C.$
Note that $\Ld (A) =\infty$ if no such $C$ exists,
that is, if $A$ does not have the completely bounded approximation
property.

\vspace{0.6\baselineskip}

There is a similar definition for von Neumann algebras, in which
$T_{\af} (a)$ is required to converge to $a$ in the weak operator
topology. (See \cite{Ha} and Section 6 of \cite{CHa}.)
There is also a definition of $\Ld (G)$ for a locally compact group
$G,$ using completely bounded norms of multipliers of $G$
which converge to $1$ uniformly on compact sets;
see \cite{Ha} and Section 1 of \cite{CHa}. We do not formally
state the definitions, but we recall the following theorems
from \cite{Ha} (restated as Propositions 6.1 and 6.2 of \cite{CHa}):

\vspace{0.6\baselineskip}

{\bf   4.3.4 Theorem.}
Let $\Gm$ be a discrete group, and let $C^*_{\rm r} (\Gm)$ and
$W^* (\Gm)$ be its reduced \ca\  and von Neumann algebra
respectively.
Then
$\Ld (\Gm) = \Ld (C^*_{\rm r} (\Gm)) = \Ld ( W^* (\Gm)).$

\vspace{0.6\baselineskip}

{\bf   4.3.5 Theorem.}
Let $G$ be a second countable locally compact group, and let
$\Gm$ be a lattice in $G.$ Then $\Ld (\Gm) = \Ld (G).$

In Section 6 of \cite{CHa}, Cowling and Haagerup exhibit type II$_1$
factors $M_n$ with $\Ld (M_n) = 2 n - 1.$
Using the same results on groups, we exhibit simple
\ca s with the same values of $\Ld.$

\vspace{0.6\baselineskip}

{\bf   4.3.6 Proposition.}
Let $\Gm_n^0$ be as in Corollary 6.6 of \cite{CHa}. Then
$A_n = C^*_{\rm r} ( \Gm_n^0)$ is a simple separable unital \ca\  which
satisfies $\Ld (A_n) = 2 n - 1.$

\vspace{0.6\baselineskip}

We recall that $\Gm_n^0$ is the quotient by its center of a particular
lattice $\Gm_n$ in the simple Lie group ${\rm Sp} (n, 1).$

\vspace{0.6\baselineskip}

{\em Proof of Proposition   4.3.6:}
It is shown in the proof of Corollary 6.6 of \cite{CHa}
that $\Ld (\Gm_n^0) = 2 n - 1.$ (This follows from the computation
$\Ld ( {\rm Sp} (n, 1)) = 2 n - 1,$ which is the main result
of \cite{CHa}, together with Theorem   4.3.5  above and
Proposition 1.3 (c) of \cite{CHa}.) Therefore $\Ld (A_n) = 2 n - 1$
by Theorem   4.3.4. Clearly $A_n$ is separable and unital.
Simplicity of $A_n$ follows from Theorem 1 of \cite{BCD}, applied
to the quotient of ${\rm Sp} (n, 1)$ by its center,
because (as observed in the introduction to \cite{BCD})
lattices satisfy the density hypothesis of that theorem. \QED

\vspace{0.6\baselineskip}

The algebras $A_n$ are not purely infinite, and their K-theory
seems to be unknown. So we will tensor them with $\OA{2}.$
For this, we need the following result.

\vspace{0.6\baselineskip}

{\bf   4.3.7. Lemma.}
Let $A$ be any \ca, and let $B$ be unital and nuclear. Then
$\Ld (A \otimes B) = \Ld (A).$

\vspace{0.6\baselineskip}

For von Neumann algebras, it is known \cite{SS} that
$\Ld (M \otimes N) = \Ld (M) \Ld (N).$ We presume, especially
in view of Remark 3.5 of \cite{SS}, that the analgous statement
is true for \ca s as well. However, the special case in the lemma
is sufficient here.

\vspace{0.6\baselineskip}

{\em Proof of Lemma   4.3.7:}
If $S : A_1 \to A_2$ and $T: B_1 \to B_2$ are
completely bounded, then the map
$S \otimes_{\rm min} T : A_1 \otimes_{\rm min} B_1 \to
       A_2 \otimes_{\rm min} B_2$
is completely bounded, and satisfies
$\| S \otimes_{\rm min} T \|_{\rm cb} =
   \| S \|_{\rm cb} \| T \|_{\rm cb}$
by Theorem 10.3 of \cite{Pl}. 
In Definition   4.3.3, one need only consider elements $a$ of a dense
subset, and so it follows that
$\Ld (A \otimes_{\rm min} B) \leq \Ld (A) \Ld (B)$ for any
\ca s $A$ and $B.$ For $B$ nuclear, we have $\Ld (B) = 1,$
so $\Ld (A \otimes B) \leq \Ld (A).$

For the reverse inequality, let
$R_{\af} : A \otimes B \to A \otimes B$
be finite rank operators such that $\| R_{\af} (x) - x \| \to 0$ for all
$x \in A \otimes B.$
Choose any state $\om$ on $B,$ and define
$T_{\af} : A \to A$ by
$T_{\af} (a) = (\id_B \otimes \om) \circ R_{\af} (a \otimes 1).$
Theorem 10.3 of \cite{Pl} implies that
$\| T_{\af} \|_{\rm cb} \leq \| R_{\af} \|_{\rm cb}.$
Also, clearly $\| T_{\af} (a) - a \| \to 0$ for all $a \in A.$
So $\Ld (A) \leq \Ld (A \otimes B).$
\QED

\vspace{0.6\baselineskip}

{\bf   4.3.8 Theorem.}
There exist infinitely many mutually nonisomorphic separable exact
unital purely infinite simple \ca s $B$ satisfying
$K_* (B) = 0$ and $\OIA{B} \cong B.$ In particular, these algebras
are approximately divisible in the sense of \cite{BKR}.

\vspace{0.6\baselineskip}

{\em Proof:} 
Let $A_n = C^*_{\rm r} ( \Gm_n^0)$ as in Proposition   4.3.6.
Set $B_n = \OT{A_n}.$ Clearly $B_n$ is separable and unital.
Furthermore, $B_n$ is purely infinite simple by the proof of
Corollary 4.2.7.
We have $\OIA{B_n} \cong B_n$ because $\OIA{\OA{2}} \cong \OA{2}.$
The algebras $B_n$ are mutually nonisomorphic because
$\Ld (B_n) = 2 n - 1,$ by the previous lemma and Proposition   4.3.6.

It remains to check exactness. The proof of Corollary 3.12 of
\cite{DH} shows that if $\Ld (A)$ is finite, then $A$ has the
slice map property (as defined, for example, in Remark 9 of
\cite{Ws1}, where it is called Property S), and this
property implies exactness (see, for example, Section 2.5 of
\cite{Ws2}). \QED

\vspace{0.6\baselineskip}

Our third result is based on the theorem of Junge and Pisier
that for $n \geq 3$ the collection of $n$-dimensional operator spaces
is not separable in the completely bounded analog of the
Banach-Mazur distance.

\vspace{0.6\baselineskip}

{\bf 4.3.9 Definition.} (\cite{JP})
Let $E$ and $F$ be operator spaces of the same finite dimension.
Then
\[
d_{\rm cb} (E, F) = \inf \{ \|T\|_{\rm cb} \|T^{-1}\|_{\rm cb} :
 T \,\, {\rm \, is \, a \, linear \, bijection \, from}
 \,\, E  \,\,{\rm to}  \,\, F\},
\]
and $\dt_{\rm cb} (E, F) = \log ( d_{\rm cb} (E, F) ).$

\vspace{0.6\baselineskip}

{\bf 4.3.10 Theorem.} (Theorem 2.3 of \cite{JP})
Let ${\rm OS}_n$ be the set of all complete isometry classes of
$n$-dimensional operator spaces. Let $n \geq 3.$ Then
$({\rm OS}_n, \dt_{\rm cb})$ is an inseparable metric space.

\vspace{0.6\baselineskip}

{\bf 4.3.11 Theorem.}
Let $G_0$ and $G_1$ be countable abelian groups, and let $g \in G_0.$
Then there exist uncountable many mutually nonisomorphic
separable unital \pisca s $A,$ each with $K_0 (A) \cong G_0$ in such
a way that $[1] \mapsto g$ and $K_1 (A) \cong G_1,$ and each
satisfying $\OIA{A} \cong A.$

\vspace{0.6\baselineskip}

{\em Proof:}
If $A$ is a separable \ca, then the set of
(complete isometry classes of) $n$-dimensional operator subspaces
of $A$ is separable (by Proposition 2.6 (a) of \cite{JP}).
By the previous theorem, it therefore suffices to show that
if $E$ is a finite dimensional operator space then there exists
a \ca\  $B$ having the properties claimed in the theorem and
such that $E$ is completely isometric to a subspace of $B.$

Since $E$ is a finite dimensional operator space, it is a
subspace of a separable \ca\ $A.$ Represent $A$ on a separable
Hilbert space $H$ with infinite multiplicity, and follow this
representation with the quotient map from $L (H)$ to the Calkin
algebra $Q.$ This gives a completely isometric embedding of
$E$ in $Q.$ For convenience, we identify $E$ with its image.
Let $u \in Q$ be the image of the unilateral shift; note that
$[u]$ generates $K_1 (Q)$ and that $K_0 (Q) = 0.$
Let $B_0 = C^* (E, 1, u) \subset Q.$
We now construct by induction an increasing sequence
$B_0 \subset B_1 \subset B_2 \subset \cdots \subset Q$ of
separable \ca s such that $B_{2 n + 1}$ is simple and such that
every nonzero \pj\  in $B_{2 n - 1}$ is \mvn\  to $1$ in
$B_{2 n},$ every selfadjoint element of $B_{2 n - 1}$ is a
limit of selfadjoint elements of $B_{2 n}$ with finite spectrum,
and every unitary in $U (B_{2 n - 1}) \cap U_0 (Q)$ is homotopic
to $1$ in $B_{2 n}.$

Given $B_{2 n},$ we choose $B_{2 n + 1}$ to be any separable
simple \ca\  with $B_{2 n} \subset B_{2 n + 1} \subset Q.$
Such a subalgebra exists by Proposition 2.2 of \cite{Bl0} and
the simplicity of $Q.$
Given $B_{2 n - 1},$ we note that it suffices to have the
required elements of $B_{2 n}$ only for countable dense subsets
$S_1$ of the nonzero \pj s in $B_{2 n - 1},$
$S_2$ of the selfadjoint elements in $B_{2 n - 1},$ and $S_3$ of
the unitaries in $U (B_{2 n - 1}) \cap U_0 (Q).$ 
For each $p \in S_1,$ since $p$ is \mvn\  to $1$ in $Q,$ we can
choose an isometry $v \in Q$ such that $v^* v = 1$ and $v v^* = p.$
Let $T_1$ be the set of all these as $p$ runs through $S_1.$ For
each $a \in S_2,$ since $Q$ has real rank zero, there is a sequence
$(b_n)$ in $Q$ consisting of selfadjoint elements with finite spectrum
such that $b_n \to a.$ Let $T_2$ be the set of all terms of all
such sequences as $a$ runs through $S_2.$
For each $u \in S_3,$ since $u \in U_0 (Q),$ there is a unitary
path $t \mapsto v (t)$ in $Q$ with $v (0) = 1$ and $v (1) = u.$
Let $T_3$ consist of all $v (t)$ for $t \in [0, 1] \cap \Q$ as
$u$ runs through $S_3.$ Then take $B_{2 n}$ to be the
C*-subalgebra of $Q$ generated by $B_{2 n - 1}$ and
$T_1 \cup T_2 \cup T_3.$ This subalgebra is separable because
$B_{2 n - 1}$ is separable and $T_1 \cup T_2 \cup T_3$ is countable.

Now set $B = \overline{\bigcup_{n = 0}^{\infty} B_n}.$ Then
$B$ is simple because it is the direct limit of the simple \ca s
$B_{2 n + 1}.$ From the construction of $B_{2 n},$ it is clear
that $B$ is unital and separable, contains the operator space $E,$
has real rank zero, that all nonzero \pj s in $B$ are
\mvn\  to $1,$ and that $U (B) \cap U_0 (Q) \subset U_0 (B).$
The third and fourth properties imply that $B$ is purely infinite and
$K_0 (B) = 0.$ The last property implies that $K_1 (B) \to K_1 (Q)$
is injective. But this map is also surjective, since $B_0$ contains
a unitary whose class generates $K_1 (Q).$ So $K_1 (B) \cong \Z.$

Taking $A = \OIA{B}$ (which has the same K-theory by the
K\"{u}nneth formula \cite{Sc1}), we obtain the statement of the
theorem for the special case $G_0 = 0,$ $g = 0,$ and $G_1 = \Z.$
For the general case, choose (by Theorem 4.2.5) a \kfalg\ $D$ satisfying
the Universal Coefficient Theorem and such that
$K_0 (D) \cong G_1$ and $K_1 (D) \cong G_0.$
(We don't actually need $D$ to be
purely infinite here, but it must be in the bootstrap category of
\cite{Sc1}.)
Then $D \otimes B$ is purely infinite and simple, and has the
right K-theory by the K\"{u}nneth formula, except that $[1] = 0.$
Choose a \pj\  $p \in D \otimes B$ such that
the isomorphism $K_0 (D \otimes B) \cong G_0$ sends $[p]$ to $g.$
Then the \ca\  $A = \OIA{p (D \otimes B) p}$ satisfies all the
conditions of the theorem and contains the given operator space $E.$
\QED

\vspace{0.6\baselineskip}

{\bf 4.3.12 Remark.}
Simplicity and pure infiniteness
of $\overline{\bigcup_{n = 0}^{\infty} B_n}$
in the proof above can also be arranged by the method of
the proof of Proposition 1.3 of \cite{DR}.
Versions of the construction here have been used many times before.

\vspace{0.6\baselineskip}

{\bf 4.3.13 Remark.}
The invariant used here, the set of finite dimensional operator
spaces contained in $A,$ does not distinguish between any two
separable exact \pisca s. (Any separable exact \ca\  embeds in $\OA{2}$
by Theorem 2.9 of \cite{KP}, and $\OA{2}$ embeds in any \pisca.)
Therefore, for given K-theory, at most one of the \ca s proved
above to be nonisomorphic can be exact.

\vspace{0.6\baselineskip}

\end{document}